\documentclass[3p]{elsarticle}
\usepackage{graphics}
\usepackage{graphicx}
 \usepackage{epsfig}

\usepackage{amssymb}
\usepackage{stmaryrd}
\usepackage{dsfont}
 \usepackage{amsthm}
  \usepackage{amsmath}
\usepackage{tikz-cd}

 \usepackage{float}
 \usepackage{newfloat}
\usepackage{comment}
\usepackage{ccicons}
\usepackage{epstopdf}

\usepackage{multirow}
\usepackage{color}
\usepackage{hhline}
\usepackage{booktabs}
\usepackage{rotating}
\usepackage{hyperref}
\usepackage{bm}

\definecolor{nverde}{RGB}{0,61,0} %nverde
\definecolor{cr1}{RGB}{200,0,0}
\definecolor{cr2}{RGB}{0,0,200}
\definecolor{cr12}{RGB}{100,0,100}

\newcommand{\halb}{\frac{1}{2}}

\renewcommand{\v}{\mathbf{v}}

\newcommand{\A}{\mathbf{A}}
\newcommand{\G}{\mathbf{G}}
\newcommand{\J}{\mathbf{J}}

\newcommand{\B}{\mathbf{B}}

\newcommand{\x}{\mathbf{x}}
\newcommand{\vc}{|\Omega_c|}

\newcommand{\f}{\mathbf{f}}

\newcommand{\q}{\mathbf{q}}

\newcommand{\devG}{\mathring{G}}
\newcommand{\kp}{n}
\newcommand{\ISigma}{\Sigma^{-1}}

\DeclareFloatingEnvironment[fileext=lop]{Diagram}

\graphicspath{ {./figures/} }

\journal{Journal of Computational Physics}

\allowdisplaybreaks

\begin{document}

\begin{frontmatter}

\title{A structure-preserving semi-implicit four-split scheme for continuum mechanics}

\author[UniTN]{M. Dumbser \corref{cor1}}
\ead{michael.dumbser@unitn.it}
\cortext[cor1]{Corresponding author}

\author[Inria,UniTN]{A. Thomann}
\ead{andrea.thomann@inria.fr}

\author[UniVR]{M. Tavelli}
\ead{maurizio.tavelli@univr.it}

\author[UniSMB,UniFE]{W. Boscheri}
\ead{walter.boscheri@univ-smb.fr}

\address[UniTN]{Laboratory of Applied Mathematics, DICAM, University of Trento, via Mesiano 77, 38123 Trento, Italy}

\address[Inria]{Universit\'e de Strasbourg, CNRS, Inria, IRMA, Strasbourg, F-67000, France}

\address[UniVR]{Department of Engineering for Innovation Medicine, University of Verona, Strada le Grazie 15, 37134, Verona, Italy}

\address[UniSMB]{Laboratoire de Math\'ematiques UMR 5127 CNRS
	Universit\'e Savoie Mont Blanc,
	73376 Le Bourget du Lac, France}

\address[UniFE]{Department of Mathematics and Computer Science, University of
	Ferrara, 44121 Ferrara, Italy}

\begin{abstract}
In this work we introduce a novel structure-preserving vertex-staggered semi-implicit four-split discretization of a unified first order hyperbolic formulation of continuum mechanics that is able to describe at the same time fluid and solid materials in one and the same mathematical model. The governing PDE system goes back to pioneering work of Godunov, Romenski, Peshkov and collaborators. Previous structure-preserving discretizations of this system allowed to respect the curl-free properties of the distortion field and of the specific thermal impulse in the absence of source terms and were also able to properly deal with the low Mach number limit with respect to the adiabatic sound speed. However, the evolution of the thermal impulse and the distortion field were still discretized explicitly, thus requiring a rather severe CFL stability restriction on the time step based on the shear sound speed and on the finite, but potentially large, speed of heat waves. Instead, the new four-split semi-implicit scheme presented in this paper has a CFL time step restriction based only on the magnitude of the velocity field of the continuum. For this purpose, the governing PDE system is split into four subsystems: i) a convective subsystem, which is the only one that is treated explicitly; ii) a heat subsystem, iii) a subsystem containing momentum, distortion field and specific thermal impulse; iv) a pressure subsystem. The last three subsystems ii)-iv) are all discretized implicitly, hence the time step is only limited by a rather mild CFL condition based on the magnitude of the velocity field. The method is consistent with the low Mach number limit of the equations, with the stiff relaxation limits and it maintains an exactly curl-free distortion field and thermal impulse in the case of linear source terms or in their absence. We show several numerical results for classical benchmark problems that allow to assess the performance of the scheme in different asymptotic limits of the governing equations, including the fluid and solid limit.
\end{abstract}


\end{frontmatter}

%%\linenumbers
%\tableofcontents
%% main text
%% The Appendices part is started with the command \appendix;
%% appendix sections are then done as normal sections
%% \appendix

% % % % % % % % % % % % % % % % % % % % % % % % % % % % % %
%         Introduction
% % % % % % % % % % % % % % % % % % % % % % % % % % % % % %

\section{Introduction}\label{sec.intro}

Continuum mechanics is mathematically modeled by time-dependent nonlinear systems of partial differential equations based on the conservation of mass, momentum and total energy, along with the temporal evolution of the material's deformation gradient. These governing equations describe different physical phenomena such as transport, diffusion, shear deformation and friction. It is worth to notice that the time scales associated with each process differ significantly. For example, diffusion occurs over much smaller time scales than advection, and pressure waves propagate much faster than material interfaces or contact discontinuities. Therefore, numerical schemes capable of handling multiple time scales simultaneously are crucial for high fidelity simulations of real-world applications.

Explicit numerical schemes are very effective for dealing with phenomena like shock wave propagation in fluid dynamics or elastic deformations in solids, while implicit methods are preferred to simulate slow diffusion processes, induced by viscous or thermal effects, since they avoid the severe parabolic stability restriction and thus enlarge the maximum admissible time step substantially.  However, multiple time scales can coexist and can arise during one simulation without being predicted in advance, thus needing the design of numerical methods that are able to dynamically adapt to the regime under consideration.

Among the extensive research activity carried out in the recent past for investigating an alternative strategy to treat problems with multiple time scales, a successful idea consists of splitting the processes in fast and slow time scales and treating them with different numerical techniques. More specifically, the terms which might yield a possible stiffness of the governing model undergo an implicit discretization ensuring numerical stability, while the remaining terms are discretized explicitly keeping robustness and shock-capturing properties. Consequently, different type of subsystems which stem from the original model, associated to fast and slow scale phenomena respectively, have to be solved. There are mainly two classes of schemes that permit to treat split subsystems. The first basic class is represented by semi-implicit methods \cite{MunzPark,KleinMach,Casulli1984,Casulli1990,BosFil2016}, which are based on a suitable time linearization of the possibly nonlinear stiff terms, hence leading to the construction of simple linearly implicit schemes \cite{BP2021}.
% and which are usually limited to first order of accuracy in time.

The second class is given by higher order implicit-explicit (IMEX) Runge-Kutta methods \cite{AscRuuSpi,BP2017,BosRus,PR_IMEX,BPR2017} or, more in general, by the so-called partitioned schemes \cite{Hofer}. IMEX schemes achieve the formal order of accuracy under a time step stability constraint independent of the values of the fast scale, and there are stiffly accurate methods available that are proven to satisfy the Asymptotic Preserving (AP) property, meaning that the limit model is consistently reproduced at the discrete level \cite{JINAP1999, JP2001, KLARAP1999}.

Most of the IMEX or semi-implicit schemes have been effectively applied to two subsystems, mainly focusing on fluid dynamics and kinetic equations. In particular, the development of the so-called all speed or all Mach number solvers \cite{MunzPark,Degond2,DumbserCasulli2016,TavelliDumbser2017,Avgerinos2019,Thomann2020} for the compressible Euler and Navier-Stokes equations has been an active field of research as well as the design of AP schemes for relaxation systems like the BGK or Boltzmann model \cite{PR_IMEX,dimarco2012,BD2021_FVBoltz,xiong2015_BGK_DG,bisi2022}. In the first case, the stiffness of the governing equations is generated in the acoustic flux terms by the Mach number regime of the flow, which tends to the incompressible behavior in the low Mach number limit. In the second case, the fast scale phenomena are controlled by a relaxation parameter in the source terms of the system like the Knudsen number.

Dealing with more than two subsystems is still not very common because of the increasing complexity of the numerical method that has necessarily to deal with more linear (or even nonlinear) systems to be solved globally on the entire computational domain. For instance, three subsystems arise when the viscous terms of the compressible Navier-Stokes equations are separated from the advection and pressure contribution, see the three-split schemes in \cite{TavelliDumbser2017,BosTav22}. First, the explicit convective subsystem is solved, then the implicit viscous subsystem is considered, and finally the pressure subsystem is discretized implicitly according to the flux splitting scheme proposed in \cite{ToroVazquez} for the Euler equations of compressible gas dynamics. The resulting scheme is suitable for low Mach number flows, i.e. it is asymptotic preserving, with a time step stability condition that is only based on the material flow speed, thus avoiding any restriction of the timestep by the parabolic terms or the acoustic waves. Another example of three-split methods can be found in the context of magnetized plasma flows, where the fast scales associated to both the acoustic and the Alfv{\'e}n Mach number are separated from the transport part of the system. In \cite{Fambri_3splitMHD,dematte2024} the three-split schemes lead to a nonlinear system associated to the magnetic field that is solved by a fixed point iteration method. Differently, in \cite{BosTho_3splitMHD} a new semi-implicit technique was designed so that the resulting three-split method involves only linear systems to be solved.

The primal objective of this work is to design a novel four-split scheme for the \textcolor{black}{Godunov-Peshkov-Romenski (GPR)} model of continuum mechanics \cite{PeshRom2014} that accounts for both fluid and solid mechanics within one and the same set of equations. It covers a very wide range of phenomena, spanning from ideal and viscous heat conducting fluids to elastic and elasto-plastic solids. 
Indeed, this is possible thanks to the presence of stiff relaxation sources in the evolution equation of the deformation gradient and of the thermal impulse, which retrieve the Navier-Stokes-Fourier model in the stiff relaxation limit \cite{GPRmodel}. For vanishing source terms, the mechanical behavior of ideal elastic solids is reproduced, involving both shear and thermal stresses. 

Our aim is to develop a numerical method whose stability condition on the time step is only dictated by the classical CFL condition based on the velocity field of the medium and not by the other fast waves present in the system, namely acoustic waves, shear waves and heat waves. In order to obtain a method that works uniformly for any material under different flow regimes, we propose to split the original model into four subsystems which involve: i) convection; ii) temperature (heat conduction); iii) deformation gradient and velocity with mechanical and thermal stress; iv) pressure. Our novel scheme is conceived such that the subsystems i) and iv) exactly collapse to the two-split method of \cite{ToroVazquez}, which has been extended to the GPR model of continuum mechanics in \cite{SIGPR} and \cite{ChiocchettiSI}. Subsystem iii) is nonlinear due to the source term in the deformation gradient equation, therefore we rely on a time linearization technique, while the remaining subsystems ii) and iv) are linear and the associated coefficient matrices are symmetric and positive definite, thus allowing the efficient conjugate gradient solver to be used. The conceptualization of the new splitting is inspired by Hamiltonian mechanics \cite{GodRom2003}, where only the Euler-Lagrange equations derive from Hamilton's principle of stationary action, while the remaining geometry equations are consequences of the definitions. The coupling of the subsystems follows this guideline introduced in \cite{MaxwellGLM,CompatibleDG1}, which will be more detailed in the sequel. The new scheme is asymptotic preserving in the low acoustic Mach number limit, meaning that a consistent discretization of the incompressible Navier-Stokes equation is retrieved, and it is also asymptotic preserving with respect to the Navier-Stokes-Fourier stress tensor for stiff relaxation sources.

Furthermore, the GPR model is also endowed with two stationary differential constraints on the curl of the inverse of the deformation gradient and the thermal impulse. From the numerical viewpoint, this implies the construction of structure-preserving spatial discretizations which can exactly mimic the classical vector calculus identities $\nabla \cdot (\nabla \times \mathbf{a})=0$  and $\nabla \times \nabla \psi = \mathbf{0}$, with $\mathbf{a} \in \mathds{R}^3$ and $\psi \in \mathds{R}$. The literature of div-curl and curl-grad preserving operators is very vast, so we limit ourselves to recall some of the main contributions in the context of Cartesian meshes, as adopted in this work. The need of exactly discrete div-curl operators is very well known for the numerical solution of the Maxwell and MHD equations, see for instance the schemes proposed in \cite{Yee66,DeVore,BalsaraSpicer1999,Balsara2004,GardinerStone,balsarahlle2d,ADERdivB}. Mimetic finite difference operators have been developed in \cite{HymanShashkov1997,JeltschTorrilhon2006,Torrilhon2004,Margolin2000,Lipnikov2014,Carney2013}, while compatible finite element schemes can be found, for example, in \cite{Nedelec1,Nedelec2,Cantarella,Hiptmair,Monk,Arnold,Alonso2015,CAMPOSPINTO2016,Zampa1,Zampa2}. In \cite{SPDG2023,Ern2023,perrier2024} discontinuous Galerkin div-curl and curl-grad preserving methods are forwarded, while methods based on a discrete de Rahm complex and Compatible Discrete Operators (CDO) have been presented in \cite{bonelle2015,DiPietro2023}. A common strategy of the aforementioned works is the use of staggered meshes, where the variables of the governing equations are not discretized at the same location on the computational grid. Typically, edge-based staggering is adopted \cite{Yee66,Margolin2000}, meaning that there is an interplay between the main and the dual grid, so that the algebraic operator is satisfied only on one grid. This strategy has been adopted in \cite{SIGPR} for the design of curl-free operators for the GPR model. Recently in \cite{Barsukow2024,Sidilkover2025}, a vertex-staggered approach has proven to be very effective for the construction of compatible operators even on unstructured meshes. In \cite{CompatibleDG1}, a quite general and simple framework is proposed, combining Discontinuous Galerkin and Finite Element discretizations in such a way as to achieve the classical vector calculus identities as special case of a discrete version of the Schwarz theorem.

The second objective of this work is thus to exactly satisfy the curl involutions which are present in the GPR model for vanishing source terms. Differently form \cite{SIGPR}, the resulting structure-preserving scheme is discretized at the aid of a vertex-staggered mesh only, hence remarkably simplifying the discrete differential operators. More specifically, no edge-based staggering is used here for the discretization of the velocity field, thus the same compatible div-curl and curl-grad operators only act on the main (cell centered) and the dual (vertex centered) grid. Eventually, our structure-preserving div-curl operator allows to recover exactly at the fully discrete level also the quadratic convergence in the low Mach number limit.

The rest of the paper is organized as follows. In Section \ref{sec.pde} we briefly recall the GPR model of continuum mechanics and we present the splitting strategy with the associated subsystems. In Section \ref{sec.method} we present the new four-split structure-preserving finite volume scheme. Computational results for a wide range of mechanical regimes are shown in Section \ref{sec.results}.
The paper closes with Section \ref{sec.conclusions}, in which we give some concluding remarks and an outlook to future work.

% % % % % % % % % % % % % % % % % % % % % % % % % % % % % %
%         splitting
% % % % % % % % % % % % % % % % % % % % % % % % % % % % % %

\section{Governing PDE system and splitting}\label{sec.pde}

The unified first order hyperbolic model of continuum mechanics of \textcolor{black}{Godunov, Peshkov and Romenski, \cite{Rom1998,PeshRom2014,GPRmodel,HTCGPR}, hereafter also simply denoted by GPR model}, reads as follows: 
\begin{subequations}\label{eqn.GPR}
	\begin{align}
		& \frac{\partial \rho}{\partial t}+\frac{\partial (\rho v_k)}{\partial
			x_k}=0,\label{eqn.conti}\\[2mm]
		&\frac{\partial \rho v_i}{\partial t}+\frac{\partial \left(\rho v_i v_k + p \, \delta_{ik} +
			\sigma_{ik} + \omega_{ik} \right)}{\partial x_k}=0, \label{eqn.momentum}\\[2mm]
		&\frac{\partial A_{i k}}{\partial t}+\frac{\partial (A_{im} v_m)}{\partial x_k} +
		v_m \left(\frac{\partial A_{ik}}{\partial x_m}-\frac{\partial A_{im}}{\partial x_k}\right)
		=-\dfrac{ \alpha_{ik} }{{\theta}_1(\tau_1)},\label{eqn.deformation}\\[2mm]
		&\frac{\partial J_k}{\partial t}+\frac{\partial \left( J_m v_m + T\delta_{ki} \right)}{\partial x_k} +
		v_m \left(\frac{\partial J_{k}}{\partial x_m}-\frac{\partial J_{m}}{\partial x_k}\right)  =
		-\dfrac{\beta_k}{{\theta}_2(\tau_2)}, \label{eqn.heatflux}\\[2mm]
		& \frac{\partial \mathcal{E}}{\partial t}+\frac{\partial \left( v_k \mathcal{E} + v_i (p \, \delta_{ik}
			+ \sigma_{ik} + \omega_{ik} ) + q_k \right)}{\partial x_k}=0. \label{eqn.energy}
	\end{align}
\end{subequations}
\textcolor{black}{Throughout this paper we make use of the Einstein summation convection over two repeated indices.}
The system is thermodynamically compatible and satisfies the following entropy inequality:
\begin{equation}
	\frac{\partial \rho S}{\partial t}+\frac{\partial \left( \rho S v_k  + \beta_k \right)}{\partial x_k} = \frac{\alpha_{ik} \alpha_{ik}}{T \, {\theta}_1(\tau_1)} + \frac{\beta_k \beta_k}{T \, {\theta}_2(\tau_2)} \geq 0. \label{eqn.entropy}
\end{equation}
Throughout this paper we assume that the temperature is positive, i.e. $T > 0$.
The state vector is given by $\mathbf{q}= (\rho, \rho v_i, \mathcal{E}, A_{ik}, J_k)^T$, with the mass density $\rho$, the velocity field $v_i$, the distorsion field $A_{ik}$, the specific thermal impulse $J_k$ and the total energy density $\mathcal{E} = \rho E = \mathcal{E}_1  +\mathcal{E}_2 , + \mathcal{E}_3 + \mathcal{E}_4 $
containing four different contributions, namely internal energy, kinetic energy, energy due to the elastic deformation of the medium and a contribution due to the presence of the specific thermal impulse,
\begin{equation}
	\mathcal{E}_1 = \frac{\rho^\gamma}{\gamma-1} e^{S/c_v}, \quad
	\mathcal{E}_2 = \halb \rho v_i v_i,
	\quad
	\mathcal{E}_3 = \frac{1}{4} \rho c_s^2 \mathring{G}_{ij} \mathring{G}_{ij},
	\quad
	\mathcal{E}_4 = \halb c_h^2 \rho J_i J_i.
\end{equation}
Here, the metric tensor is denoted by $\mathbf{G}$ with $
{G}_{ik} = A_{ji} A_{jk} $ and its trace-free (deviatoric) part $\mathring{\mathbf{G}}$ is given by
$\mathring{G}_{ik} = {G}_{ik} - \frac{1}{3} \, G_{mm} \delta_{ik}
$.
The thermodynamic dual variables or main field read  $\mathbf{p} = \partial_{\q} \mathcal{E}=\{p_i\} = \left( r,
v_i, T, \alpha_{ik}, \beta_k \right)^T$ with
\begin{equation}
	r = \partial_{\rho} \mathcal{E}, %=\frac{\partial \mathcal{E}}{\partial \rho},
	\qquad
	v_i = \partial_{\rho v_i} \mathcal{E}, %=\frac{\partial \mathcal{E}}{\partial (\rho v_i)},
	\qquad
	T = \partial_{\rho S} \mathcal{E}, %=\frac{\partial \mathcal{E}}{\partial (\rho S)},
	\qquad
	\alpha_{ik} = \partial_{A_{ik}} \mathcal{E}, %= \frac{\partial \mathcal{E}}{\partial A_{ik}},
	\qquad
	\beta_{k} = \partial_{J_{k}} \mathcal{E}. %= \frac{\partial \mathcal{E}}{\partial J_{k}}.
\end{equation}
The pressure is given by the relation
$
p =  \rho^2 \partial_\rho E,
$
and the two stress tensors read
\begin{equation}
	\sigma_{ik} = A_{ji} \partial_{A_{jk}} \mathcal{E} =
	A_{ji} \alpha_{jk} = \rho c_s^2 G_{ij} \mathring{G}_{jk},
	\qquad
	\omega_{ik} = J_i \partial_{J_{k}} \mathcal{E} = J_i \beta_{k} = \rho c_h^2 J_i J_k.
\end{equation}
The heat flux is given by
\begin{equation}
	q_k = \partial_{\rho S} \mathcal{E} \, \partial_{J_k} \mathcal{E} = T \beta_k = \rho c_h^2 T J_k,
\end{equation}
while $\theta_1(\tau_1)>0$ and $\theta_2(\tau_2)>0$ are two  functions of $\mathbf{q}$ that depend on the two relaxation times $\tau_1>0$ and $\tau_2>0$ and are defined as
\begin{equation}\label{eqn:theta1}
	\theta_1 = \frac{1}{3}  \rho \tau_1 \, c_s^2 \, \left| \mathbf{A} \right|^{-\frac{5}{3}},
	\qquad
	\theta_2 = \rho c_h^2 \tau_2.
\end{equation}
The determinants of $\A$ and $\G$ are related to the density of the medium and the reference density $\rho_0$ by the relations
\begin{equation}
	|\A| = \frac{\rho}{\rho_0}, \qquad |\G| = \left( \frac{\rho}{\rho_0} \right)^2. 
\end{equation}
A formal asymptotic analysis of the model \cite{GPRmodel}  showed that for $\tau_1 \to 0$  the stress tensor $\sigma_{ik}$ tends to 
\begin{equation}
	\sigma_{ik} = -\frac{1}{6} \rho_0 c_s^2 \tau_1 \left( \partial_k v_i + \partial_i v_k
	- \frac{2}{3} \left( \partial_m v_m\right) \delta_{ik} \right). 
	\label{eqn.asymptoticlimit.stress}
\end{equation}
Instead, for $\tau_2 \to 0$ and the $\theta_2$ used in this paper, which is different from \cite{GPRmodel} and which leads to a \textit{linear} relaxation source term in the equation for $\J$, the heat flux $q_k$ tends to 
\begin{equation}
	 q_k = -\rho T c_h^2 \tau_2 \partial_k T, 
	\label{eqn.asymptoticlimit.heatflux}
\end{equation}
hence retrieving the compressible Navier-Stokes-Fourier equations,  with shear viscosity $ \mu = \frac{1}{6} \rho_0 c_s^2 \tau_1 $ and thermal conductivity \mbox{$ \lambda = \rho T c_h^2 \tau_2 $}.

In the numerical scheme, we will \textcolor{black}{also make use of the governing equation for the metric tensor as auxiliary quantity, since it will allow us to obtain a linearly implicit scheme in the strain relaxation source term that is asymptotically consistent with the Navier-Stokes equations in the stiff relaxation limit when $\tau_1 \to 0$.}  It can be derived from \eqref{eqn.deformation} and reads
\begin{equation}
\frac{\partial G_{i k}}{\partial t}
+ v_m \frac{\partial G_{ik} }{\partial x_m}
+ G_{im} \frac{\partial v_m }{\partial x_k}
+ G_{mk} \frac{\partial v_m }{\partial x_i}
=-\dfrac{ 2 \sigma_{ik} }{\theta_1(\tau_1)},
\label{eqn.pde.metric}
\end{equation}
while it is sometimes more convenient to use the equation for the specific thermal impulse rewritten as
\begin{equation}
	\frac{\partial J_{k}}{\partial t}
	+ v_m \frac{\partial J_{k} }{\partial x_m}
	+ J_{m} \frac{\partial v_m }{\partial x_k}
	+ \frac{\partial T }{\partial x_k}
	=-\dfrac{ \beta_{k} }{{\theta}(\tau_2)}.
	\label{eqn.pde.J}
\end{equation}
In the following we introduce three characteristic Mach numbers which induce the different scales in the system.
The first is the classical acoustic Mach number
\begin{equation}
	M_a = \frac{\| \v \|}{c_0}, \qquad c_0 = \sqrt{\frac{\gamma p}{\rho}},
\end{equation}
with $\v$ the flow speed and $c_0$ the adiabatic sound speed; second, the so-called shear Mach number
\begin{equation}
	M_s = \frac{\| \v \|}{c_s},
\end{equation}
with respect to the shear sound speed  and finally the so-called heat Mach number
\begin{equation}
	M_h = \frac{\| \v \|}{c_T}, \qquad c_T = c_h \sqrt{\frac{T}{c_v}}
\end{equation}
with respect to the heat propagation speed $c_T$.
In the following we will introduce the splitting of system \eqref{eqn.GPR} in four subsystems.

\subsection{Convective subsystem}
The first subsystem is the convective subsystem, which contains the mass flux and the transport of momentum, distortion field, specific thermal impulse, as well as the transport of kinetic energy, deformation energy and the contribution of the thermal impulse to the total energy.
It reads as follows
\begin{subequations}\label{eqn.GPR.s1}
	\begin{align}
		& \frac{\partial \rho}{\partial t}+\frac{\partial (\rho v_k)}{\partial
			x_k}=0,\label{eqn.conti.1}\\[2mm]
		&\frac{\partial \rho v_i}{\partial t}+\frac{\partial \left(\rho v_i v_k  \right)}{\partial x_k}=0, \label{eqn.momentum.1}\\[2mm]
		&\frac{\partial A_{i k}}{\partial t}  +
		v_m  \frac{\partial A_{ik}}{\partial x_m}
		=0,\label{eqn.deformation.1}\\[2mm]
		&\frac{\partial J_k}{\partial t}
		 + v_m \frac{\partial J_{k}}{\partial x_m} = 0, \label{eqn.heatflux.1}\\[2mm]
		& \frac{\partial \mathcal{E}}{\partial t}+\frac{\partial \left( v_m \mathcal{E}_{2,3,4}  \right)}{\partial x_m}=0, \label{eqn.energy.1}
	\end{align}
\end{subequations}
with $\mathcal{E}_{2,3,4} = \mathcal{E}_{2}+\mathcal{E}_{3}+\mathcal{E}_{4}$.
The convective subsystem in direction $x_1$ has the eigenvalues $\lambda_1=0$ and $\lambda_i = v_1$, $i \in \left\{2, 3, \cdots, 17 \right\}$ thus only exhibits material waves. The subsystem is weakly hyperbolic since a full set of eigenvectors cannot be established.

\subsection{Temperature subsystem}
The temperature subsystem contains the coupling of the specific thermal impulse with the temperature and the heat flux of the total energy equation. Since in addition only the density is a contributing variable, all further variables have been omitted. Thus it reads
\begin{subequations}\label{eqn.GPR.s2}
	\begin{align}
		&\frac{\partial \rho}{\partial t} = 0, \\[2mm]
		&\frac{\partial J_k}{\partial t}
		+  \frac{\partial T}{\partial x_k} = -\frac{1}{\theta_2} \beta_k, \label{eqn.heatflux.2} \\[2mm]
		& \frac{\partial \mathcal{E}}{\partial t}+\frac{\partial q_k }{\partial x_k}=0. \label{eqn.energy.1}
	\end{align}
\end{subequations}
Since with our choice of energies we have $\mathcal{E}_1 = \rho c_v T$, $\mathcal{E}_4 = \halb \rho c_h^2 J_k J_k$ and $q_k = \rho T c_h^2 J_k$ the eigenvalues of the temperature subsystem in the $x_1$ direction are
$\lambda_{1,5} = \mp c_h \sqrt{T/c_v} = \mp c_T$ and $\lambda_{2,3,4}=0$.
Hence, this subsystem is only connected to the heat Mach number $M_h$.
For physically reasonable states $\rho >0$, $T>0$ a full set of linearly independent eigenvectors exists, hence the temperature subsystem is hyperbolic.
In terms of the primitive variables $J_k$ and $T$ and setting $\theta_2 = \tau_2 \rho c_h^2$ it takes the simple form
\begin{subequations}\label{eqn.GPR.s2b}
	\begin{align}
		&\frac{\partial \rho}{\partial t} = 0, \\[2mm]
		&\frac{\partial J_k}{\partial t}
		+  \frac{\partial T}{\partial x_k} = -\frac{1}{\tau_2} J_k, \label{eqn.heatflux.2b} \\[2mm]
		& \frac{c_v}{T c_h^2} \frac{\partial T}{\partial t}+\frac{\partial J_k }{\partial x_k}=0. \label{eqn.energy.1b}
	\end{align}
\end{subequations}

\subsection{G-J-v subsystem}
The next subsystem concerns the coupling between momentum equation via the two stress tensors with the distortion matrix and specific thermal impulse. Here, we have replaced \eqref{eqn.deformation} by \eqref{eqn.pde.metric}.
Due to the appearance of the stress tensors in the total energy equation, it is added to this subsystem while the other variables have been omitted:
\begin{subequations}\label{eqn.GPR.s3}
	\begin{align}
	&\frac{\partial \rho v_i}{\partial t}
	+\frac{\partial \left(\sigma_{ik} + \omega_{ik} \right)}{\partial x_k}=0, \label{eqn.momentum.3}\\[2mm]
	& \frac{\partial J_{k}}{\partial t}
	+ J_{m} \frac{\partial v_m }{\partial x_k}
	=-\dfrac{ J_{k} }{\tau_2},
	\label{eqn.heatflux.3} \\[2mm]
	&	\frac{\partial G_{i k}}{\partial t}
	+ G_{im} \frac{\partial v_m }{\partial x_k}
	+ G_{mk} \frac{\partial v_m }{\partial x_i}
	=-\dfrac{ 2 \sigma_{ik} }{\theta_1(\tau_1)}, \label{eqn.metric.3}  \\[2mm]
	& \frac{\partial \mathcal{E}}{\partial t}+\frac{\partial v_i (\sigma_{ik}+\omega_{ik}) }{\partial x_k}=0. \label{eqn.energy.3}
	\end{align}
\end{subequations}
Analytic expressions of the eigenvalues of the G-J-v subsystem are not available, however with the quasi-2D setting $G = \text{diag}(G_{kk})$, $J_3 = 0$, we obtain the eigenvalues
and a full set of eigenvectors.
Due to the coupling of $J_k$ and $G_{ik}$ this subsystem is connected both to the shear and the heat Mach numbers $M_s$ and $M_h$ respectively.

\subsection{Pressure subsystem}
The final subsystem is the classical pressure subsystem already contained in the Toro-V\'azquez splitting \cite{ToroVazquez} and which includes the coupling of the momentum with the pressure and specific enthalpy. It is given by
\begin{subequations}\label{eqn.GPR.s4}
	\begin{align}
		&\frac{\partial \rho}{\partial t} = 0, \\[2mm]
		&\frac{\partial \rho v_i}{\partial t}+\frac{\partial  p    }{\partial x_i}=0, \label{eqn.momentum.4}\\[2mm]
		& \frac{\partial \mathcal{E} }{\partial t}+\frac{\partial \left( v_k ( \mathcal{E}_1 +  p ) \right)}{\partial x_k}=0. \label{eqn.energy.4}
	\end{align}
\end{subequations}
Since it is identical with the one introduced in \cite{MunzPark,ToroVazquez} the eigenvalues in the $x_1$ direction are
$\lambda_{1,5}=\halb \left( v_1 \mp a \right) $  and $\lambda_{2,3,4}=0$, with
$a^2 = v_1^2 + 4 c_0^2$ and $c_0^2 = \gamma p / \rho$. This subsystem is connected only to the acoustic Mach number $M_a$.
One can show that a full set of linearly independent eigenvectors exists for physically reasonable states $\rho >0$, $p>0$, hence the pressure subsystem is hyperbolic.

This concludes the description of the four subsystems which form the basis of the new 4-split scheme introduced in this paper.

% % % % % % % % % % % % % % % % % % % % % % % % % % % % % %
%        Numerical method.
% % % % % % % % % % % % % % % % % % % % % % % % % % % % % %

\section{Numerical method}\label{sec.method}

\begin{figure}[!h]
	\centering
	\includegraphics[width=0.45\textwidth]{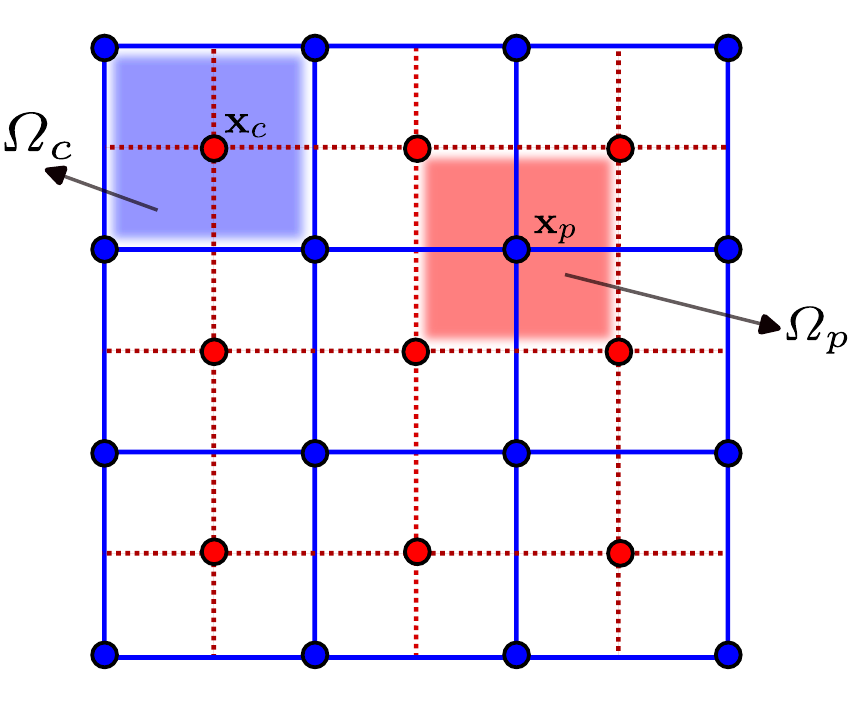}
	\caption{Vertex-based staggered mesh, with the primary cells $\Omega_c$ and the dual cells $\Omega_p$.}
	\label{fig:Grid}
\end{figure}
In the following, we assume summation over repeated indices according to the usual Einstein summation convention.
Let us consider a two-dimensional domain $\Omega$ which is partitioned into primal uniform Cartesian grid cells $\Omega_c$ centred around the nodes $\x_c$. Furthermore, we introduce dual uniform Cartesian grid cells $\Omega_p$ which are centred around the vertices of the primal grid $\x_p$, see also Figure \ref{fig:Grid}.
We hence distinguish between discrete scalar fields $\phi^p=\phi(\x^p)$ and vector fields $\A^p=\A(\x^p)$ defined on the vertices of the primal grid and scalar and vector fields defined on the centers of the primal grid $\phi^c=\phi(\x^c)$ and $\A^c=\A(\x^c)$, respectively.

%For discrete scalar fields $\phi^p=\phi(\x^p)$ and vector fields $\A^p=\A(\x^p)$ defined on the vertices of the primal grid and for scalar and vector fields defined on the centers of the primal grid (the vertices of the dual grid)  $\phi^c=\phi(\x^c)$ and $\A^c=\A(\x^c)$ the following discrete version of the gradient ($\v = \nabla \phi$), divergence ($\psi = \nabla \cdot \A$) and curl ($\B = \nabla \times \A$) operators and their duals are naturally defined as follows, where in the following we assume summation over repeated indices, according to the usual Einstein summation convention:

Based on this notation, the discrete versions of gradient ($\v = \nabla \phi$), divergence ($\psi = \nabla \cdot \A$) and curl ($\B = \nabla \times \A$) operators and their duals are naturally defined as follows
\begin{alignat}{2}
	& v^c_k=\partial^{cp}_k \phi^p   = \frac{1}{|\Omega_c|} \sum_{p \in \Omega_c}  n^{cp}_k  \phi^p,  \qquad
	&&v^p_k=\partial^{pc}_k \phi^c   = \frac{1}{|\Omega_p|} \sum_{c \in \Omega_p}  n^{pc}_k  \phi^c,
	\\
	& \psi^c = \partial^{cp}_k A^p_k   = \frac{1}{|\Omega_c|} \sum_{k} \sum_{p \in \Omega_c}  n^{cp}_k  A^p_k,  \qquad
	&&\psi^p = \partial^{pc}_k A^c_k   = \frac{1}{|\Omega_p|} \sum_{k} \sum_{c \in \Omega_p}  n^{pc}_k  A^c_k,
	\\
	&B^c_i= \epsilon_{ijk} \partial^{cp}_j A^p_k   = \frac{\epsilon_{ijk}}{|\Omega_c|} \sum_{k} \sum_{p \in \Omega_c}  n^{cp}_j   A^p_k,  \qquad
	&&B^p_i= \epsilon_{ijk} \partial^{pc}_j A^c_k   = \frac{\epsilon_{ijk}}{|\Omega_p|} \sum_{k} \sum_{c \in \Omega_p}  n^{pc}_j   A^c_k.
\end{alignat}
On the regular Cartesian meshes considered in this paper we have 
$ |\Omega_c| = |\Omega_p| = \Delta x \Delta y$ and the corner normal vectors $n_k^{cp}$ are simply 
$ n_k^{pc} = \halb ( \pm \Delta y, \pm \Delta x)^T$, with $n_k^{cp} = -n_k^{pc}$. For the definition of the corner normals on general meshes see \cite{Maire2007,Despres2009,Maire2011,Maire2020,HTCLagrange,HTCLagrangeGPR}. 

In \cite{SIGPR,MaxwellGLM} it was already shown that the discrete operators introduced above satisfy the discrete vector calculus identities 
\begin{equation}
	\epsilon_{ijk} \partial_j^{pc} \partial_k^{cq} \phi^q = 0, \qquad 
	\epsilon_{ijk} \partial_i^{pc} \partial_j^{cq} A_k^q = 0, 			
\end{equation}
which are discrete analogues of $\nabla \times \nabla \phi = 0$ and $\nabla \cdot \nabla \times \A= 0$. Since it is very easy to check their validity we do not repeat the calculations here.

In our scheme the density, the pressure, the temperature, the distortion field and the specific thermal impulse are defined in the cells $c$, while the momentum is defined in the vertices $p$. When needed, simple arithmetic averaging between the grid location is carried out to transfer data from one grid to the other one.

In the following we present the fully discrete semi-implicit four-split scheme using the discrete operators given above.
Therein, the steps of the numerical scheme are given by the order of the subsystems given above, i.e. we first solve the convective, then the temperature subsystem, followed by the G-J-v and pressure sub-systems.

\subsection{Discrete convective subsystem}

The convective subsystem is evolved via an explicit cell-centered finite volume scheme that reads
\begin{equation}
	\q^{*,c} = \q^{n,c} - \frac{\Delta t}{\vc} \sum_{p \in \mathcal{P}_c} n_k^{pc} \,   
	\f^{p}_k,
	\label{eqn.fv}
\end{equation}
 with a nodal numerical flux tensor $\f_k^{p}$ of the Rusanov-type
\begin{equation}
	\f_k^{p} = \bar{\f}_k^{p} 
	- \halb h^p s^p_{\max} \partial_k^{pc} q^{n,c},
\label{eqn.rusanov}
\end{equation}
where $h^p = \frac{4|\Omega_p|}{|\partial \Omega_p|} = \frac{2\Delta x \Delta y}{\Delta x + \Delta y}$ is a characteristic length scale at the vertex, 
$$\bar{\f}_k^{p} = \frac{1}{|\mathcal{C}_p|} \sum \limits_{c \in \mathcal{C}_p}  \f_k^c $$
is the average flux around a node and  
$s^p_{\max} = \max \limits_c(|\v^c|)$ is the maximum convective speed at the vertex.
Furthermore, $\mathcal{P}_c$ denotes the set of points $p$ of cell $c$ and $\mathcal{C}_p$ denotes the set of cells $c$ around a point $p$. 
Since it is the only subsystem which is solved explicitly, we obtain immediately a CFL condition that is independent of the Mach number, 
\begin{equation}
    \Delta t = \textnormal{CFL} \min \limits_{p} \left( \frac{h^p}{s_{\max}^{p}} 
    \right), 
\end{equation}
for a CFL coefficient $\textnormal{CFL} \leq \halb$.

\subsection{Discrete temperature subsystem}
The temperature subsystem \eqref{eqn.GPR.s2}, taking into account the dependency of the total energy $\mathcal{E}$ on the temperature and thermal impulse yields a non-conservative auxiliary subsystem \eqref{eqn.GPR.s2b} coupling $T$ and $J_k$. 
Following the above given convention, the density and temperature are discretized in the cell center, while an auxiliary thermal impulse defined on the nodes $J_k^{**,p}$ is introduced. 
Following the definition of the respective discrete compatible operators, we obtain the following discretization
\begin{subequations}
	\begin{align}
    J_k^{**,p} &= J_k^{*,p} - \Delta t \, \partial^{pc}_k T^{**,c} - \frac{\Delta t}{\tau_2} J_k^{**,p}, \label{eqn.heatflux.disc.2} \\
    \frac{c_v}{T^{n,c} c_h^2} T^{**,c} &= \frac{c_v}{T^{n,c} c_h^2}  T^{*,c} - \Delta t \,
    \partial^{cp}_k \, J_k^{**,p}. \label{eqn.energy.disc.2}
    \end{align}
\end{subequations}
Therein, quantities denoted by $\ast$ are the values after the convective step and $\ast\ast$ denotes the newly updated variables.
Note that in this subsystem the effect of the relaxation with respect to heat conduction on the temperature is taken into account. 
Substitution of \eqref{eqn.heatflux.disc.2} into \eqref{eqn.energy.disc.2} yields the following discrete scalar wave equation for the temperature $T^{\ast\ast,c}$
\begin{equation}
	  \left( 1+\frac{\Delta t}{\tau_2} \right) \frac{c_v}{T^{n,c} c_h^2} \, T^{**,c}  - \Delta t^2
	\partial^{cp}_k \left( \partial^{pa}_k T^{**,a} \right)
	= \left( 1+\frac{\Delta t}{\tau_2} \right) \frac{c_v}{T^{n,c} c_h^2} \, T^{*,c} - \Delta t \,
		\partial^{cp}_k \left( J_k^{*,p} \right). \label{eqn.wave.T}
\end{equation}
The system is symmetric and positive definite and can therefore be solved easily via a classical matrix-free conjugate gradient method.
Once $T^{**,c}$ is obtained, the thermal impulse and energy update are 
\begin{subequations}
	\begin{align}\label{eqn.heat.disc.final}
		J^{**,c}_k &= \left(\frac{\tau_2}{\tau_2 + \Delta t }\right)(J^{*,c}_k - \Delta t \, \partial_k^{cp} T^{**,p}), \\
		\mathcal{E}^{**,c} &= \mathcal{E}^{*,c} - \Delta t \, \partial_k^{cp} q_k^{**,p}.\label{eqn.heat.disc.final.E}
	\end{align}
\end{subequations}
Therein, $q_k^{**,p} =  \rho^{*,p} \, c_h^2 \, T^{**,p} \, J^{**,p}$, where $T^{**,p}$ and $\rho^{*,p}$ are obtained by arithmetic averaging from $T^{**,c}$ and $\rho^{*,c}$, respectively, while the thermal impulse on the nodes is obtained by solving \eqref{eqn.heatflux.disc.2}. 

\subsection{Discrete $\devG$-J-v subsystem}
In order to obtain the $\devG-J-v$ system which is actually used in our scheme, we start from a discrete $G-J-v$ subsystem:
\begin{subequations}
	\begin{align}
		\rho^{*,p} v_i^{**,p}  &= \rho^{*,p} v_i^{*,p} - \Delta t \, \partial^{pc}_k
		\left(  \sigma_{ik}^{**,c} + \omega_{ik}^{***,c} \right), \label{eqn.mom.disc.3} \\
		G_{ik}^{**,c} &= G_{ik}^{*,c}
		 - \Delta t \, G_{im}^{n,c} \, \partial_k^{cp} v_m^{**,p}
		 - \Delta t \, G_{mk}^{n,c} \, \partial_i^{cp} v_m^{**,p}
		 -  \dfrac{ 2 \Delta t  }{ \theta_1^{n,c} } \, \sigma_{ik}^{**,c},
		 \label{eqn.metric.disc.3} \\
		J_{k}^{***,c} &= J_{k}^{**,c}
		- \Delta t \, J_{m}^{**,c} \, \partial_k^{cp} v_m^{**,p}
		-  \dfrac{  \Delta t  }{ \theta_2} \, J_{k}^{***,c},
		\label{eqn.J.disc.3}
	\end{align}
\end{subequations}
where $v_i^{**,p}, G_{ik}^{**,c}$ and $J_k^{***,c}$ are the unknown quantities, as $J_k$ has contributions from the two previous subsystems while the metric tensor and the velocity only have contributions from the first subsystem, and $\omega_{ik}^{***,c}$, $\sigma_{ik}^{**,c}$ are some proper semi-implicit discretizations of the stress tensors that will be detailed later.
Applying the deviator (trace-free) operator on \eqref{eqn.metric.disc.3} we obtain 
\begin{eqnarray}
	\devG_{ik}^{**,c}
	%&=&\devG_{ik}^{n,c}
	%                  - \Delta t\left(G_{im}^{\kp,c}\partial_k^{cp} v_m^{**,p}+G_{mk}^{\kp,c}\partial_i^{cp} v_m^{**,p}\right)
	%                  +\frac{\Delta t}{3}tr\left(G_{im}^{\kp,c}\partial_k^{cp} v_m^{**,p}+G_{mk}^{\kp,c}\partial_i^{cp} %v_m^{**,p}\right)\delta_{ik}
	%                  -\frac{2\Delta t}{\rho \theta_1}\mathring{\sigma}_{ik}^{**,c} \nonumber \\
  	%&=& \devG_{ik}^{n,c}
  	%- \Delta t\left(G_{im}^{\kp,c}\partial_k^{cp} v_m^{**,p}+G_{mk}^{\kp,c}\partial_i^{cp} v_m^{**,p}\right)
  	%+ \frac{\Delta t}{3} \left(G_{nm}^{\kp,c}\partial_n^{cp} v_m^{**,p}+G_{mn}^{\kp,c}\partial_n^{cp} v_m^{**,p}\right)\delta_{ik}
  	%-\frac{2\Delta t}{\rho \theta_1}\mathring{\sigma}_{ik}^{**,c} \nonumber \\
  	%&=& \devG_{ik}^{n,c}
%- \Delta t\left(G_{im}^{\kp,c}\partial_k^{cp} v_m^{**,p}+G_{mk}^{\kp,c}\partial_i^{cp} v_m^{**,p}\right)
%+ \frac{\Delta t}{3} G_{nm}^{\kp,c}\left(\partial_n^{cp} v_m^{**,p}+\partial_n^{cp} v_m^{**,p}\right)\delta_{ik}
%-\frac{2\Delta t}{\rho \theta_1}\mathring{\sigma}_{ik}^{**,c} \nonumber \\
	&=& \devG_{ik}^{n,c}
- \Delta t\left(G_{im}^{\kp,c}\partial_k^{cp} v_m^{**,p}+G_{mk}^{\kp,c}\partial_i^{cp} v_m^{**,p}-\frac{2}{3}\delta_{ik}G_{nm}^{\kp,c}\partial_n^{cp}v_m^{**,p}\right)
-\frac{2\Delta t}{\theta_1^{n,c}}\mathring{\sigma}_{ik}^{**,c}
\label{eqn.devG.1}
\end{eqnarray}
where we used the symmetry of $\G$ to simplify the trace contribution. We now introduce the following semi-implicit approximations for the two stress tensors
\begin{equation}
	\sigma_{ik}^{**,c} = \rho c_s^2 G_{il}^{\kp,c}\devG_{lk}^{**,c}
	\label{eqn.si.sigma}
\end{equation} 
and
\begin{equation}
	\omega_{ik}^{***,c} = \rho c_h^2 J_i^{**,c} J_k^{***,c}.
	\label{eqn.si.omega}
\end{equation} 
A semi implicit discretization of $\mathring{\sigma}$, required in \eqref{eqn.devG.1}, reads
\begin{equation}
	\mathring{\sigma}_{ik}^{**,c} = \rho^{n,c} \, c_s^2 \left(G_{il}^{\kp,c}\devG_{lk}^{**,c}-\frac{1}{3}G_{mn}^{\kp,c}\devG_{nm}^{\kp,c}\delta_{ik}\right).
	\label{eqn.si.devsigma}
\end{equation} 
Substituting the quantity \eqref{eqn.si.devsigma} into \eqref{eqn.devG.1} we obtain 
\begin{eqnarray}
	\devG_{ik}^{**,c}&=& \devG_{ik}^{n,c}
	- \Delta t\left(G_{im}^{\kp,c}\partial_k^{cp} v_m^{**,p}+G_{mk}^{\kp,c}\partial_i^{cp} v_m^{**,p}-\frac{2}{3}\delta_{ik}G_{nm}^{\kp,c}\partial_n^{cp}v_m^{**,p}\right) \nonumber \\
	&&-\frac{2\Delta t \rho^{n,c} c_s^2}{\theta_1^{n,c}} G_{il}^{\kp,c}\devG_{lk}^{**,c}
	+\frac{2\Delta t \rho^{n,c} c_s^2}{3 \, \theta_1^{n,c}} G_{mn}^{\kp,c}\devG_{nm}^{\kp,c}\delta_{ik}
	\label{eqn.devG.2}
\end{eqnarray}
and then, grouping $\devG^{**,c}$, one has 
\begin{eqnarray}
	\left(\delta_{il}+\frac{2\Delta t \rho^{n,c} c_s^2}{\theta_1^{n,c}}G_{il}^{\kp,c}\right)\devG_{lk}^{**,c} &=& \devG_{ik}^{n,c}
	- \Delta t\left(G_{im}^{\kp,c}\partial_k^{cp} v_m^{**,p}+G_{mk}^{\kp,c}\partial_i^{cp} v_m^{**,p}-\frac{2}{3}\delta_{ik}G_{nm}^{\kp,c}\partial_n^{cp}v_m^{**,p}\right)  \nonumber \\
	&&+\frac{2\Delta t \rho^{n,c} c_s^2}{3 \theta_1^{n,c}} G_{mn}^{\kp,c}\devG_{nm}^{\kp,c}\delta_{ik}.
	\label{eqn.devG.3}
\end{eqnarray}
It is now convenient to define 
\begin{equation}
	\Sigma_{il}^{\kp,c}:=\delta_{il}+\frac{2\Delta t \rho^{n,c} c_s^2}{\theta_1^{n,c}}G_{il}^{\kp,c}.
\end{equation}
$\Sigma_{il}^{\kp,c}$ is a symmetric three by three matrix that can easily be inverted using a direct method. More precisely, since $\Sigma_{il}^{\kp,c}$ is a symmetric positive definite square matrix, its determinant can be estimated as 
	\begin{equation}
		\left|\mathbf{I} + \frac{2\Delta t c_s^2}{\theta_1^{n,c}} \mathbf{G}^{n,c} \right|^{\frac{1}{3}}\geq |\mathbf{I}|^{\frac{1}{3}}+\frac{2\Delta t c_s^2}{\theta_1^{n,c}}|\mathbf{G}^{n,c}|^{\frac{1}{3}}\geq 1+\frac{2\Delta t c_s^2}{\theta_1^{n,c}}\left(\frac{\rho^{n,c}}{\rho_0}\right)^{\frac{2}{3}}>0,
	\end{equation}
and thus its inverse exists.

Multiplication of \eqref{eqn.devG.3} by $\ISigma_{li}$ leads to an expression for $\devG^{**,c}$ in terms of the only unknown $v_m^{**,p}$:
\begin{eqnarray}
	\devG_{lk}^{**,c} 
	%&=& \ISigma_{li}\devG_{ik}^{n,c}
	%- \Delta t  \ISigma_{li}\left(G_{im}^{\kp,c}\partial_k^{cp} v_m^{**,p}+G_{mk}^{\kp,c}\partial_i^{cp} %v_m^{**,p}-\frac{2}{3}\delta_{ik}G_{nm}^{\kp,c}\partial_n^{cp}v_m^{**,p}\right)  %\nonumber \\
	%+\frac{2\Delta t c_s^2}{3 \theta_1}  \ISigma_{li}G_{mn}^{\kp,c}\devG_{nm}^{\kp,c}\delta_{ik} \nonumber \\
	&=& \ISigma_{li}\devG_{ik}^{n,c}
	- \Delta t  \ISigma_{li}\left(G_{im}^{\kp,c}\partial_k^{cp} v_m^{**,p}+G_{mk}^{\kp,c}\partial_i^{cp} v_m^{**,p}-\frac{2}{3}\delta_{ik}G_{nm}^{\kp,c}\partial_n^{cp}v_m^{**,p}\right)+\frac{2\Delta t \rho^{n,c} c_s^2}{3 \theta_1^{n,c}}  \ISigma_{lk}G_{mn}^{\kp,c}\devG_{nm}^{\kp,c}. \nonumber \\
	\label{eqn.devG.4}
\end{eqnarray}
Using a similar reasoning we may write an expression for $J_k^{***,c}$ starting from \eqref{eqn.J.disc.3}
\begin{equation}
	\left(1+\frac{\Delta t}{\tau_2}\right)J_k^{***,c}=J_k^{**,c}-\Delta t J_m^{**,c}\partial_k^{cp}v_m^{**,p},
	\label{eqn.J.disc.4}
\end{equation}
and so
\begin{equation}
	J_k^{***,c}=\frac{1}{1+\frac{\Delta t}{\tau_2}}J_k^{**,c}-\frac{\Delta t}{1+\frac{\Delta t}{\tau_2}} J_m^{**,c}\partial_k^{cp}v_m^{**,p}.
	\label{eqn.J.disc.5}
\end{equation}
Using \eqref{eqn.J.disc.5} in \eqref{eqn.si.omega} we get
\begin{equation}
	\omega_{ik}^{***,c}=\frac{\rho^{n,c} c_h^2}{1+\frac{\Delta t}{\tau_2}}J_i^{**,c}J_k^{**,c}-\frac{\Delta t \rho^{n,c} c_h^2}{1+\frac{\Delta t}{\tau_2}}J_i^{**,c}J_m^{**,c}\partial_k^{cp}v_m^{**,p}
	\label{eqn.si.omega.2}
\end{equation}
We are now ready to formally substitute Eqs. \eqref{eqn.si.omega.2} and \eqref{eqn.devG.4} via \eqref{eqn.si.sigma} into the discrete momentum \eqref{eqn.mom.disc.3}
\begin{eqnarray}
	\rho^{*,p}v_i^{**,p}&=& \rho^{*,p}v_i^{*,p}-\Delta t \, \partial^{pc}_k \left(  \sigma_{ik}^{**,c} + \omega_{ik}^{***,c} \right), \nonumber \\
	&=&  \rho^{*,p}v_i^{*,p}-\Delta t \partial_k^{pc} \rho^{n,c} c_s^2 G_{il}^{\kp,c}\devG_{lk}^{**,c}-\Delta t \partial_k^{pc} \rho^{n,c} c_h^2 J_i^{**,c}J_k^{***,c} \nonumber \\
	&=& \rho^{*,p}v_i^{*,p}-\Delta t \partial_k^{pc} \rho^{n,c} c_s^2 G_{il}^{\kp,c}\ISigma_{la}\devG_{ak}^{*,c}-\frac{2\Delta t^2 \rho^{n,c} c_s^4}{3 \theta_1^{n,c}} \partial_k^{pc} G_{il}^{\kp,c} \ISigma_{lk}G_{mn}^{\kp,c}\devG_{nm}^{\kp,c}-\Delta t \partial_k^{pc} \rho^{n,c} c_h^2 J_i^{**,c}J_k^{***,c} \nonumber \\
	&& +\Delta t^2 \partial_k^{pc} \rho^{n,c} c_s^2 G_{il}^{\kp,c}\ISigma_{al}\left(G_{am}^{\kp,c}\partial_k^{cp} v_m^{**,p}+G_{mk}^{\kp,c}\partial_a^{cp} v_m^{**,p}-\frac{2}{3}\delta_{ak}G_{nm}^{\kp,c}\partial_n^{cp}v_m^{**,p}\right) \nonumber \\
	&=& \rho^{*,p}v_i^{*,p}-\Delta t \partial_k^{pc} \rho^{n,c} c_s^2 G_{il}^{\kp,c}\ISigma_{la}\devG_{ak}^{*,c}-\frac{2\Delta t^2 \rho^{n,c} c_s^4}{3 \theta_1^{n,c}} \partial_k^{pc} G_{il}^{\kp,c} \ISigma_{lk}G_{mn}^{\kp,c}\devG_{nm}^{\kp,c} \nonumber \\
	&&+\Delta t^2 \partial_k^{pc} \rho^{n,c} c_s^2 G_{il}^{\kp,c}\ISigma_{al}\left(G_{am}^{\kp,c}\partial_k^{cp} v_m^{**,p}+G_{mk}^{\kp,c}\partial_a^{cp} v_m^{**,p}-\frac{2}{3}\delta_{ak}G_{nm}^{\kp,c}\partial_n^{cp}v_m^{**,p}\right) \nonumber \\
	&&-\partial_k^{pc}\frac{\Delta t \rho^{n,c} c_h^2}{1+\frac{\Delta t}{\tau_2}}J_i^{**,c}J_k^{**,c}
	+\frac{\Delta t^2 \rho^{n,c} c_h^2}{1+\frac{\Delta t}{\tau_2}} \partial_k^{pc}J_i^{**,c}J_m^{**,c}\partial_k^{cp}v_m^{**,p}.
	\label{eq.linsys.devGJv}
\end{eqnarray}
Now, one can recognize a linear system for the only unknown, the new velocity field $v^{**,p}$. We may eventually write in a compact form the coefficient matrix using a 4-rank tensor. Indeed one can define a new tensor
\begin{equation} 
	H_{ik\,nm}^{\kp,c} = C_{ik\,nm}^{\kp,c} + \frac{\tau_2 \rho^{n,c} c_h^2}{\Delta t+\tau_2}\delta_{kn}J_i^{**,c}J_m^{**,m}\label{eqn.H.r4.1}
\end{equation}
with 
\begin{equation}
	C_{ik\,nm}^{\kp,c} = \rho^{n,c} c_s^2 G_{il}^{\kp,c}\ISigma_{al}\left(G_{am}^{\kp,c}\delta_{kn}+G_{mk}^{\kp,c}\delta_{an}-\frac{2}{3}\delta_{ak}G_{nm}^{\kp,c}\right). \label{eqn.C.r4.1} 
\end{equation}
Hence, the final vector wave equation for the discrete velocity field reads 
\begin{equation}
	\rho^{*,p}v_i^{**,p} - \Delta t^2 \partial_k^{pc} H_{ik\,nm}^{\kp,c} \partial_n^{cp}v_m^{**,p} =  b_i^{\kp,c},  
	\label{eqn.velocity.wave} 
\end{equation}
with the known right hand side 
\begin{equation}
	b_i^{\kp,c} = \rho^{*,p}v_i^{*,p}-\Delta t \partial_k^{pc} \rho^{n,c} c_s^2 G_{il}^{\kp,c}\ISigma_{la}\devG_{ak}^{*,c}-\frac{2\Delta t^2 \rho^{n,c} c_s^4}{3 \theta_1^{n,c}} \partial_k^{pc} G_{il}^{\kp,c} \ISigma_{lk}G_{mn}^{\kp,c}\devG_{nm}^{\kp,c}-\partial_k^{pc}\frac{\Delta t \rho^{n,c} c_h^2}{1+\frac{\Delta t}{\tau_2}}J_i^{**,c}J_k^{**,c}. 
\end{equation}

\subsubsection{Asymptotic limit: the viscous Navier-Stokes  stress tensor}
In this section we want to analyze the case when $\tau_1\rightarrow 0$. Without losing of generality we consider only viscous effects, i.e. $c_h=0$. When $\tau_1\rightarrow 0$ then $\theta_1^{n,c} = \rho^{n,c} \frac{\tau_1 c_s^2}{3}|\G^{n,c}|^{-\frac{5}{6}}\rightarrow 0$ and then
\begin{equation}
	\ISigma_{il}=\left(\delta_{il}+\frac{2\Delta t  \rho^{n,c} c_s^2}{\theta_1^{n,c}}G_{il}^{\kp,c}\right)^{-1}\approx \left(\frac{2\Delta t  \rho^{n,c} c_s^2}{\theta_1^{n,c}}G_{il}^{\kp,c}\right)^{-1}=\frac{\theta_1^{n,c}}{2\Delta t  \rho^{n,c} c_s^2}\left(G_{il}^{\kp,c}\right)^{-1} = \frac{1}{6} \frac{\tau_1}{\Delta t} |\G^{n,c}|^{-\frac{5}{6}} \left(G_{il}^{\kp,c}\right)^{-1},
	\label{eq.approx1.tau0}
\end{equation}
since in the limit the term $\frac{2\Delta t \rho^{n,c} c_s^2}{\theta_1^{n,c}}$ dominates $\delta_{il}$, which then results just a perturbation. Furthermore, when $\tau_1\rightarrow 0 $ the Hilbert expansion of the discrete solution 
\begin{equation}
	\G^{**,c} = \G^{**,c}_0 + \tau_1 \G^{**,c}_1 + \mathcal{O}(\tau_1^2)
\label{eqn.hilbert.G} 
\end{equation}	
inserted into \eqref{eqn.metric.disc.3} leads to the immediate result  
\begin{equation}
	\qquad \mathring{\G}^{**,c}_0 = 0, \qquad \textnormal{ and therefore } \qquad \G^{**,c}_0 = g^{**,c} \, \mathbf{I} + \mathcal{O}(\tau_1),  
	\label{eq.approx2.tau0}
\end{equation}
with $g^{**,c}=|\G^{**,c}|^{\frac{1}{3}} = (\rho^{**,c}/\rho_0)^{2/3}$, see \cite{GPRmodel} for the continuous case. Hence, from now on, we will always assume $\G^{n,c}_0 = |\G^{n,c}|^{\frac{1}{3}} \, \mathbf{I} = (\rho^{n,c}/\rho_0)^{2/3} \, \mathbf{I}$. 
%For the analysis it is useful to compute the following quantity:
%\begin{equation}
%	\frac{\rho \theta_1}{2}=\frac{\rho \theta_1 c_s^2}{6|G|^{\frac{5}{6}}}=\frac{\frac{\rho}{\rho_0}\rho_0 \tau_1 c_s^2}{6 %|G|^{\frac{5}{6}}}=\frac{|G|^{\frac{1}{2}}\mu}{|G|^{\frac{5}{6}}}=\frac{\mu}{g}
%	\label{eq.approx3.tau0}
%\end{equation}
%where $\mu=\frac{1}{6}\rho_0 c_s^2 \tau_1$ is the viscosity coefficient as obtained from the formal asymptotic analysis at the continuum %level, see again \cite{GPRmodel}.
Using the approximations \eqref{eq.approx1.tau0}-\eqref{eq.approx2.tau0} in the limit we can simplify our discrete vector wave equation for the velocity field \eqref{eq.linsys.devGJv} as follows: 
\begin{eqnarray}
	\rho^{*,p}v_i^{**,p} \!&=& \! b_i^{n,p}+\Delta t^2 \partial_k^{pc} \frac{\rho^{n,c} c_s^2 \tau_1}{6\Delta t} |\G^{n,c}|^{-\frac{5}{6}}  G_{il}^{\kp,c}\left(G_{al}^{\kp,c}\right)^{-1}\left(G_{am}^{\kp,c}\partial_k^{cp} v_m^{**,p}+G_{mk}^{\kp,c}\partial_a^{cp} v_m^{**,p}-\frac{2}{3}\delta_{ak}G_{nm}^{\kp,c}\partial_n^{cp}v_m^{**,p}\right) \nonumber \\
	&=& b_i^{n,p}+\Delta t \partial_k^{pc} \frac{1}{6} \rho^{n,c} c_s^2 \tau_1 \overbrace{|\G^{n,c}|^{-\frac{5}{6}} |\G^{n,c}|^{\frac{1}{3}}}^{|\G^{n,c}|^{-\halb} = \rho_0 /  \rho^{n,c}}  \delta_{ia}\left(\delta_{am}\partial_k^{cp} v_m^{**,p}+\delta_{mk}\partial_a^{cp} v_m^{**,p}-\frac{2}{3}\delta_{ak}\delta_{nm}\partial_n^{cp}v_m^{**,p}\right)  \nonumber \\
	&=&  b_i^{n,p}+\Delta t \partial_k^{pc} \frac{1}{6} \rho_0 c_s^2 \tau_1  \left(\partial_k^{cp} v_i^{**,p}+\partial_i^{cp}v_k^{**,p}-\frac{2}{3}\delta_{ik}\partial_m^{cp}v_m^{**,p}\right),
\end{eqnarray}
from which the classical compressible Navier-Stokes stress tensor with viscosity coefficient $\frac{1}{6} \rho_0 c_s^2 \tau_1$ becomes totally evident also at the fully discrete level.  
%Furthermore
%\begin{equation}
%	b_i^{\kp,c}= \rho^{*,p}v_i^{*,p}-\Delta t \partial_k^{pc} \rho c_s^2 G_{il}^{\kp,c}\ISigma_{la}\devG_{ak}^{*,c}-\frac{2\Delta t^2 \rho %c_s^4}{3\theta_1} \partial_k^{pc} G_{il}^{\kp,c} \ISigma_{lk}G_{mn}^{\kp,c}\devG_{nm}^{\kp,c}= \rho^{*,p}v_i^{*,p}
%\end{equation}
%since $\devG_{ij}=0$.
%which is the classical Navier-Stokes stress tensor. It worth to be noted that we automatically obtain the parabolic limit when $\tau \rightarrow 0$ also at the discrete level, furthermore the resulting system for $v^{**,c}$ is simply
%\begin{equation}
%	\rho^{*,p}v_i^{**,p}=\rho^{*,p}v_i^{*,p}+\Delta t \partial_k^{pc}\mu \left(\partial_k^{cp} %v_i^{**,p}+\partial_i^{cp}v_k^{**,p}-\frac{2}{3}\delta_{ik}\partial_m^{cp}v_m^{**,p}\right).
%\end{equation}
%that is a consistent semi-implicit discretization of the advection-viscous contribution for the compressible Navier-Stokes equations.

\subsubsection{Asymptotic limit: The Fourier heat conduction}
Regarding the formal Fourier heat conduction limit $\tau_2 \to 0$, let us consider well-prepared initial data for the thermal impulse given by $J_k^{n,c} = \tau_2 J_{k,1}^{n,c} + \mathcal{O}(\tau_2^2)$ following \cite{GPRmodel}. Then, after the convective step we have $J_k^{*,c} = \tau_2 J_{k,1}^{*,c} + \mathcal{O}(\tau_2^2)$. Note that the arithmetic averaging does not change the asymptotic expansion and thus $J_k^{*,p}$ is well-prepared as well. 
Using the Hilbert expansion, i.e. assuming in the discrete subsystem $J_k^{**,p} = J_{k,0}^{**,p} + \tau_2 J_{k,1}^{**,p} + \mathcal{O}(\tau_2^2)$, from \eqref{eqn.heatflux.disc.2} one gets the following conditions $J_{k,0}^{**,p} = 0$ and $J_{k,1}^{**,p} = - \partial_k^{cp} T^{**,c}$.
Therefore, we have $J_k^{**,p} = - \tau_2 \partial_k^{cp} T^{**,c} + \mathcal{O}(\tau^2_2)$. Inserting this relation into the discrete heat flux, we obtain 
\begin{equation}
	q_k^{**,p} = \rho^{*,p} \, T^{**,p} \, c_h^2  \tau_2  \, \partial_k^{pc} T^{**,c} + \mathcal{O}(\tau_2^2) = \lambda^p \, \partial_k^{pc} T^{**,c} + \mathcal{O}(\tau_2^2).
\end{equation} 
Moreover, in the discrete G-J-v subsystem by \eqref{eqn.J.disc.4} using the Hilbert expansion of $J_k^{***,c}$ one obtains immediately $J_k^{***,c} = \mathcal{O}(\tau_2)$.   
In the pressure subsystems presented in the sequel, the final update $J_k^{n+1,c}$ is well-prepared with the analogue argumentation. Therefore, in the stiff relaxation limit $\tau_2 \to 0$ we recover the discrete Fourier heat flux independently of $\Delta t$. 
\subsection{Discrete pressure subsystem}
Finally, with analogue argumentation as above, the discrete pressure subsystem is given by
\begin{subequations}
	\begin{align}
		(\rho v_i)^{n+1,p}  &= (\rho v_i)^{**,p} - \Delta t \, \partial^{pc}_i p^{n+1,c}  , \label{eqn.mom.disc.4} \\
		\frac{1}{\gamma - 1} p^{n+1,c} &= \frac{1}{\gamma - 1} p^{**,c} - \Delta t \,
		\partial^{cp}_i \left( h^{**,p} (\rho v_i)^{n+1,p} \right) . \label{eqn.energy.disc.4}
	\end{align}
\end{subequations}
where the total energy reduces under the ideal gas law to the contribution $\mathcal{E}_1 = p/(\gamma - 1)$.
Substitution of \eqref{eqn.mom.disc.4} into \eqref{eqn.energy.disc.4} yields the following classical discrete scalar wave equation for the pressure
\begin{equation}
	\frac{1}{\gamma - 1} p^{n+1,c}  - \Delta t^2
	\partial^{cp}_i \left( h^{**,p} \, \partial^{pa}_i p^{n+1,a} \right)
	= \frac{1}{\gamma - 1} p^{**,c} -   \Delta t \,
	\partial^{cp}_i \left( h^{**,p} \, (\rho v_i)^{**,p} \right). \label{eqn.wave.p}
\end{equation}
The coefficient matrix of this system is symmetric and positive definite, and \eqref{eqn.wave.p} can therefore be solved easily via a classical matrix-free conjugate gradient method.
The final momentum is then given by \eqref{eqn.mom.disc.4} and the scheme is completed by formally setting $\rho^{n+1,p} = \rho^{*,p}$, while the distorsion field and the thermal impulse are updated via the following compatible structure-preserving discretization that preserves curl-free fields exactly at the discrete level in the absence of source terms, by virtue of the use of compatible discrete nabla operators, see also \cite{SIGPR,MaxwellGLM}:
\begin{equation}
   A_{ik}^{n+1,c} = A_{ik}^{n,c} - \Delta t \, \partial_k^{cp} \left( v_m^{n+1,p} A_{im}^{n,p} \right) - \frac{\Delta t}{4} \sum \limits_{p \in \mathcal{P}_c} v_m^{n+1,p} \left( \partial_m^{pa} A_{ik}^{n,a} - 
   \partial_k^{pa} A_{im}^{n,a} \right) + \frac{\Delta t}{\theta_1} \alpha_{ik}^{n+1,c},   	 
	\label{eqn.SP.A}    
\end{equation}
\textcolor{black}{which is a compatible discretization of eqn. \eqref{eqn.deformation}, see also \cite{SIGPR}}.
In order to enforce the compatibility with the determinant constraint we may uniformly rescale all components of $A_{ik}^{n+1,c}$ so that 
$|A_{ik}^{n+1,c}| = \rho^{n+1,c} / \rho_0$ holds. For alternative compatible discretizations, see \cite{BoscheriGPRGCL,HTCLagrangeGPR}.
\begin{equation}
   J_{k}^{n+1,c} = J_{k}^{n,c} - \Delta t \, \partial_k^{cp} \left( v_m^{n+1,p} J_{m}^{n,p} + T^{**,p} \right) - \frac{\Delta t}{4} \sum \limits_{p \in \mathcal{P}_c}  v_m^{n+1,p} \left( \partial_m^{pa} J_{k}^{n,a} - 
\partial_k^{pa} J_{m}^{n,a} \right) + \frac{\Delta t}{\theta_2} \beta_{ik}^{n+1,c},   	 	
	\label{eqn.SP.J}    
\end{equation}
\textcolor{black}{which is a compatible discretization of eqn. \eqref{eqn.heatflux}, see also \cite{SIGPR}}.
The discrete total energy is finally updated as 
\begin{equation}
	\mathcal{E}^{n+1,c} = \mathcal{E}^{**,c} - \Delta t \, \partial_k^{cp} \left( h^{**,p} \left( \rho v_k \right)^{n+1,p} + (\sigma_{ik}^{n+1,p} + \omega_{ik}^{n+1,p}) v^{n+1,p}_i\right). 
	\label{eqn.energy}    	
\end{equation}

\subsection{Summary of the algorithm}\label{sec.summary}

In the following we summarize the complete algorithm of the proposed structure-preserving four-split scheme. 

\subsubsection{Convective subsystem}

The first step is the computation of the explicit convective terms 
\begin{equation}
	\q^{*,c} = \q^{n,c} - \frac{\Delta t}{\vc} \sum_{p \in \mathcal{P}_c} n_k^{pc} \,   
	\f^{p}_k,
	\label{sum.eqn.fv}
\end{equation}
with a nodal Rusanov-type flux: 
\begin{equation}
	\f_k^{p} = \bar{\f}_k^{p} 
	- \halb h^p s^p_{\max} \partial_k^{pc} q^{n,c}, \qquad \bar{\f}_k^{p} = \frac{1}{|\mathcal{C}_p|} \sum \limits_{c \in \mathcal{C}_p}  \f_k^c 
	\label{sum.eqn.rusanov}
\end{equation}
with $h^p = \frac{4|\Omega_p|}{|\partial \Omega_p|} = \frac{2\Delta x \Delta y}{\Delta x + \Delta y}$ a characteristic length scale and $s^p_{\max} = \max \limits_c(|\v^c|)$ the maximum convective speed at the vertex.

\subsubsection{Temperature subsystem}

Using the result of the previous step, one then needs to solve the temperature subsystem which accounts for the heat conduction, but not yet for the contribution of $\J$ to the stress tensor in the momentum equation. Find the new temperature by solving 
\begin{equation}
	\left( 1+\frac{\Delta t}{\tau_2} \right) \frac{c_v}{T^{n,c} c_h^2} \, T^{**,c}  - \Delta t^2
	\partial^{cp}_k \left( \partial^{pa}_k T^{**,a} \right)
	= \left( 1+\frac{\Delta t}{\tau_2} \right) \frac{c_v}{T^{n,c} c_h^2} \, T^{*,c} - \Delta t \,
	\partial^{cp}_k \left( J_k^{*,p} \right), \label{sum.eqn.wave.T}
\end{equation}
and update the auxiliary thermal impulse as 
\begin{subequations}
	\begin{align}\label{sum.eqn.heat.disc.final}
		J^{**,c}_k &= \left(\frac{\tilde \theta_2}{\tilde \theta_2 + \Delta t }\right)(J^{*,c}_k - \Delta t \partial_k^{cp} T^{**,p})
	\end{align}
\end{subequations}
as well as the heat flux contribution to the total energy 
\begin{equation}
	\mathcal{E}^{**,c} = \mathcal{E}^{*,c} - \Delta t \, \partial_k^{cp}q_k^{**,p}.
\end{equation}

\subsubsection{G-J-v subsystem}

We now can solve the coupled system of the discrete momentum equation, the discrete thermal impulse equation and the discrete equation for $\mathring{\G}$, i.e. the deviator of $\G$. Using a discrete equation for $\mathring{\G}$ is a \textit{major key ingredient} of the method presented in this paper. Find the new velocity vector by solving the discrete vector wave equation for the velocity 
\begin{equation}
	\rho^{*,p}v_i^{**,p} - \Delta t^2 \partial_k^{pc} H_{ik\,nm}^{\kp,c} \partial_n^{cp}v_m^{**,p} =  b_i^{\kp,c}, 
	\label{sum.eqn.velocity.wave} 
\end{equation}
with 
\begin{equation}
	H_{ik\,nm}^{\kp,c} = C_{ik\,nm}^{\kp,c} + \frac{\theta_2 \rho c_h^2}{\Delta t+\theta_2}\delta_{kn}J_i^{**,c}J_m^{**,m}, 
	\qquad 
		C_{ik\,nm}^{\kp,c} = \rho c_s^2 G_{il}^{\kp,c}\ISigma_{al}\left(G_{am}^{\kp,c}\delta_{kn}+G_{mk}^{\kp,c}\delta_{an}-\frac{2}{3}\delta_{ak}G_{nm}^{\kp,c}\right),
		\label{sum.eqn.H.r4.1}
\end{equation}
and 
\begin{equation}
	b_i^{\kp,c} = \rho^{*,p}v_i^{*,p}-\Delta t \partial_k^{pc} \rho c_s^2 G_{il}^{\kp,c}\ISigma_{la}\devG_{ak}^{*,c}-\frac{2\Delta t^2 \rho c_s^4}{3\theta_1} \partial_k^{pc} G_{il}^{\kp,c} \ISigma_{lk}G_{mn}^{\kp,c}\devG_{nm}^{\kp,c}-\partial_k^{pc}\frac{\Delta t \rho c_h^2}{1+\frac{\Delta t}{\theta_2}}J_i^{**,c}J_k^{**,c}.
\end{equation}

\subsubsection{Pressure subsystem}

The fourth and last subsystem that needs to be solved is the pressure subsystem. The discrete pressure wave equation reads
\begin{equation}
	\frac{1}{\gamma - 1} p^{n+1,c}  - \Delta t^2
	\partial^{cp}_i \left( h^{**,p} \, \partial^{pa}_i p^{n+1,a} \right)
	= \frac{1}{\gamma - 1} p^{**,c} -   \Delta t \,
	\partial^{cp}_i \left( h^{**,p} \, (\rho v_i)^{**,p} \right), \label{sum.eqn.wave.p}
\end{equation}
which for $M_a \to 0$ reduces to the pressure Poisson equation for the incompressible Euler and Navier-Stokes equations and which thus allows to update the momentum in a  manner that is compatible with the low Mach number limit: 
\begin{equation}
		(\rho v_i)^{n+1,p}  = (\rho v_i)^{**,p} - \Delta t \, \partial^{pc}_i p^{n+1,c}.  \label{sum.eqn.mom} 
\end{equation}

\subsubsection{Compatible update of the distortion field, of the thermal impulse and of the total energy}

In the absence of source terms the distortion field and the thermal impulse must remain curl-free for all times if they were initially curl-free. This is achieved via \textcolor{black}{a compatible discretization of the original equations for $A_{ik}$ and $J_k$, i.e. eqns. \eqref{eqn.deformation} and \eqref{eqn.heatflux}, see also \cite{SIGPR}}:  
\begin{equation}
   A_{ik}^{n+1,c} = A_{ik}^{n,c} - \Delta t \, \partial_k^{cp} \left( v_m^{n+1,p} A_{im}^{n,p} \right) - \frac{\Delta t}{4} \sum \limits_{p \in \mathcal{P}_c} v_m^{n+1,p} \left( \partial_m^{pa} A_{ik}^{n,a} - 
\partial_k^{pa} A_{im}^{n,a} \right) + \frac{\Delta t}{\theta_1} \alpha_{ik}^{n+1,c},   	 
	\label{sum.eqn.SP.A}   	 
\end{equation}
\begin{equation}
   J_{k}^{n+1,c} = J_{k}^{n,c} - \Delta t \, \partial_k^{cp} \left( v_m^{n+1,p} J_{m}^{n,p} + T^{**,p} \right) - \frac{\Delta t}{4} \sum \limits_{p \in \mathcal{P}_c}  v_m^{n+1,p} \left( \partial_m^{pa} J_{k}^{n,a} - 
\partial_k^{pa} J_{m}^{n,a} \right) + \frac{\Delta t}{\theta_2} \beta_{ik}^{n+1,c}.   	 	 	 	
	\label{sum.eqn.SP.J}   	 
\end{equation}
Finally, the new discrete total energy reads 
\begin{equation}
	\mathcal{E}^{n+1,c} = \mathcal{E}^{**,c} - \Delta t \partial_k^{cp} \left( h^{**,p} \left( \rho v_k \right)^{n+1,p} + (\sigma_{ik}^{n+1,p} + \omega_{ik}^{n+1,p}) v^{n+1,p}_i\right). 
	\label{sum.eqn.energy}   	 	
\end{equation}
\textcolor{black}{The proposed scheme is a finite volume method and is therefore consistent with the integral form of the conservation laws. Indeed, it conserves mass, momentum and total energy exactly at the discrete level. Unfortunately, we were not yet able to prove a discrete entropy inequality for this scheme, which is clearly a shortcoming. Further research in this direction will be the subject of future work. }

% % % % % % % % % % % % % % % % % % % % % % % % % % % % % %
%            Numerical experiments                        %
% % % % % % % % % % % % % % % % % % % % % % % % % % % % % %

\section{Numerical results}\label{sec.results}

\subsection{Taylor-Green vortex at low Mach number}

\begin{figure}[!htbp]
	\begin{center}
		\begin{tabular}{cc} 
			\includegraphics[width=0.45\textwidth]{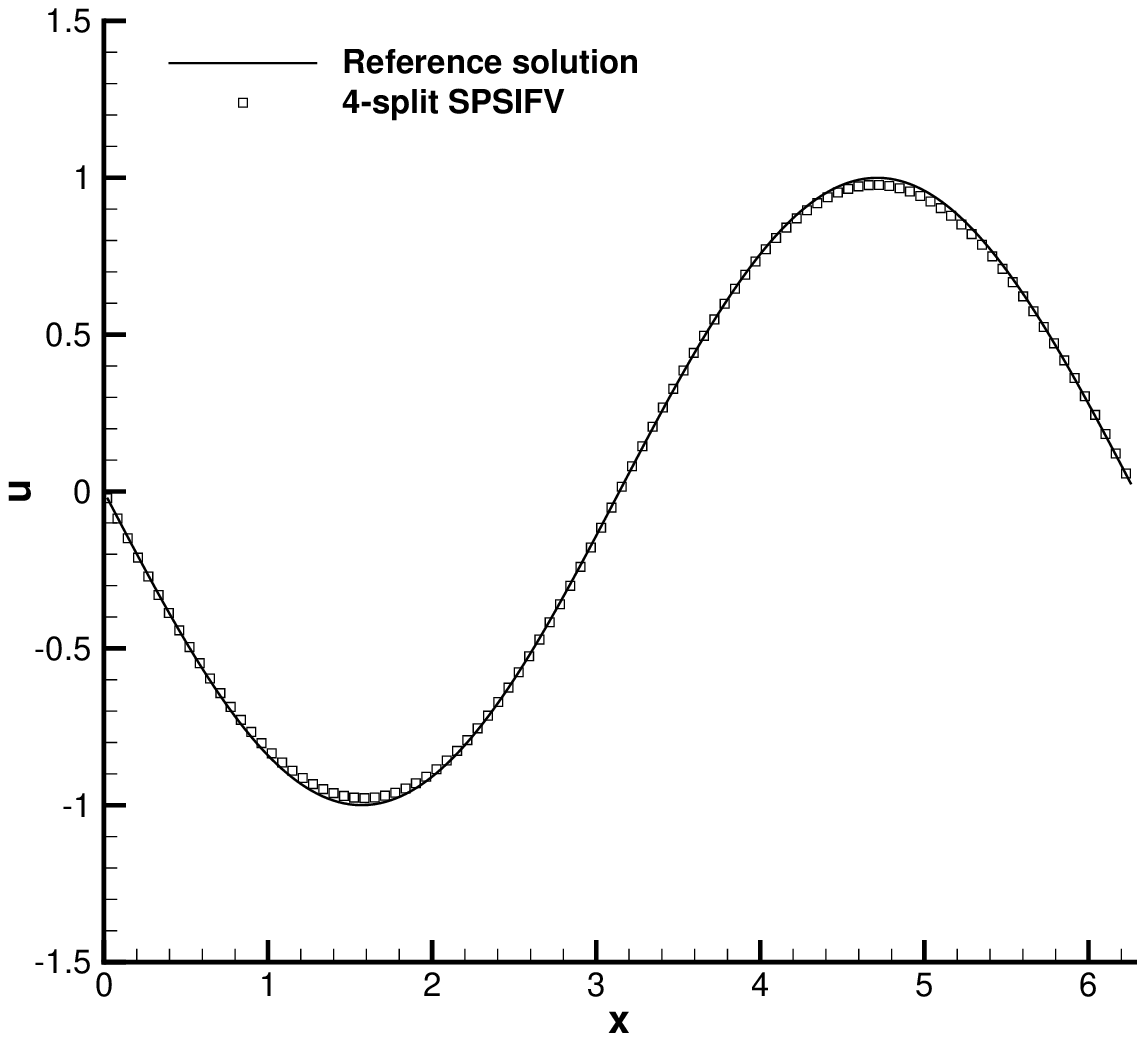}   & 
			\includegraphics[width=0.45\textwidth]{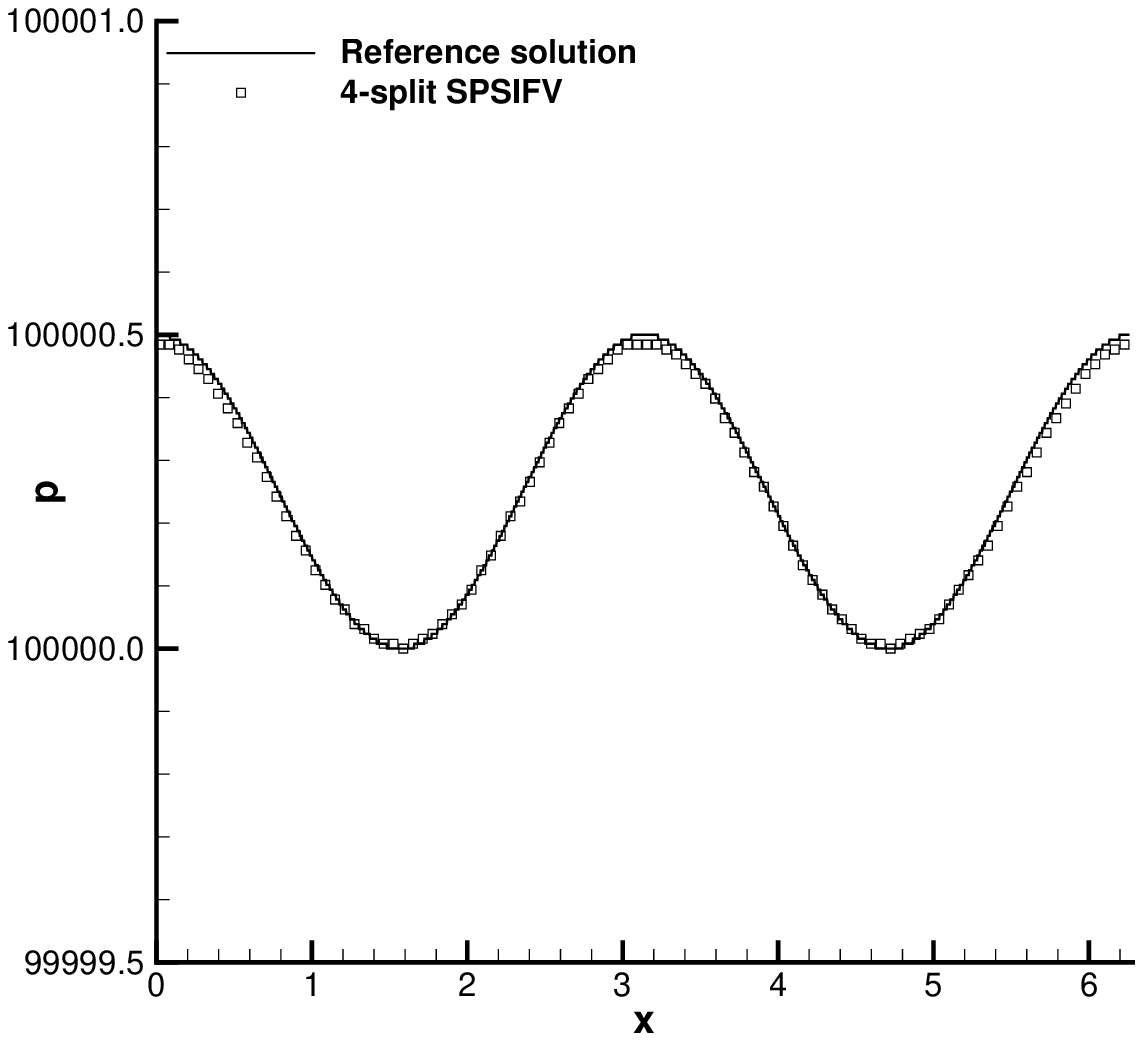}  
		\end{tabular} 
		\caption{Numerical solution of the GPR model for the Taylor-Green vortex problem in the low Mach number limit and in the stiff relaxation limit with $\tau_1 = 10^{-8}$ and $\tau_2 = 10^{-10}$ at time $t=1.0$ using the new semi-implicit 4-split scheme. 1D cuts along the $x$ axis and comparison with the exact solution of the incompressible Navier-Stokes equations for the velocity components $u$ (left) and the pressure $p$ (right).  } 
		\label{fig.tgv}
	\end{center}
\end{figure}

The two-dimensional Taylor-Green vortex
\begin{eqnarray}
	\rho(x,y,t) &=& \rho_0, \label{eq:TG_rho} \\
	u(x,y,t)&=&\sin(x)\cos(y)e^{-2\nu t}, \label{eq:TG_0} \\
	v(x,y,t)&=&-\cos(x)\sin(y)e^{-2\nu t}, \label{eq:TG_1} \\
	p(x,y,t)&=& p_0 + \frac{1}{4}(\cos(2x)+\cos(2y))e^{-4\nu t},
	\label{eq:TG_2}
\end{eqnarray}
is an exact solution of the incompressible Navier-Stokes equations, i.e. it is a solution of the compressible equations in the low Mach number limit $M_a \to 0$.
The computational domain is set to $\Omega=[0,2\pi]^2$ with periodic boundaries everywhere.
We solve the GPR model with the new 4-split structure preserving semi-implicit finite volume (SPSIFV) scheme presented in this paper up to a final time of $t=1.0$ on a computational grid composed of $200 \times 200$ cells. The following parameters are used here:
$\gamma=1.4$, $\rho_0=1$, $\tau_1=10^{-8}$, $\tau_2=10^{-10}$, $c_v=c_p/\gamma$ with $c_p=1004$, $c_s=1000$, $c_h=100$. The initial condition for the velocity and for the pressure is \eqref{eq:TG_0}--\eqref{eq:TG_2}, respectively, with $p_0=10^5$, hence the
acoustic Mach number in this test problem is $M_a=0.0027$, while the shear Mach number is $M_s=10^{-3}$. The distortion field is initially set to $\mathbf{A}=\mathbf{I}$ and the initial specific thermal impulse is $\mathbf{J}=0$. 
The numerical solution together with a reference solution is given in Figure \ref{fig.tgv}. We observe a good agreement of our numerical solution with the exact solution of the incompressible Navier-Stokes equations. 

We now show numerical evidence that the proposed semi-implicit 4-split scheme has the expected behaviour towards the low Mach number limit, i.e. that $\rho \to \rho_0$ for $M_a \to 0$. We therefore solve the Taylor-Green vortex problem again with increasing background pressure $p_0$ on a \textit{fixed} mesh of $128 \times 128$ cells and using the same time step based on the convective CFL condition for all Mach regimes. 
It is well-known \cite{KlaMaj,KlaMaj82,Klein2001,KleinMach,MunzPark,MunzDumbserRoller} that the density fluctuations should scale with $M_a^2$, which is a direct consequence of the scaling of the divergence of the velocity, that should also scale with $M_a^2$. In Table \ref{tab.Machconv} we report the $L^2$ and $L^\infty$ errors of the density and the $L^\infty$ error of the divergence of the velocity field as a function of $p_0$ and the acoustic Mach number $M_a$. One can clearly observe the expected second order convergence in terms of $M_a$ for the density and divergence errors in all norms. 
\begin{table}  
		\caption{Convergence results of the new 4-split semi-implicit FV scheme obtained for the GPR model in the low Mach number limit. The error norms refer to the density $\rho$ in $L^2$ and $L^\infty$ norm as well as to the divergence of the velocity in $L^\infty$ norm at a final time of $t=0.1$. The observed convergence orders refer to the acoustic Mach number $M_a$ and the density errors and clearly show second order convergence in the Mach number. } 	
		\begin{center} 
			\renewcommand{\arraystretch}{1.1}
			\begin{tabular}{cccccccc} 
				\hline
				$p_0$ & $M_a$ & ${L^2}(\rho)$ & ${L^\infty}(\rho)$ &  ${L^\infty}(\nabla \cdot \mathbf{v})$ & $\mathcal{O}^2_M(\rho)$ & $\mathcal{O}^\infty_M(\rho)$ 
				& $\mathcal{O}^\infty_M(\nabla \cdot \mathbf{v})$ \\ 
				\hline
                    1.00E+02 &	8.45E-02  &	9.1182E-05 & 2.0115E-04	& 4.4641E-03 &	   &     &     \\ 
                    1.00E+03 &	2.67E-02  &	9.7651E-06 & 2.1414E-05	& 2.8035E-04 & 1.9 & 1.9 & 2.4 \\
                    1.00E+04 &	8.45E-03  &	9.7860E-07 & 2.1482E-06	& 2.8032E-05 & 2.0 & 2.0 & 2.0 \\
                    1.00E+05 &	2.67E-03  &	9.7860E-08 & 2.1482E-07	& 2.8034E-06 & 2.0 & 2.0 & 2.0 \\
                    1.00E+06 &	8.45E-04  &	9.7860E-09 & 2.1493E-08	& 2.8185E-07 & 2.0 & 2.0 & 2.0 \\
                    1.00E+07 &	2.67E-04  &	9.7858E-10 & 2.1480E-09	& 2.8029E-08 & 2.0 & 2.0 & 2.0 \\
                    1.00E+08 &	8.45E-05  &	9.7835E-11 & 2.1461E-10	& 2.8003E-09 & 2.0 & 2.0 & 2.0 \\
                    1.00E+09 &	2.67E-05  &	9.7653E-12 & 2.1265E-11	& 2.7896E-10 & 2.0 & 2.0 & 2.0 \\
                    1.00E+10 &	8.45E-06  &	1.0015E-12 & 2.1417E-12	& 2.7711E-11 & 2.0 & 2.0 & 2.0 \\
                    1.00E+11 &	2.67E-06  &	9.7653E-14 & 2.1292E-13	& 2.7865E-12 & 2.0 & 2.0 & 2.0 \\
				\hline           
			\end{tabular}
		\end{center}
		\label{tab.Machconv}
	\end{table}

\subsection{Riemann problems in the fluid and solid limit}

Here we solve a set of 1D Riemann problems in the fluid and solid limit of the GPR
model. In the fluid limit for $\tau_1, \tau_2 \to 0$ the reference solution is given by the exact solution of the 
Riemann problem of the compressible Euler equations \cite{toro-book}. In the solid limit we compare with numerical results published previously \cite{SIGPR,HTCGPR}. The purpose of these tests is to show the ability of our new method to deal also with high Mach number flows and shock waves. The computational domain is given by 
$\Omega = [-0.5,+0.5]^2$ and is discretized with a computational grid of $1000 \times 10$ elements, apart from RP3 for which we employ a mesh of $2000 \times 10$ due to the large number of different waves present in this test.  
We use periodic boundary conditions in $y$ direction and Neumann boundary conditions in the $x$ direction. 
For RP1 and RP2 the parameters of the GPR model are set to $\rho_0 = 1$, $\gamma = 1.4$ and $c_s = 100$, $c_h = 10$, 
$c_v=1$, while for RP3 and RP4 we use $\rho_0 = 1$, $\gamma = 1.4$ and $c_s = 1$, $c_h = 1$, $c_v=1$. 

The initial data for $\rho$, $u$ and $v$ and $p$ are summarized in Table \ref{tab.ic.ideal}, where also the relaxation times $\tau_1$ and $\tau_2$ are provided. The remaining state variables are $w=0$, $\A = \mathbf{I}$ and $\J = 0$. 

RP1 is the classical Sod shock tube while RP2 is the well-known Riemann problem of Lax. RP3 is a Riemann problem in the solid limit, while RP4 is the same Riemann problem as RP3 but in the fluid limit. 

\begin{figure}[!htbp]
	\begin{center}
		\includegraphics[width=0.32\textwidth]{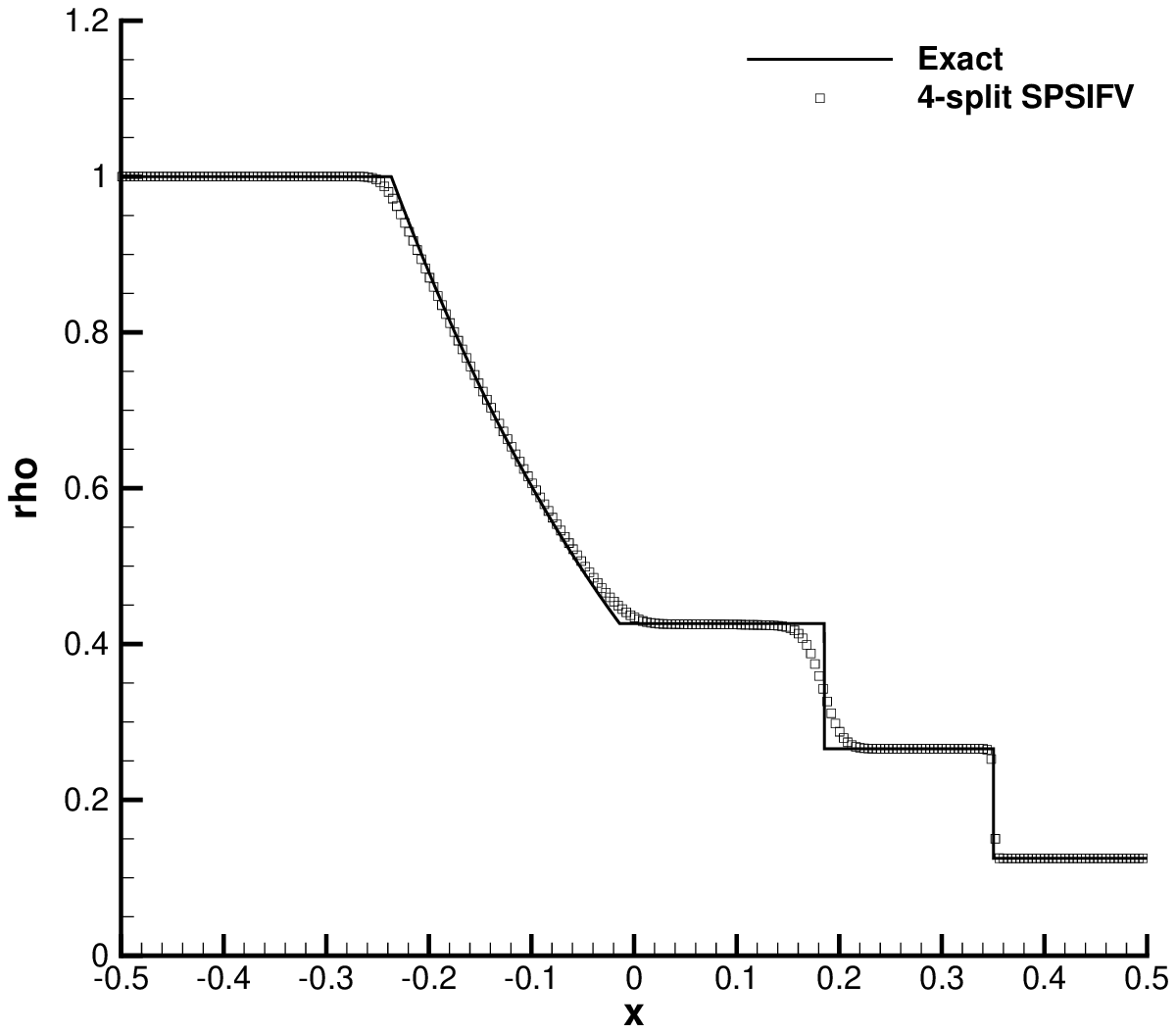}   
		\includegraphics[width=0.32\textwidth]{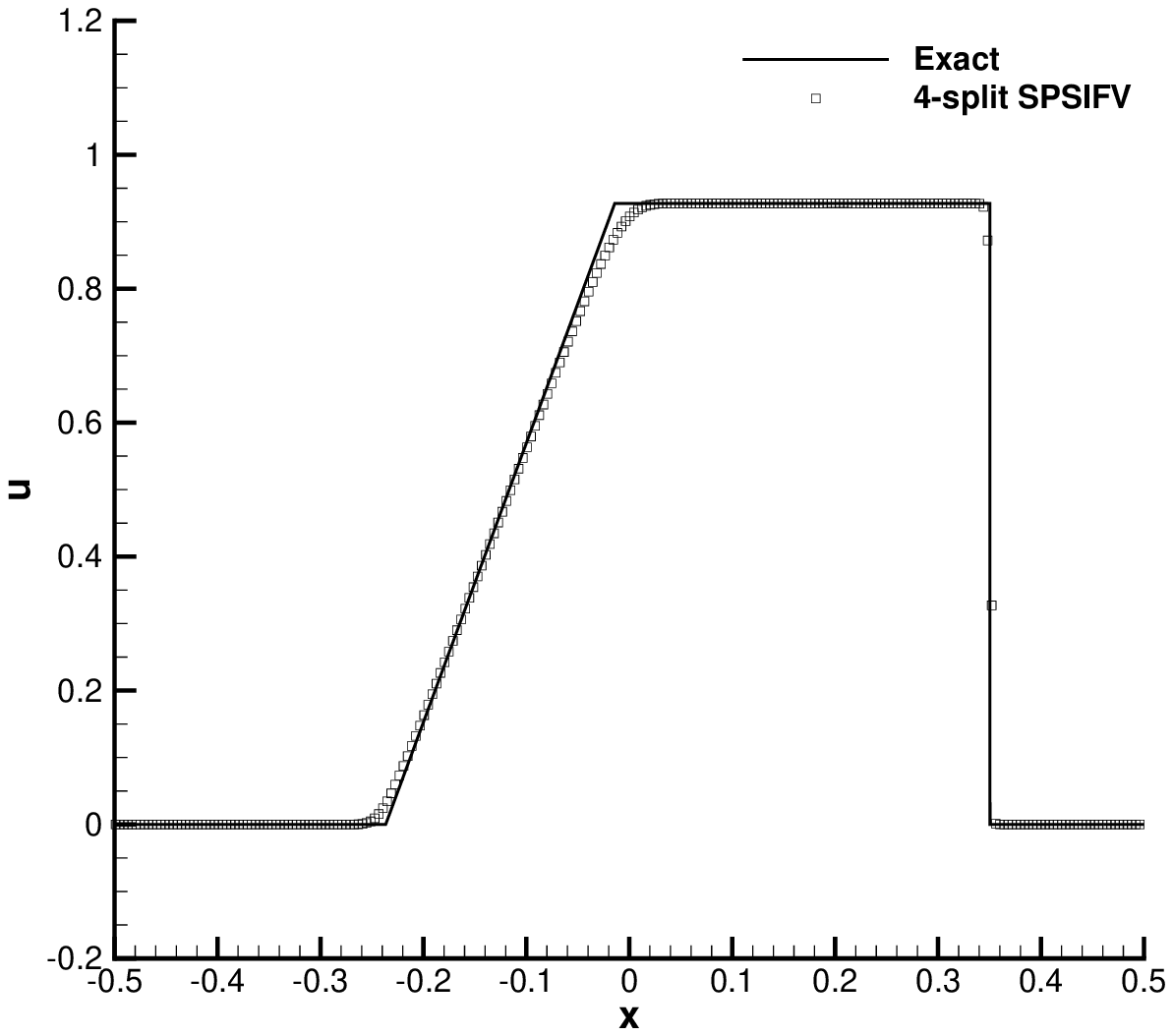}   
		\includegraphics[width=0.32\textwidth]{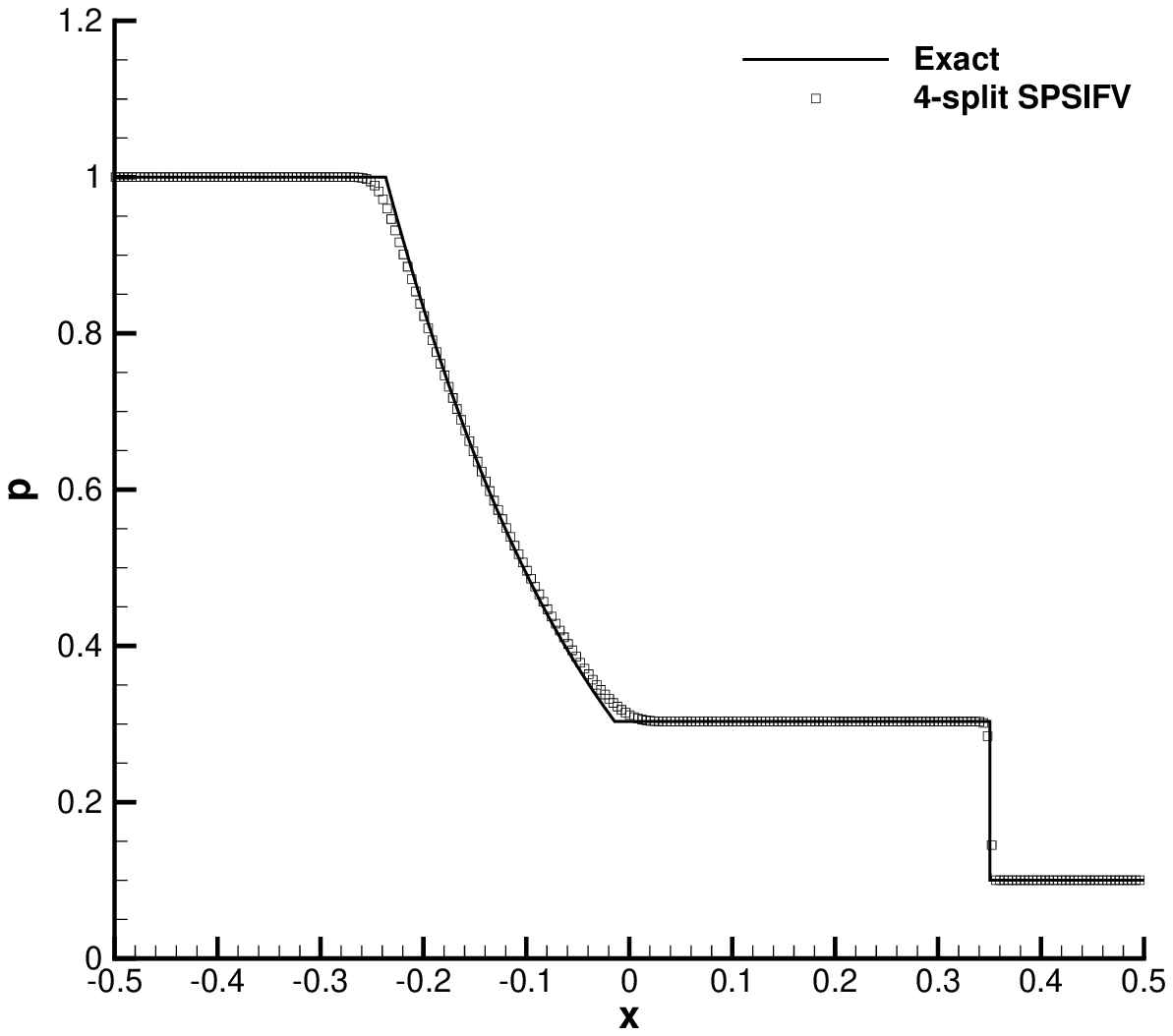}   
		\caption{Exact solution of the compressible Euler equations and numerical solution of the GPR model in the stiff relaxation limit ($\tau_1 = 10^{-10}, \tau_2 = 10^{-12}$) for Riemann problem RP1 (Sod shock tube) obtained with the new semi-implicit 4-split scheme. The density $\rho$, the velocity component $u$ and the pressure $p$ are shown at a final time of $t=0.2$.} 
		\label{fig.rp1}
	\end{center}
\end{figure}

\begin{figure}[!htbp]
	\begin{center}
		\includegraphics[width=0.32\textwidth]{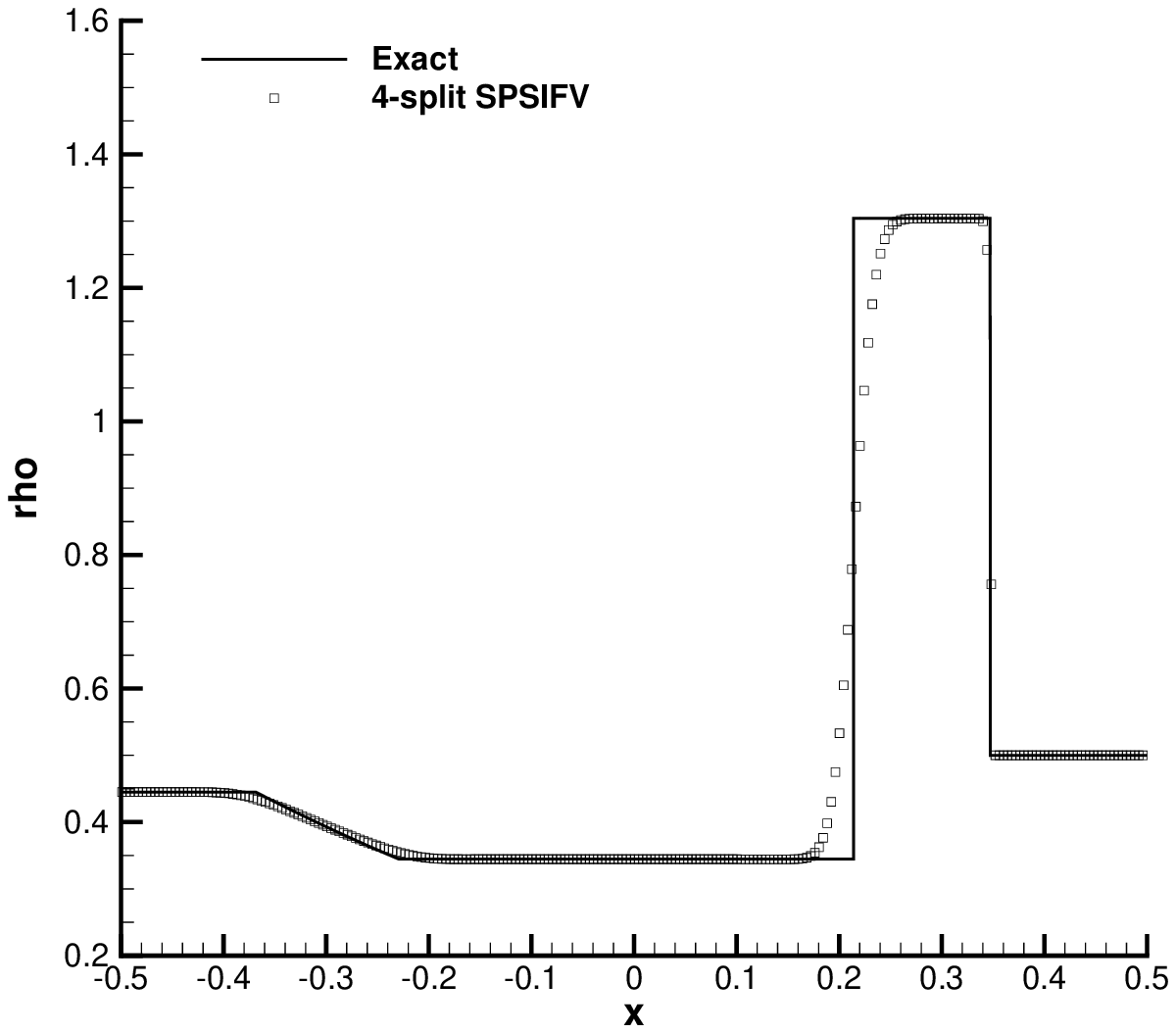}   
		\includegraphics[width=0.32\textwidth]{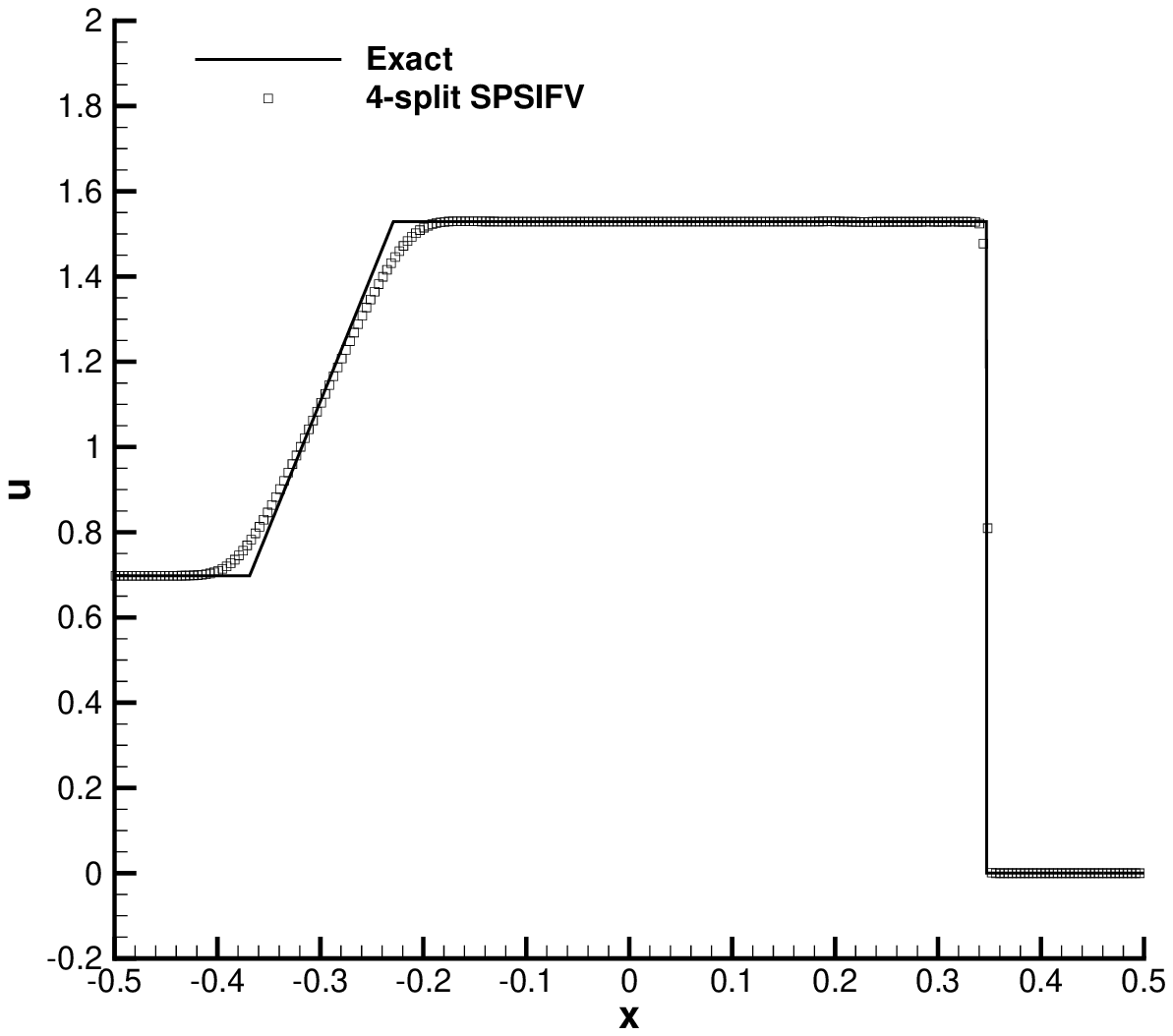}   
		\includegraphics[width=0.32\textwidth]{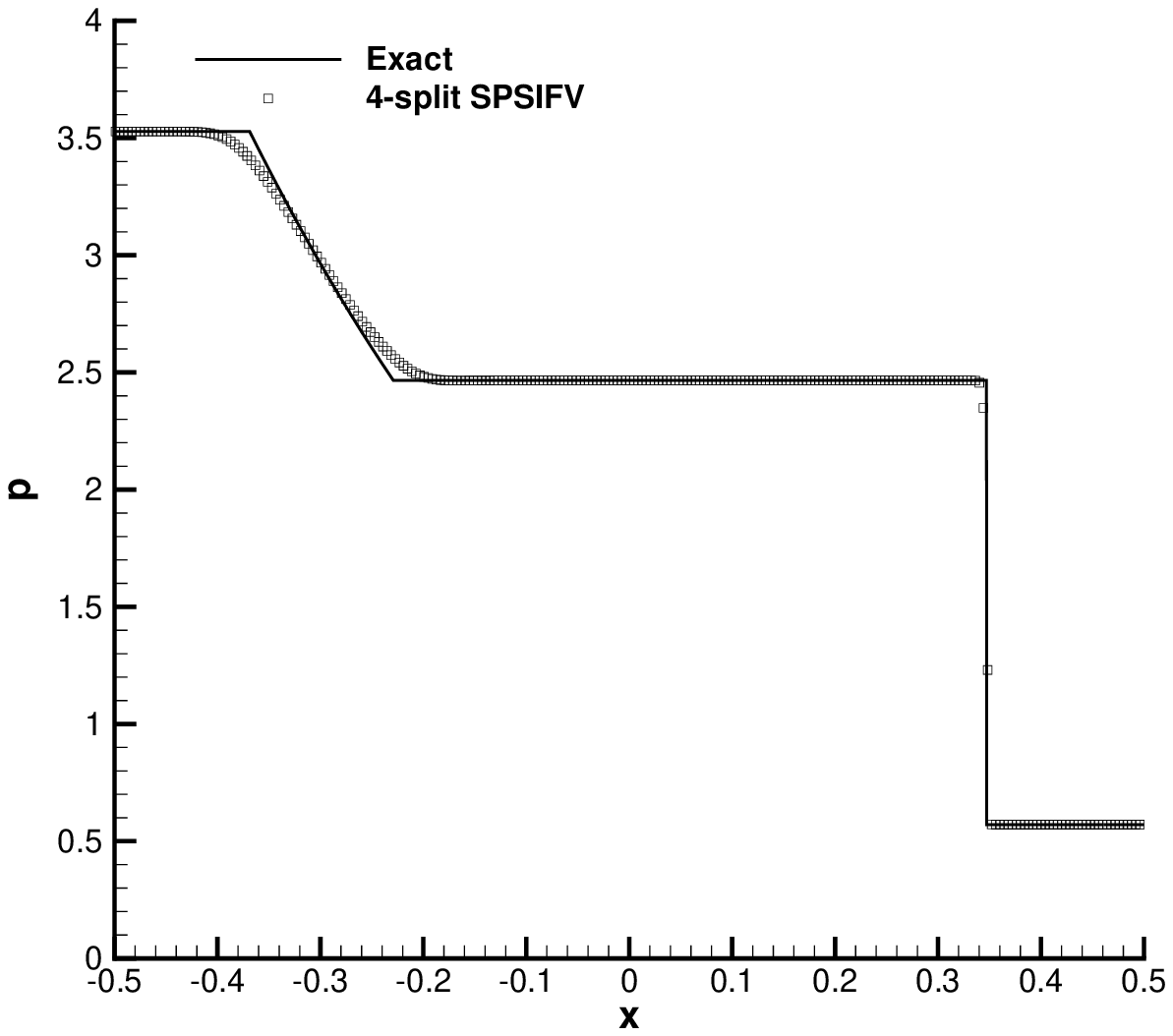}   
		\caption{Exact solution of the Euler equations and numerical solution of the GPR model in the stiff relaxation limit ($\tau_1 = 10^{-10}, \tau_2 = 10^{-12}$) for Riemann problem RP2 (Lax shock tube). The density $\rho$, the velocity component $u$ and the pressure $p$ are shown at a final time of $t=0.14$.} 
		\label{fig.rp2}
	\end{center}
\end{figure}

\begin{figure}[!htbp]
	\begin{center}
		\begin{tabular}{ccc}
		\includegraphics[width=0.32\textwidth]{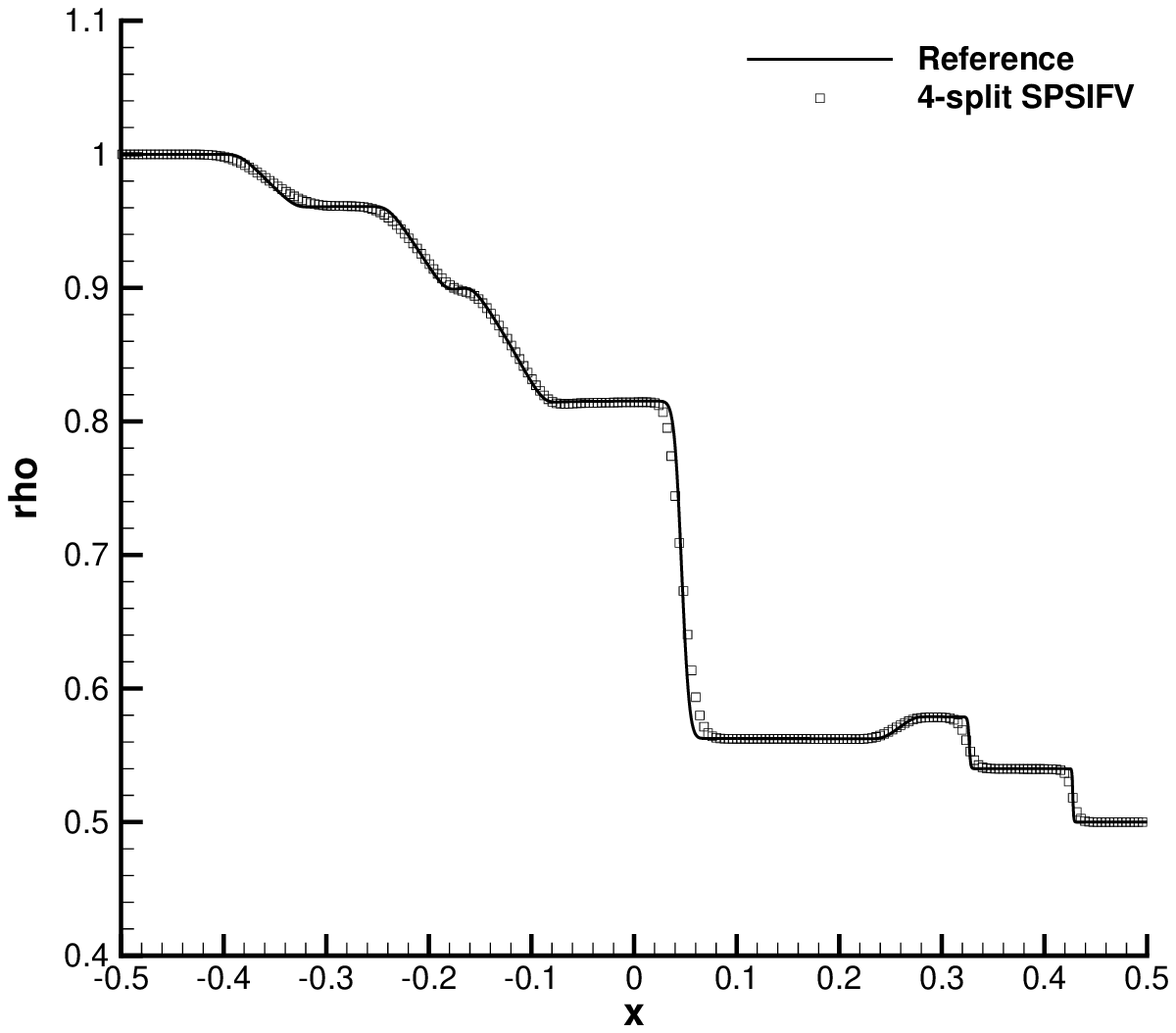}   &
		\includegraphics[width=0.32\textwidth]{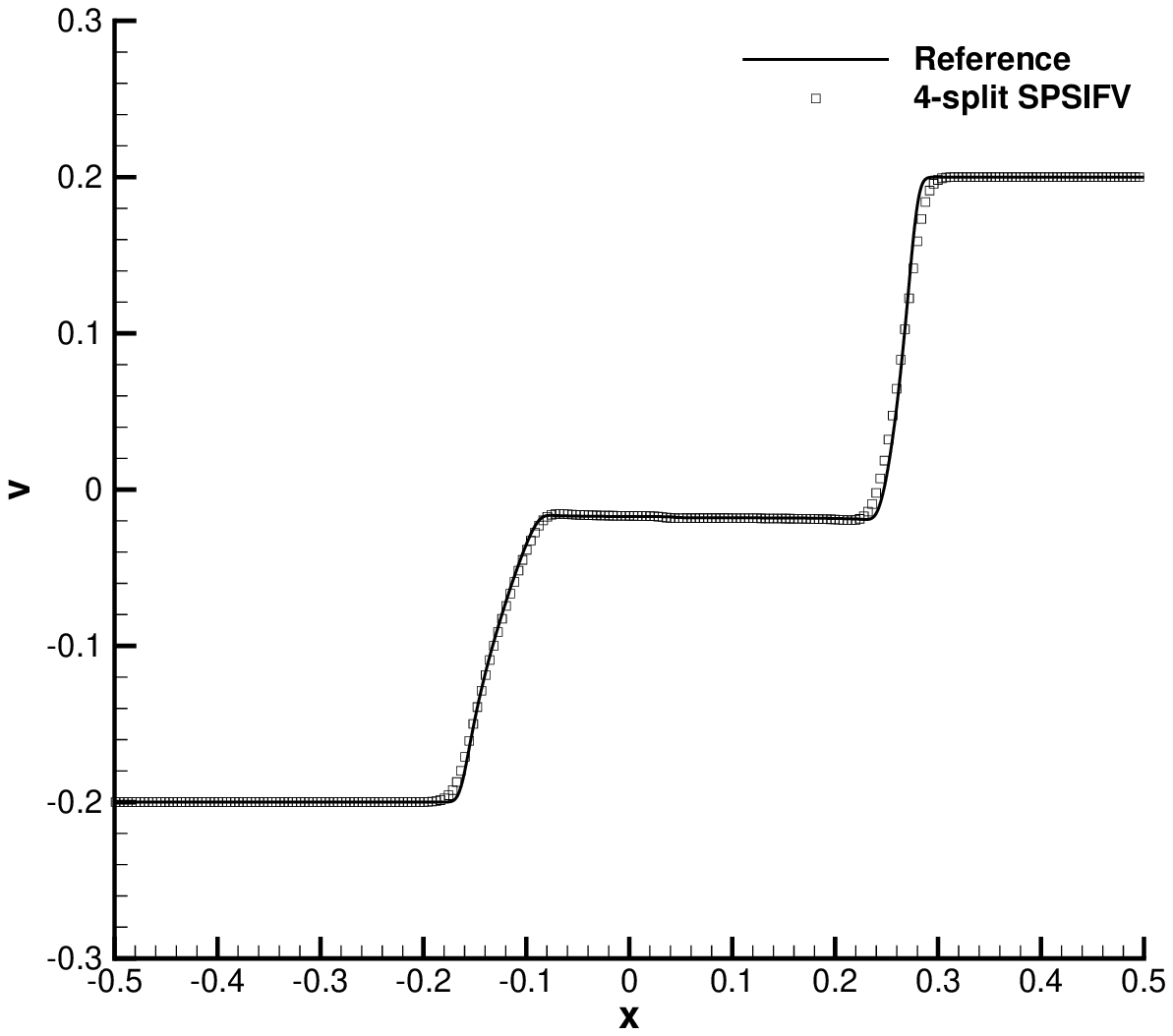}     &
		\includegraphics[width=0.32\textwidth]{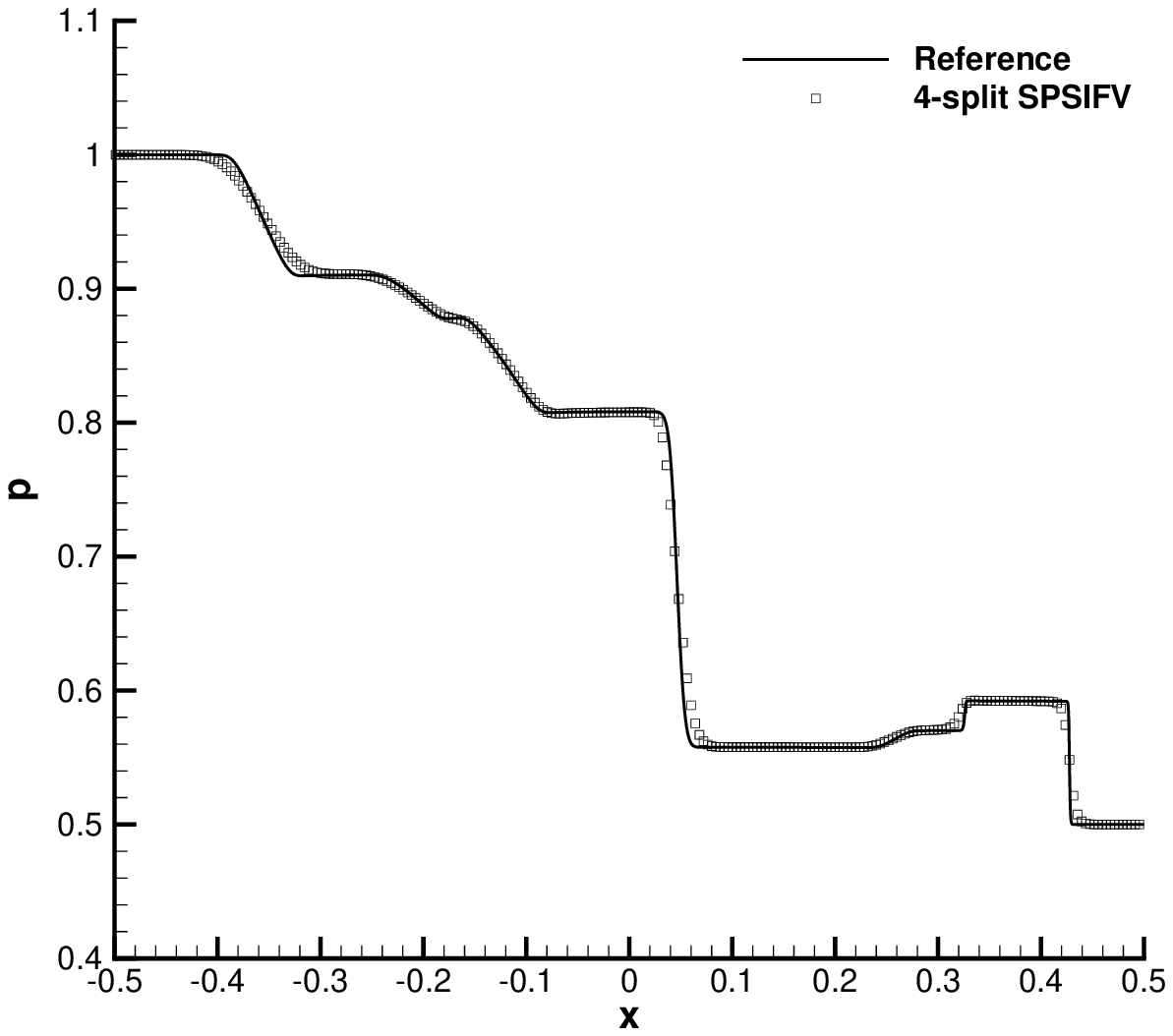}     \\ 
		\includegraphics[width=0.32\textwidth]{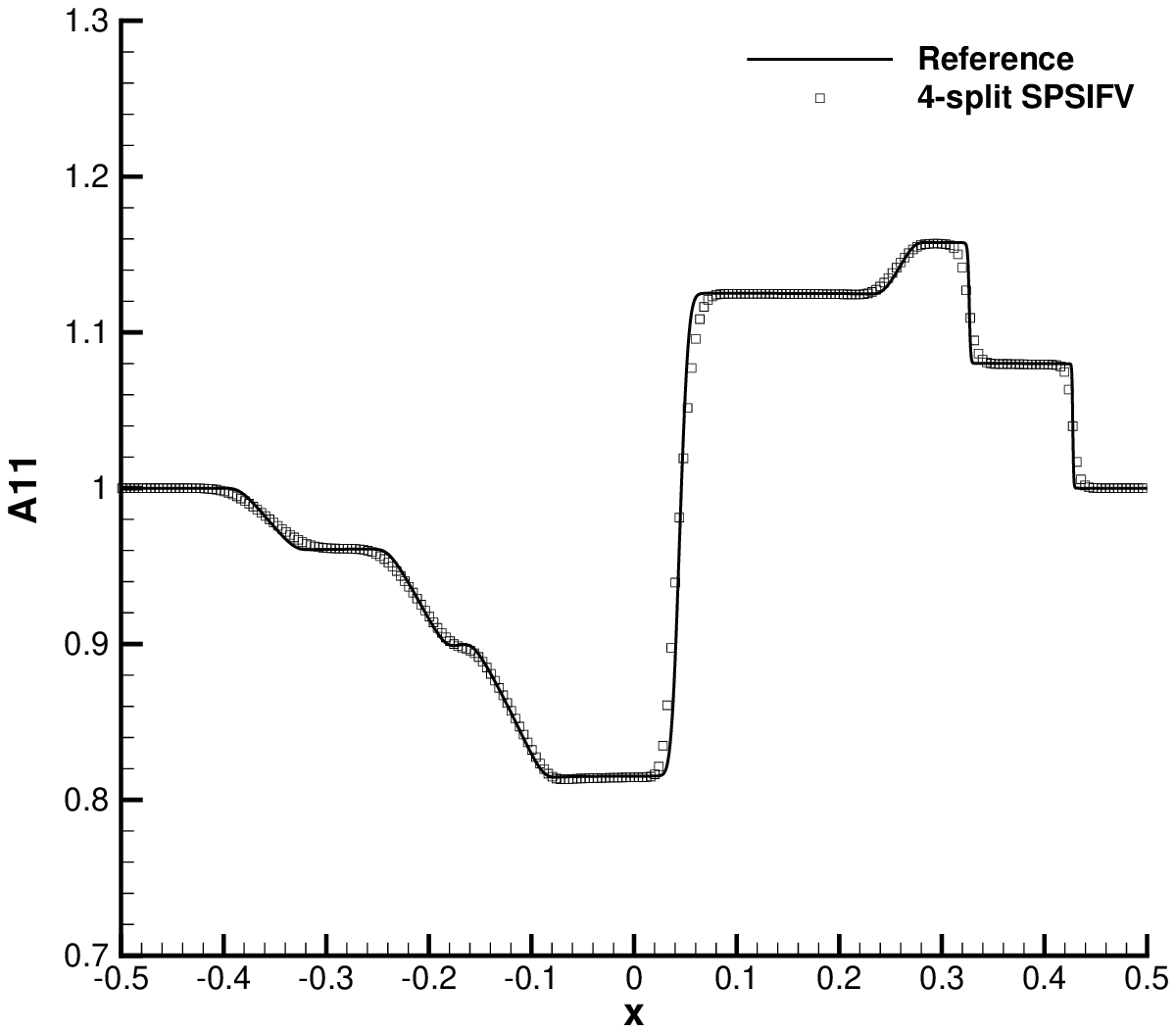}   &
		\includegraphics[width=0.32\textwidth]{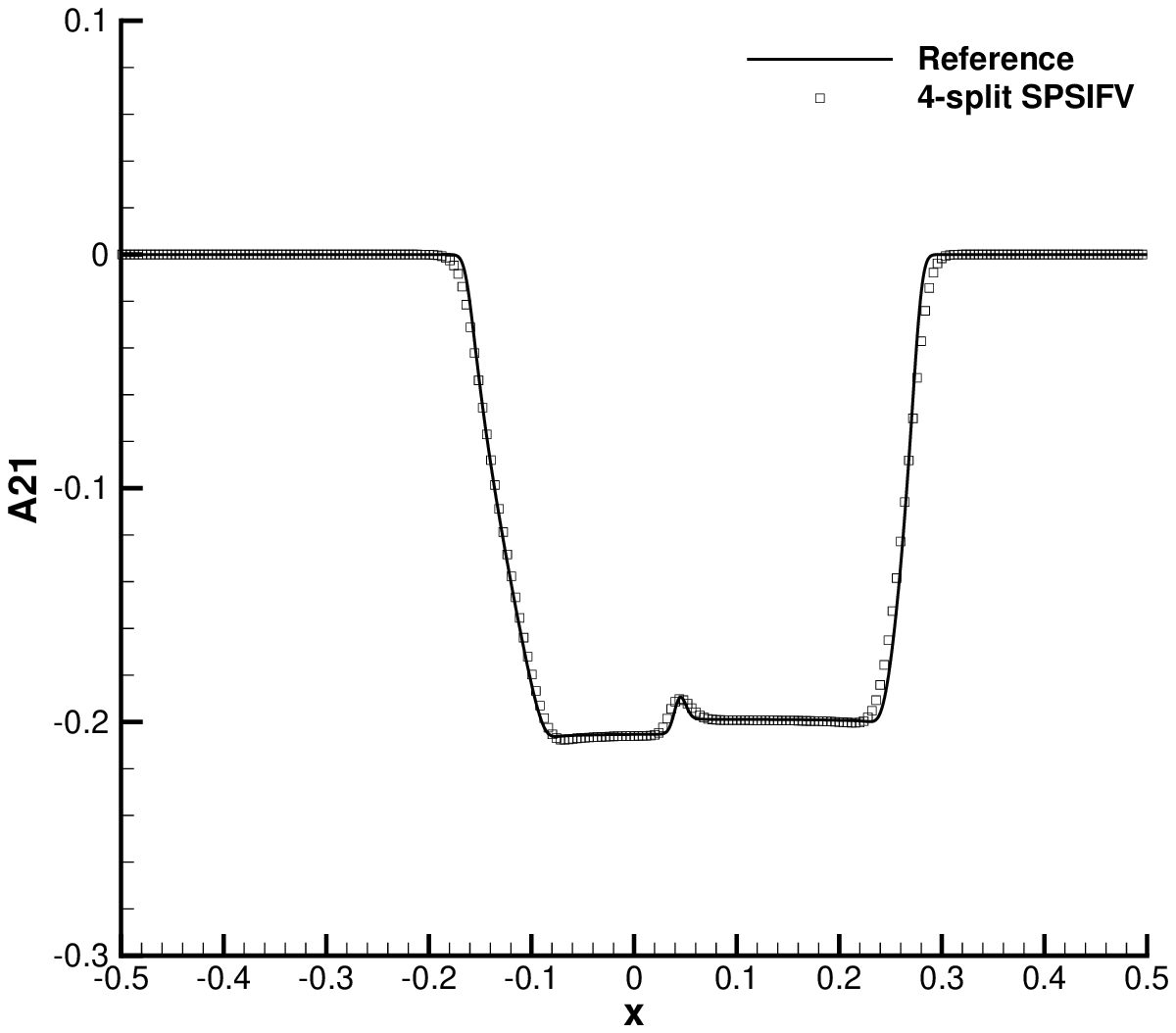}   &
		\includegraphics[width=0.32\textwidth]{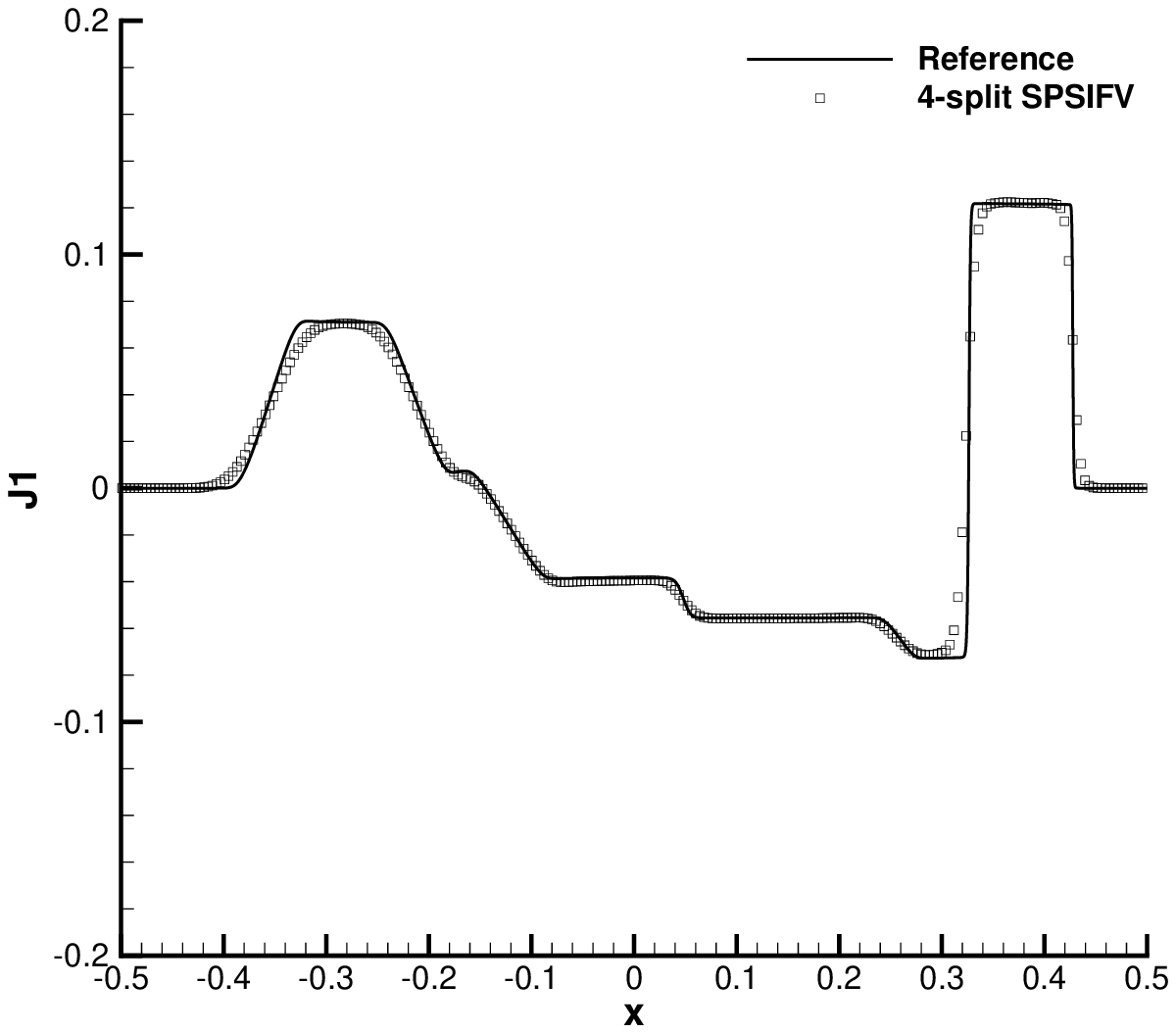}   
		\end{tabular}    
		\caption{Reference solution and numerical solution of the homogeneous GPR model without source terms ($\tau_1 = \tau_2 = 10^{20}$) for Riemann problem RP3 at a final time of $t=0.2$. Top row: density $\rho$, the velocity component $v$ and the pressure $p$. Bottom row: distorsion field components $A_{11}$, $A_{21}$ and thermal impulse component $J_1$. 
		One can note seven waves that are contained in the homogeneous part of the GPR model.}  
		\label{fig.rp3}
	\end{center}
\end{figure}

\begin{figure}[!htbp]
	\begin{center}
		\includegraphics[width=0.32\textwidth]{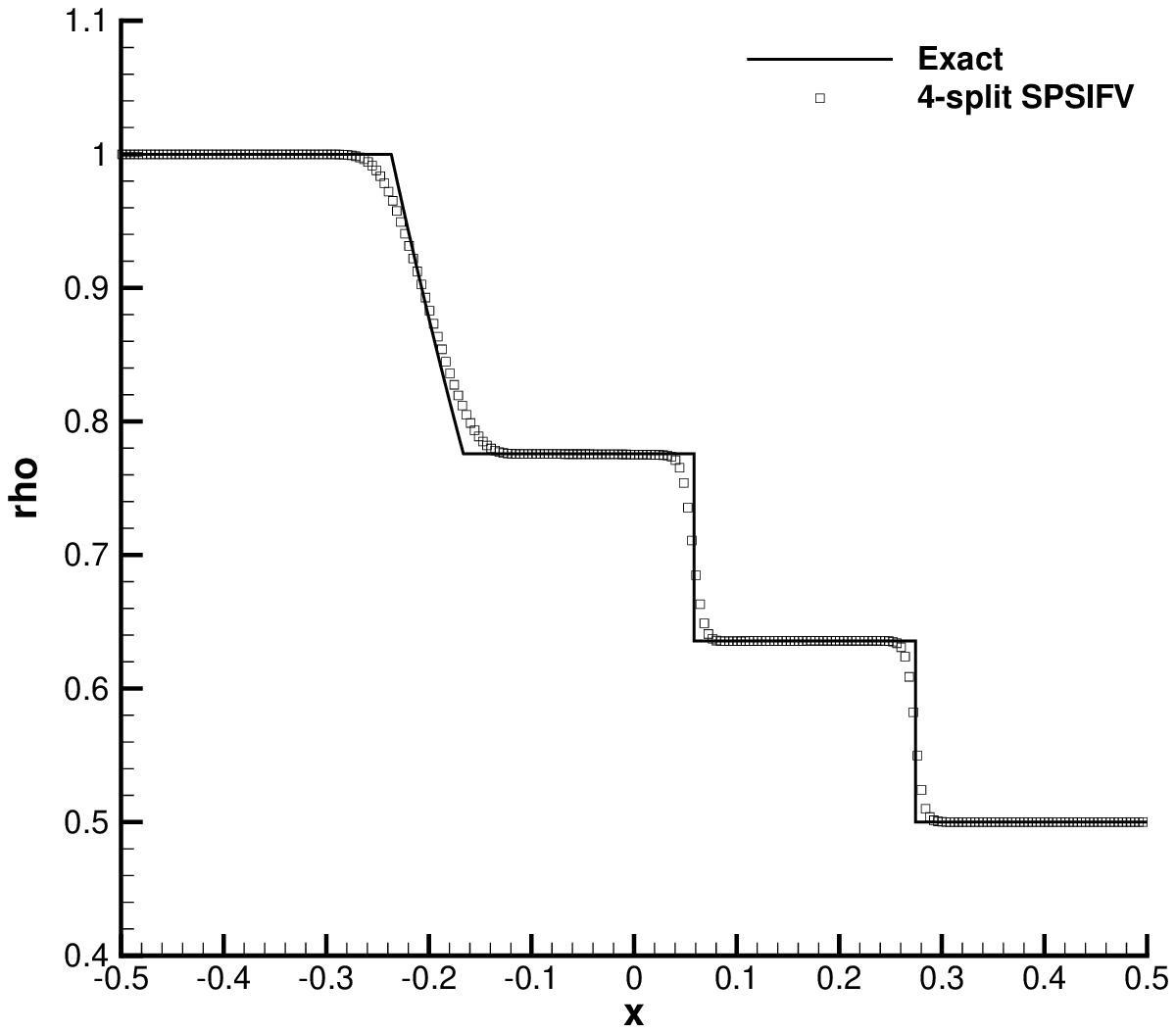}   
		\includegraphics[width=0.32\textwidth]{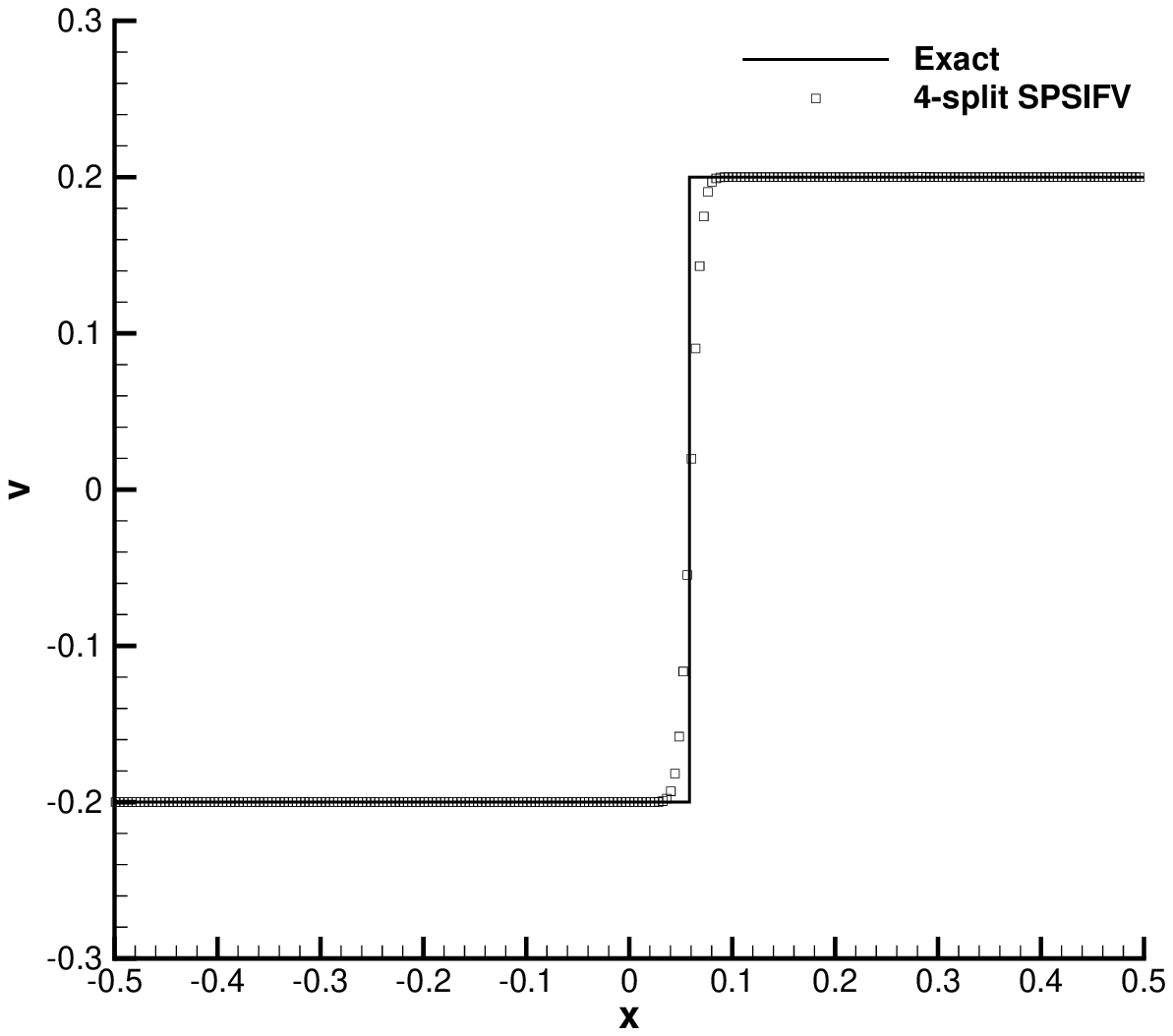}   
		\includegraphics[width=0.32\textwidth]{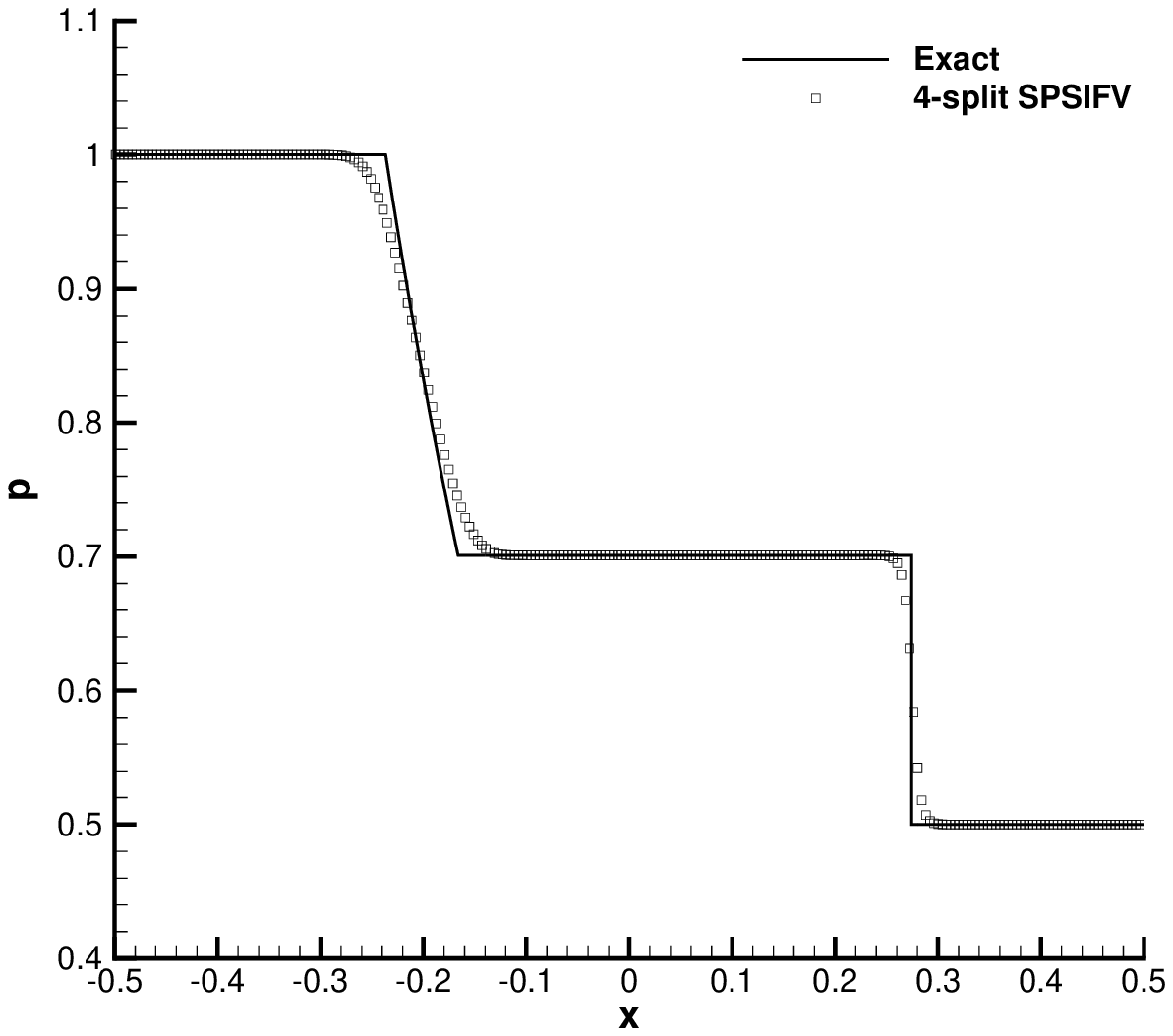}   
		\caption{Exact solution of the Euler equations and numerical solution of the GPR model in the stiff relaxation limit ($\tau_1 = 10^{-10}, \tau_2 = 10^{-12}$) for Riemann problem RP4. The density $\rho$, the velocity component $u$ and the pressure $p$ are shown at a final time of $t=0.2$.} 
		\label{fig.rp4}
	\end{center}
\end{figure}

The computational results obtained with the new semi-implicit four-split scheme are shown 
in Figures \ref{fig.rp1} - \ref{fig.rp4}. We observe a good agreement between the numerical 
and the reference solution. 

\begin{table}[!htbp]
	\renewcommand{\arraystretch}{1.25}
	\caption{Initial states left (L) and right (R) for density $\rho$, velocity $\mathbf{v}$ and pressure $p$  
		for a set of Riemann problems solved on the domain $\Omega=[-\frac{1}{2},+\frac{1}{2}]^2$ using the new 4-split 
		FV scheme. 
		The Riemann problems include the fluid limit (RP1-RP3) as well as the solid limit (RP4). The parameters $c_h$ and $c_s$ as well as the relaxation times $\tau_1$ and $\tau_2$ are also specified. In all cases we set $\gamma=1.4$. } 
	\begin{center} 
		\begin{tabular}{ccccccccccccc} 
			\hline
			RP & $\rho_L$ & $u_L$ & $v_L$ & $p_L$ & $\rho_R$ & $u_R$ & $v_R$ & $p_R$ & $c_s$ & $c_h$ & $\tau_1$ & $\tau_2$  \\ 
			\hline
			RP1 &  1.0      &  0.0       &  0.0 & 1.0     & 0.125      &  0.0        &  0.0 & 0.1      & 100 & 10 & $10^{-10}$ & $10^{-12}$ \\
			RP2 &  0.445    &  0.698     &  0.0 & 3.528   & 0.5        &  0.0        &  0.0 & 0.571    & 100 & 10 & $10^{-10}$ & $10^{-12}$ \\ 
			RP3 &  1.0      &  0.0       & -0.2 & 1.0     & 0.5        &  0.0        & +0.2 & 0.5      & 1 & 1 & $10^{20}$ & $10^{20}$ \\ 
			RP4 &  1.0      &  0.0       & -0.2 & 1.0     & 0.5        &  0.0        & +0.2 & 0.5      & 1 & 1 & $10^{-10}$ & $10^{-12}$ \\ 
			\hline
		\end{tabular}
	\end{center} 
	\label{tab.ic.ideal}
\end{table}

\subsection{Simple shear layer}

\begin{figure}[!htbp]
	\begin{center}
		\begin{tabular}{cc} 
			\includegraphics[width=0.38\textwidth]{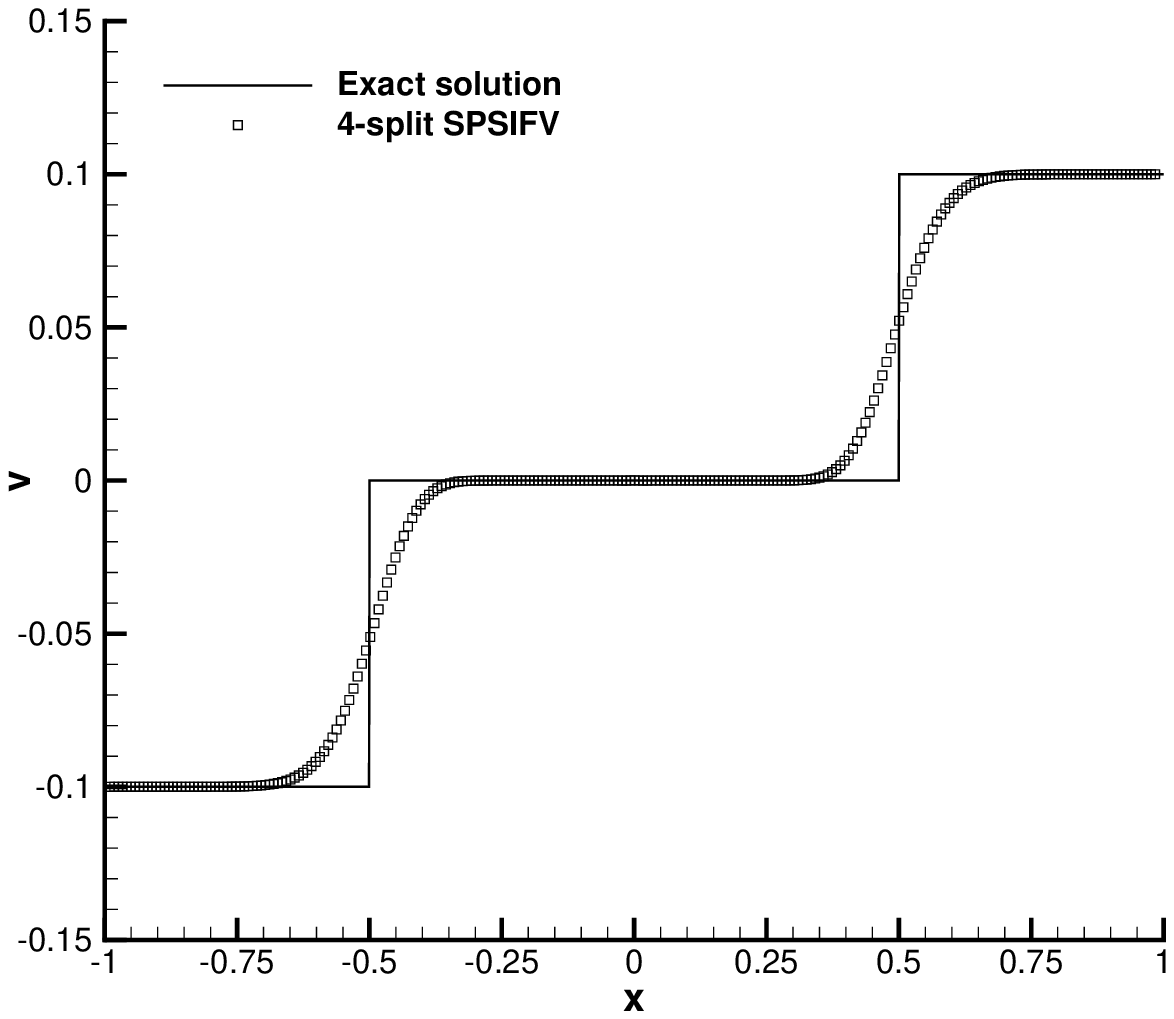}   & 
			\includegraphics[width=0.38\textwidth]{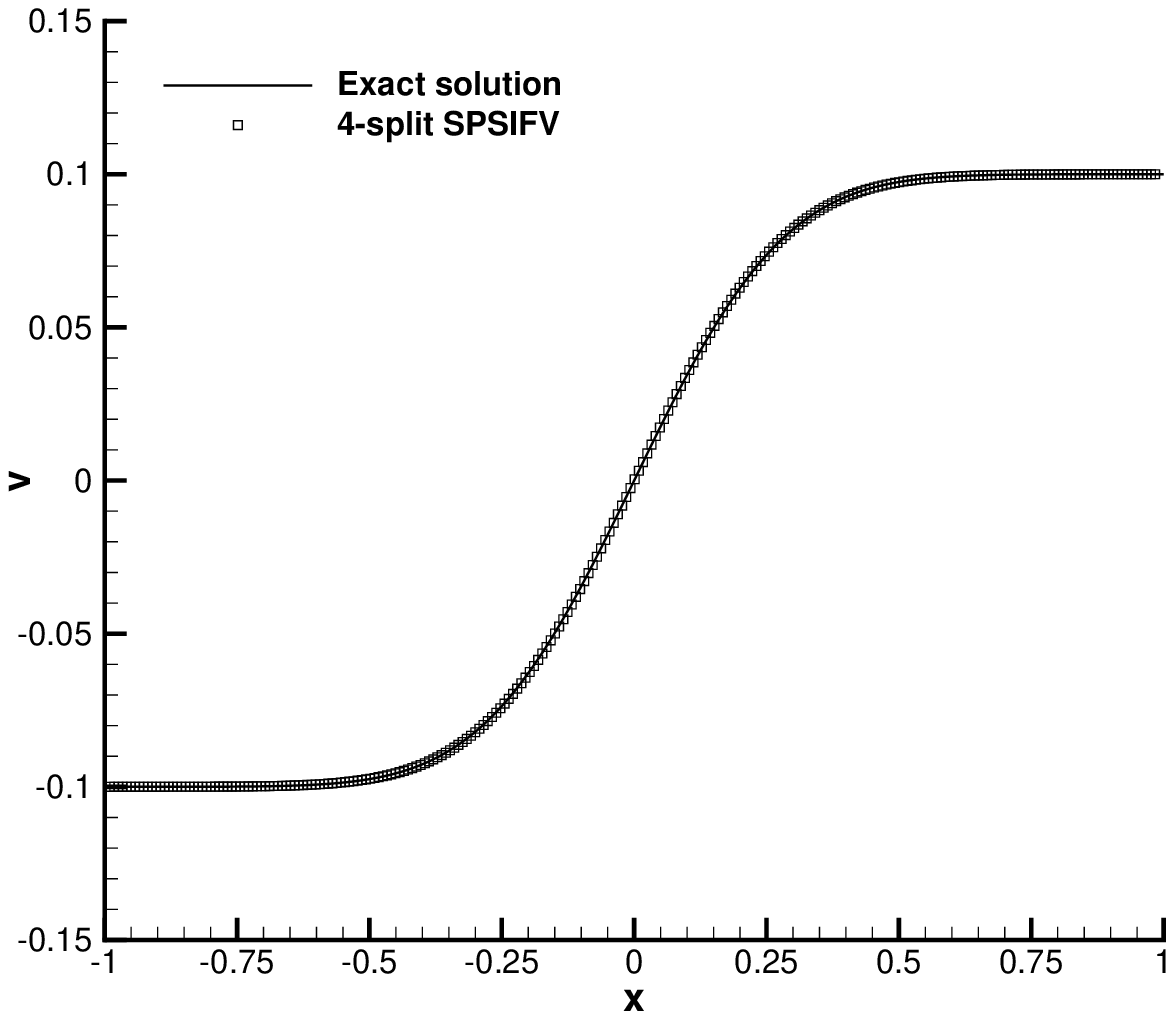}  \\  
			\includegraphics[width=0.38\textwidth]{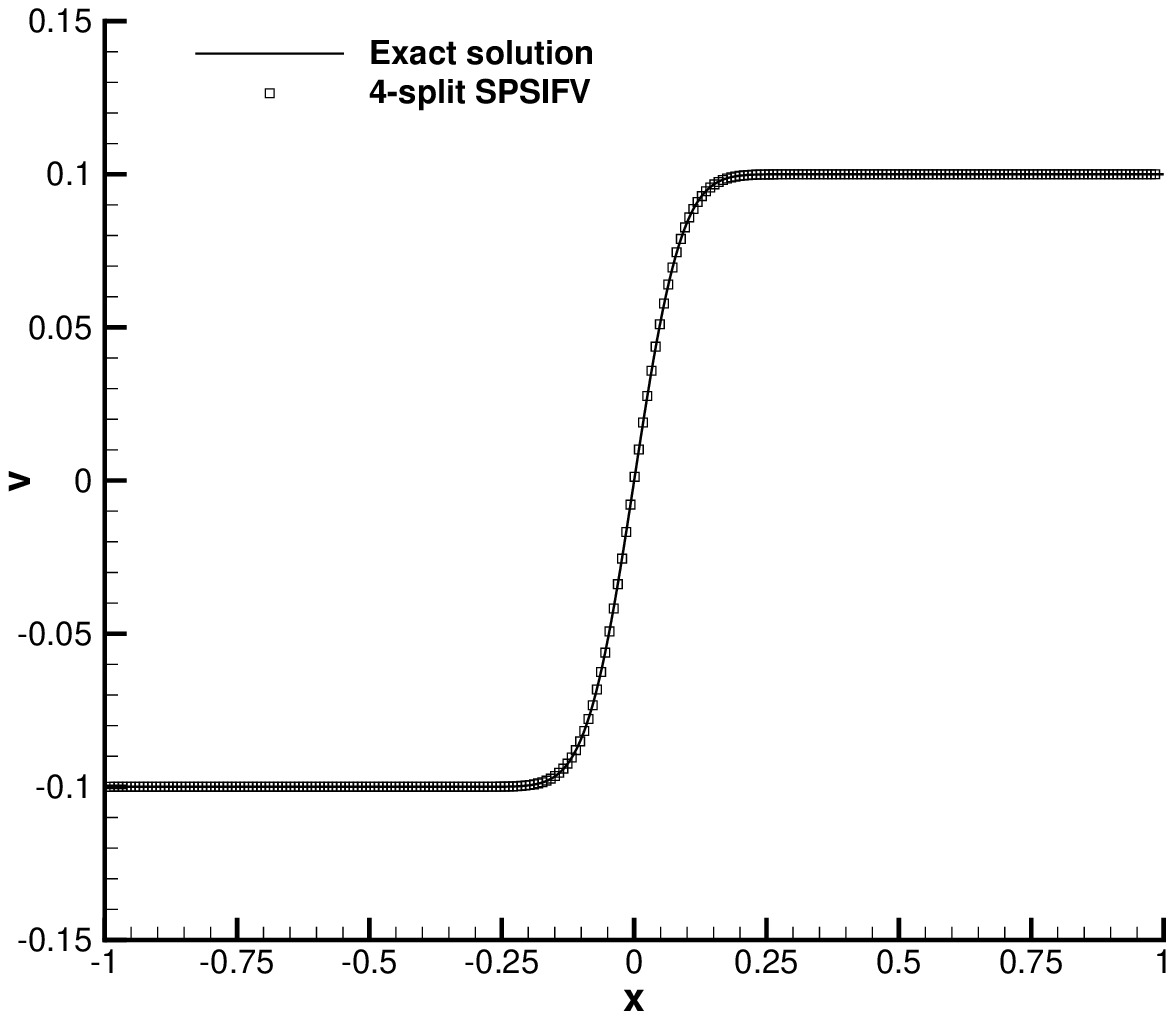}   & 
			\includegraphics[width=0.38\textwidth]{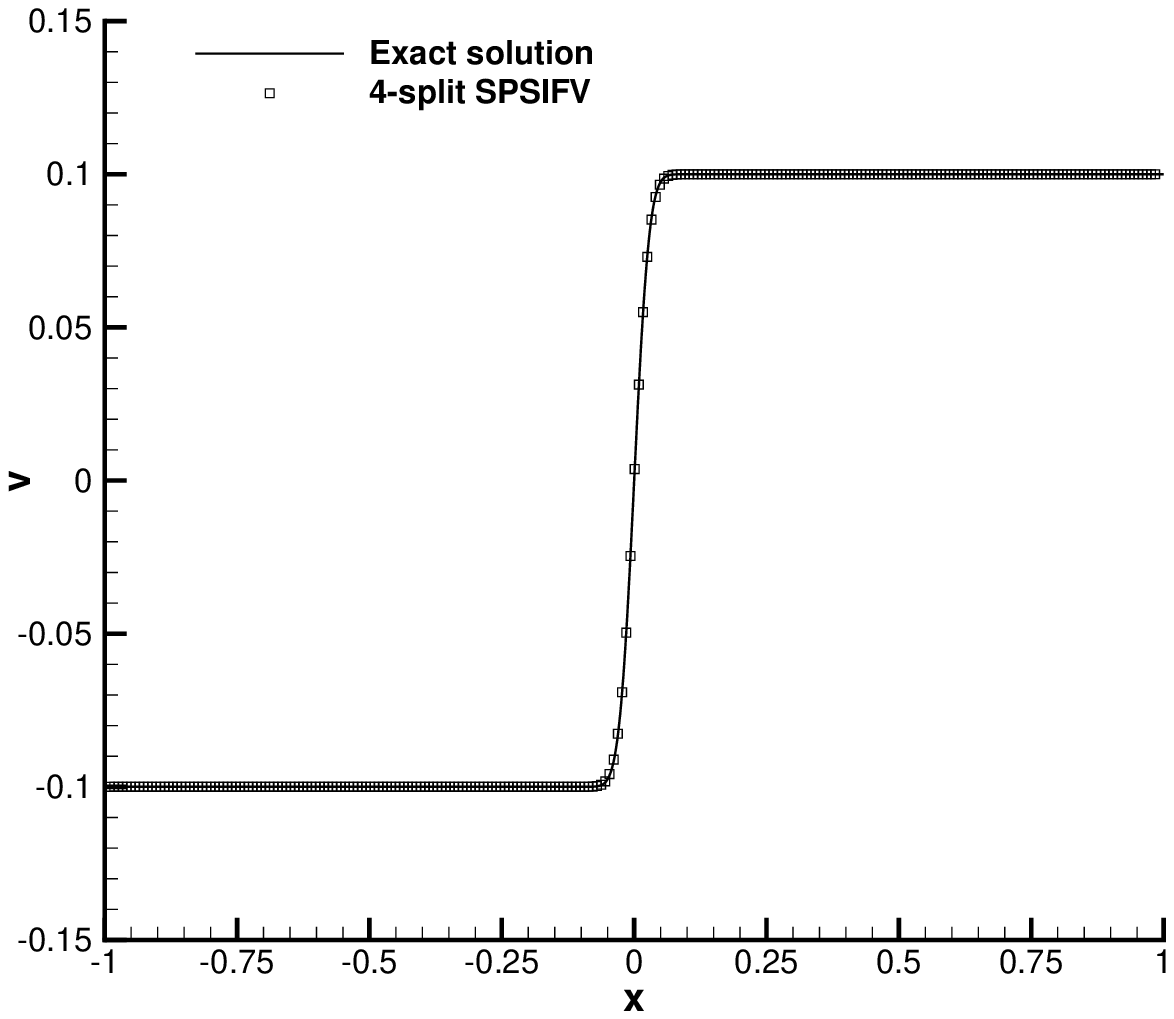}  
		\end{tabular} 
		\vspace*{-4mm}
		\caption{Numerical solution obtained with the new 4-split semi-implicit FV scheme at the corresponding final times for the GPR model applied to a simple shear flow problem. Results for the solid (top left) and for fluids with different viscosities: $\mu=10^{-1}$ (top right), $\mu=10^{-2}$ (bottom left) and $\mu=10^{-3}$ (bottom right).  } 
		\label{fig.shear}
	\end{center}
\end{figure}

In this section we test our novel 4-split semi-implicit finite volume scheme on a simple shear flow in fluids and solids. The computational domain in two space dimensions is given by $\Omega = [-1,+1] \times [-1,+1]$ with initial data that read:  
$\rho=1$, $v_1=v_3=0$, $p=10^5$, $\A=\mathbf{I}$, $\mathbf{J}=\mathbf{0}$, while the velocity 
component $v_2$ is $v_2 = -v_0$ for $x<0$
and $v_2=+v_0$ for $x\geq 0$, with $v_0=0.1$. The other parameters are $\gamma=1.4$, $c_v=1$, $\rho_0=1$, $c_s=10^3$, $c_h=10^2$ and $\tau_2=10^{-12}$. We impose periodic boundary conditions in the $y$ direction and Neumann boundaries in the $x$ direction. The calculations are carried out on a computational mesh of $1000 \times 10$ control volumes up to the final time of $t=0.25$ for the fluid case and up to $t=5 \cdot 10^{-4}$ for the solid case. The corresponding acoustic Mach number of this test case is of the order $M_a \approx 10^{-3}$, while the shear Mach number is $M_s = 10^{-4}$. 
In the Navier-Stokes limit the reference solution is given by the exact solution of the incompressible Navier-Stokes equations for the first problem of Stokes, see e.g. \cite{GPRmodel,SIGPR}.  
For the solid limit ($\tau_1 \to \infty$), this initial condition leads to two shear waves traveling to the left and right, respectively, with speed $c_s$. The corresponding Riemann problem for the linearized equations can be solved exactly. The time step in all cases is set to $\Delta t = 5 \cdot 10^{-3}$ in the fluid case and $\Delta t = 10^{-5}$ in the solid case.   
A comparison of the numerical results with the reference solution for different values of $\mu$ is depicted in Fig. \ref{fig.shear}. An excellent agreement between numerical solution and reference solution can be observed for all cases, despite the very low Mach numbers. 

\subsection{Lid-driven cavity}

The lid-driven cavity \cite{Ghia1982} is a classical benchmark for incompressible Navier-Stokes solvers. However, it can also be used to validate compressible flow solvers in the low Mach number regime, see e.g. \cite{TavelliDumbser2017,DumbserCasulli2016} and it has also been successfully used as benchmark problem for numerical methods applied to the GPR model in \cite{GPRmodel,SIGPR,HTCGPR,HTCAbgrall}.
The computational domain used in this paper is $\Omega = [0,1] \times [0,1]$. The initial condition is $\rho=1$, $\mathbf{v}=0$, $p=10^8$,  $\A=\mathbf{I}$ and $\mathbf{J}=0$. Furthermore we set $\gamma=1.4$, $c_v = 10^5$, $c_s = 10^3$, $\rho_0=1$, $\tau_2 = 10^{-14}$ and $c_h=100$. The viscosity coefficient is set to $\mu=10^{-2}$, hence the Reynolds number of the flow is $Re=100$. The flow inside the cavity is generated by the moving lid at the upper boundary, whose velocity is set to $\mathbf{v}=(1,0,0)$. On all other boundaries, a no-slip wall boundary condition with $\mathbf{v}=0$ is imposed. With the chosen initial and boundary conditions, the acoustic Mach number of this test problem is $M_a=2.7 \cdot 10^{-4}$ and the shear Mach number is $M_s = 10^{-3}$ with respect to the lid velocity.

\begin{figure}[!htbp]
	\begin{center}
			\includegraphics[trim=10 10 10 10,clip,width=0.47\textwidth]{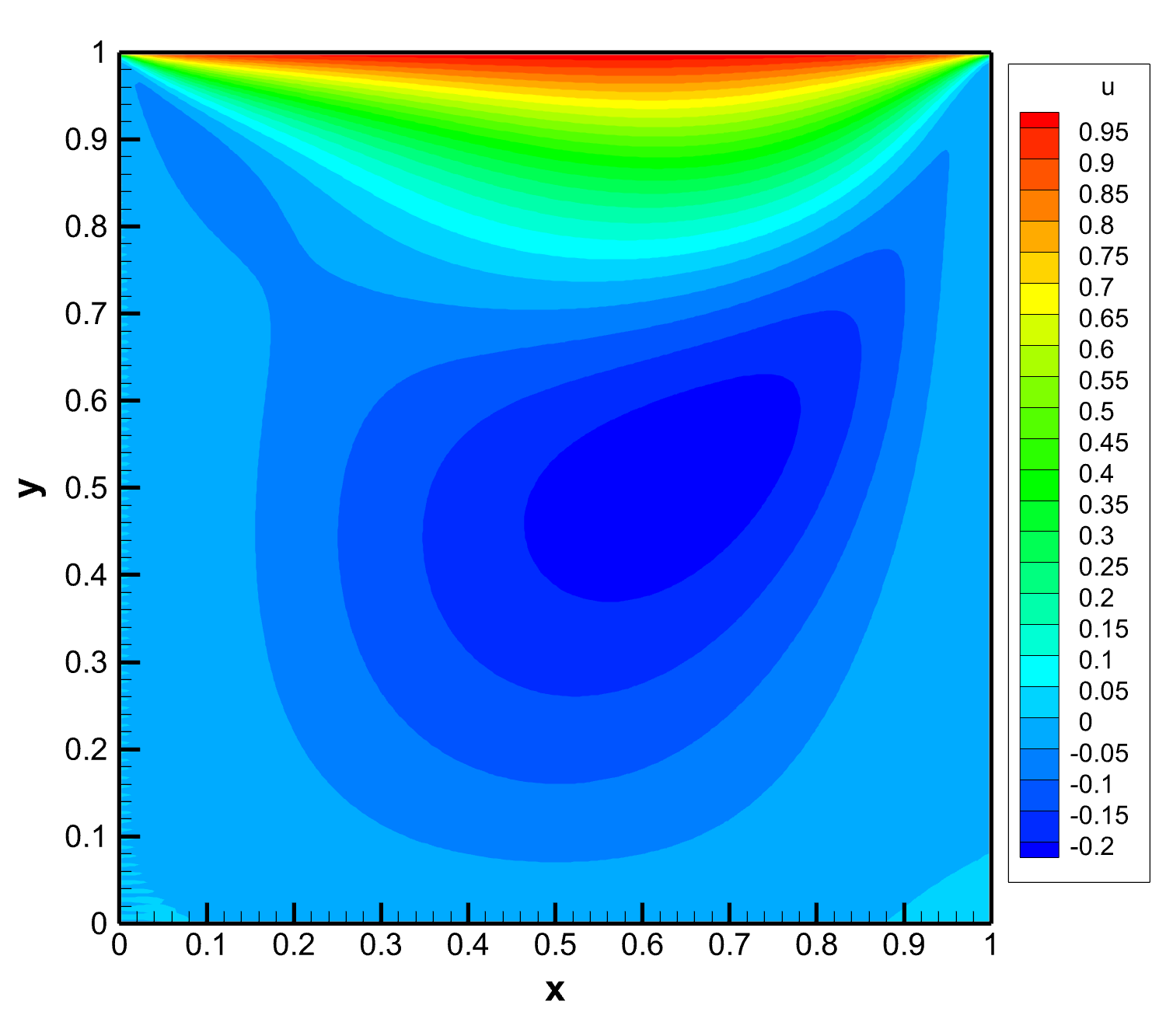} \\ \vspace{2cm}
			\includegraphics[width=0.47\textwidth]{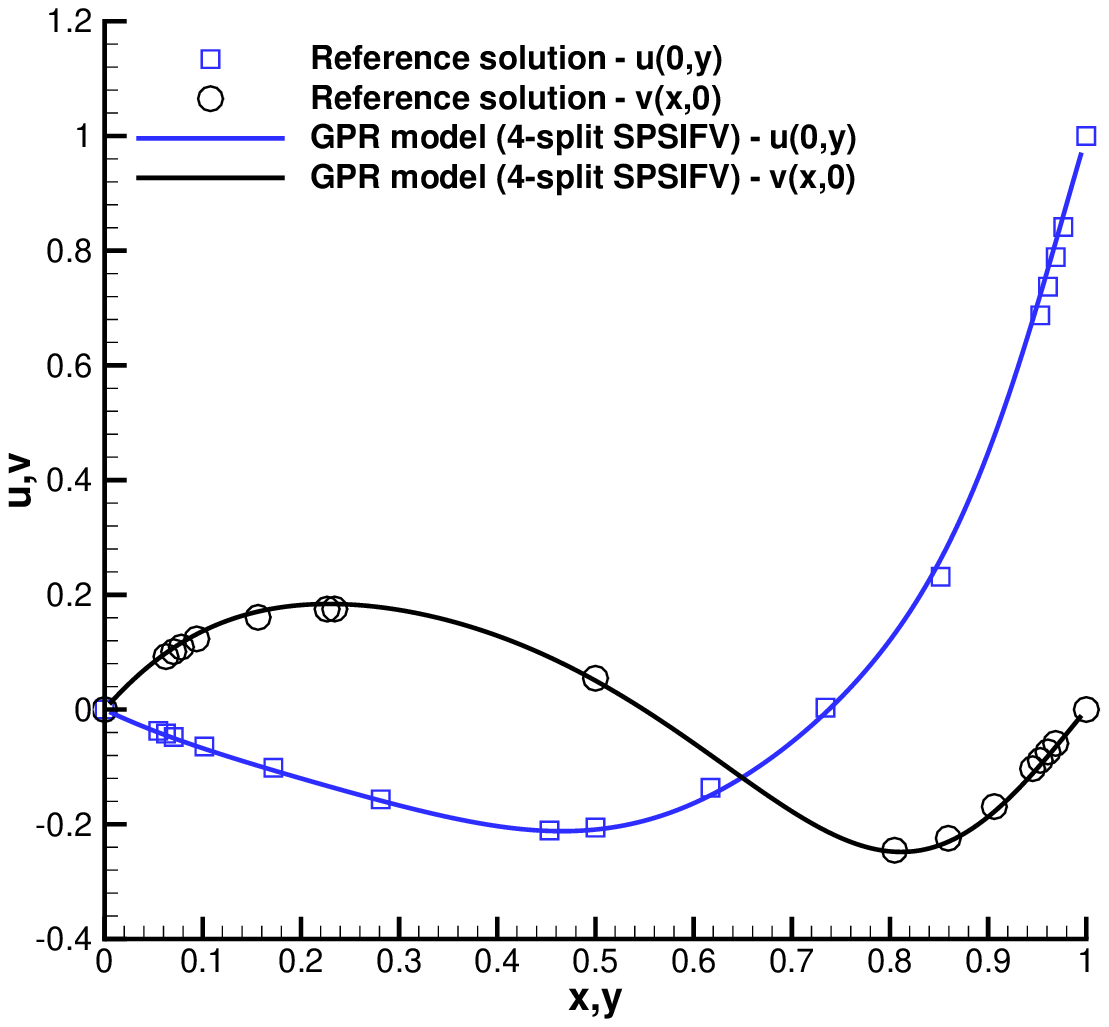}\vspace{-2cm}
		\caption{Lid-driven cavity problem at $t=10$ for Mach number $M=2.7 \cdot 10^{-4}$, Reynolds number $Re=100$ and $M_s=10^{-3}$. Color contours of the velocity component $u$ (left) and
			comparison of the velocity components $u$ and $v$ on 1D cuts along the $x$ and $y$ axis with the reference solution of Ghia \textit{et al.} \cite{Ghia1982} (right). }
		\label{fig.cavity}
	\end{center}
\end{figure}

The new structure-preserving 4-split finite volume scheme is run until a final time of $t=10$ using a computational grid composed of $200 \times 200$ elements. A direct comparison of our computational results with the reference solution of Ghia \textit{et al.} \cite{Ghia1982} obtained for the incompressible Navier-Stokes equations is depicted in Fig. \ref{fig.cavity}. One can observe a good agreement between the  numerical solution obtained with the new structure-preserving 4-split scheme for the GPR model and the incompressible Navier-Stokes reference solution.
Compared to \cite{GPRmodel} and \cite{SIGPR} we emphasize that our time step restriction is only based on the fluid velocity and not on the adiabatic sound speed, the shear sound speed $c_s$ or the heat wave speed $c_T$.

\textcolor{black}{Compared to a classical explicit FV scheme our timestep is 10,000 times larger, and compared to the structure-preserving semi-implicit scheme introduced in \cite{SIGPR}, in which shear and heat waves were still treated explicitly, the timestep of our new scheme is 1000 times larger. }

\subsection{Solid rotor}

\begin{figure}[!htbp]
	\begin{center}
		\begin{tabular}{ccc}
			\includegraphics[trim=10 10 10 10,clip,width=0.3\textwidth]{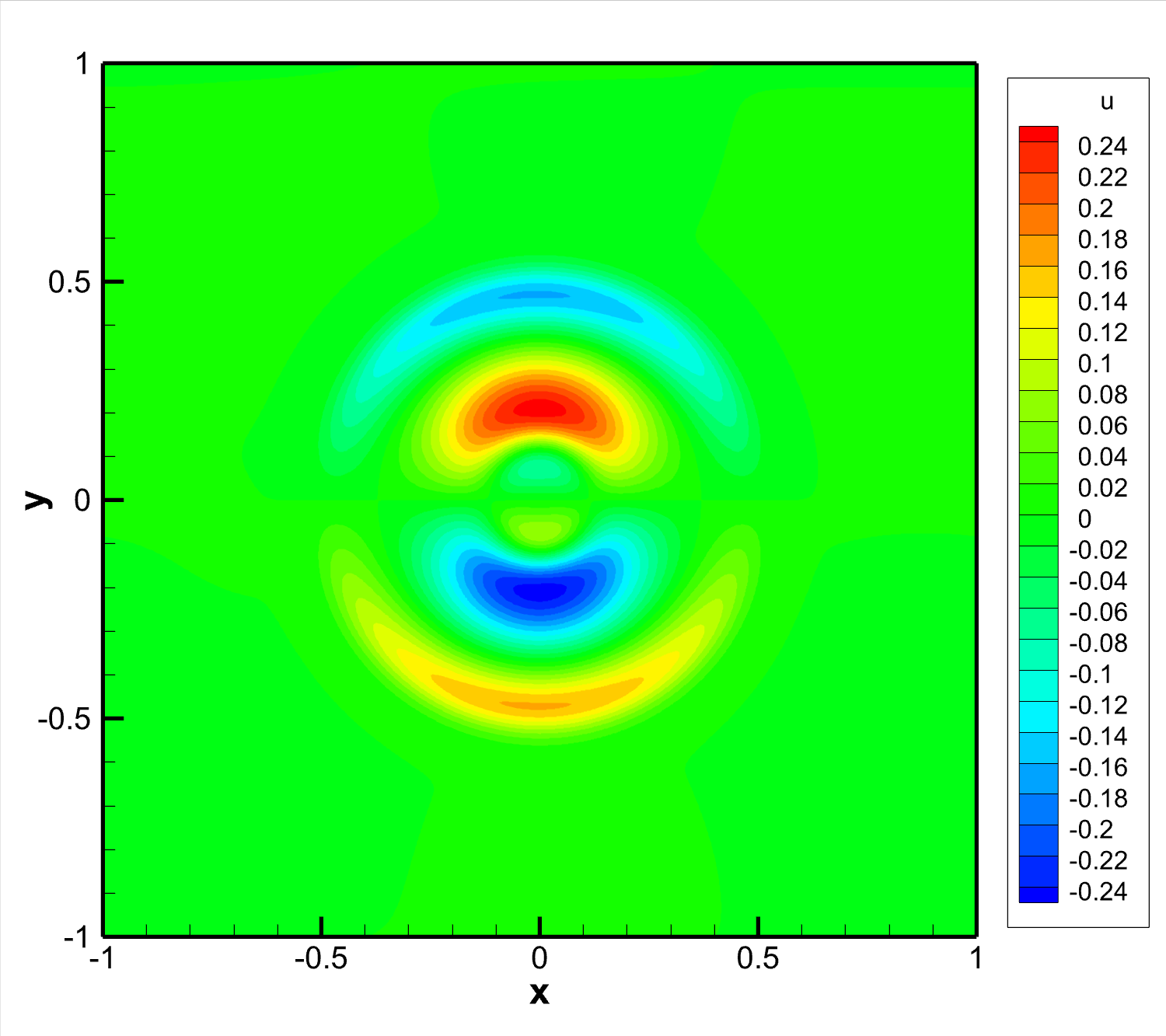}  &
			\includegraphics[trim=10 10 10 10,clip,width=0.3\textwidth]{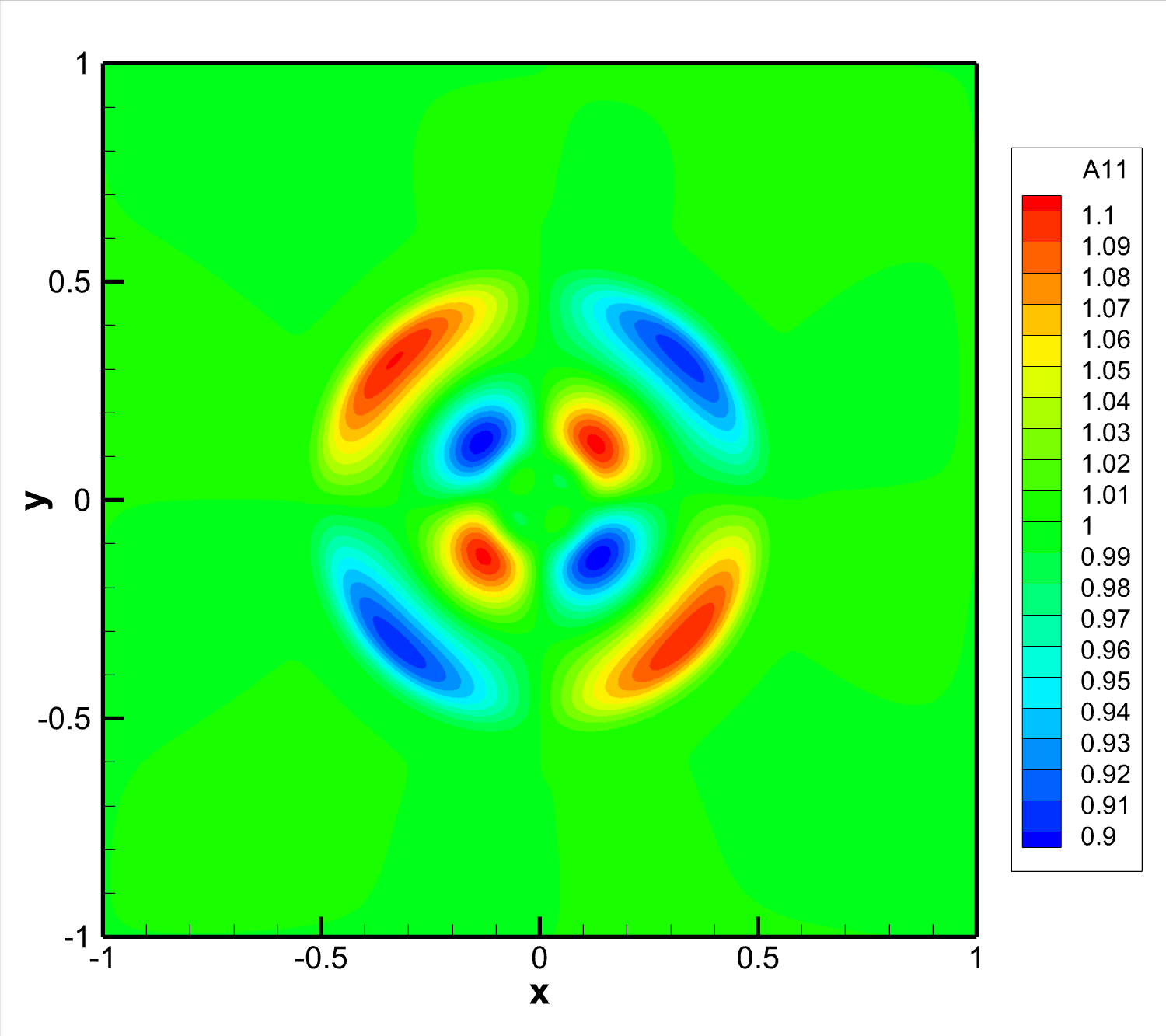}  &
			\includegraphics[trim=10 10 10 10,clip,width=0.3\textwidth]{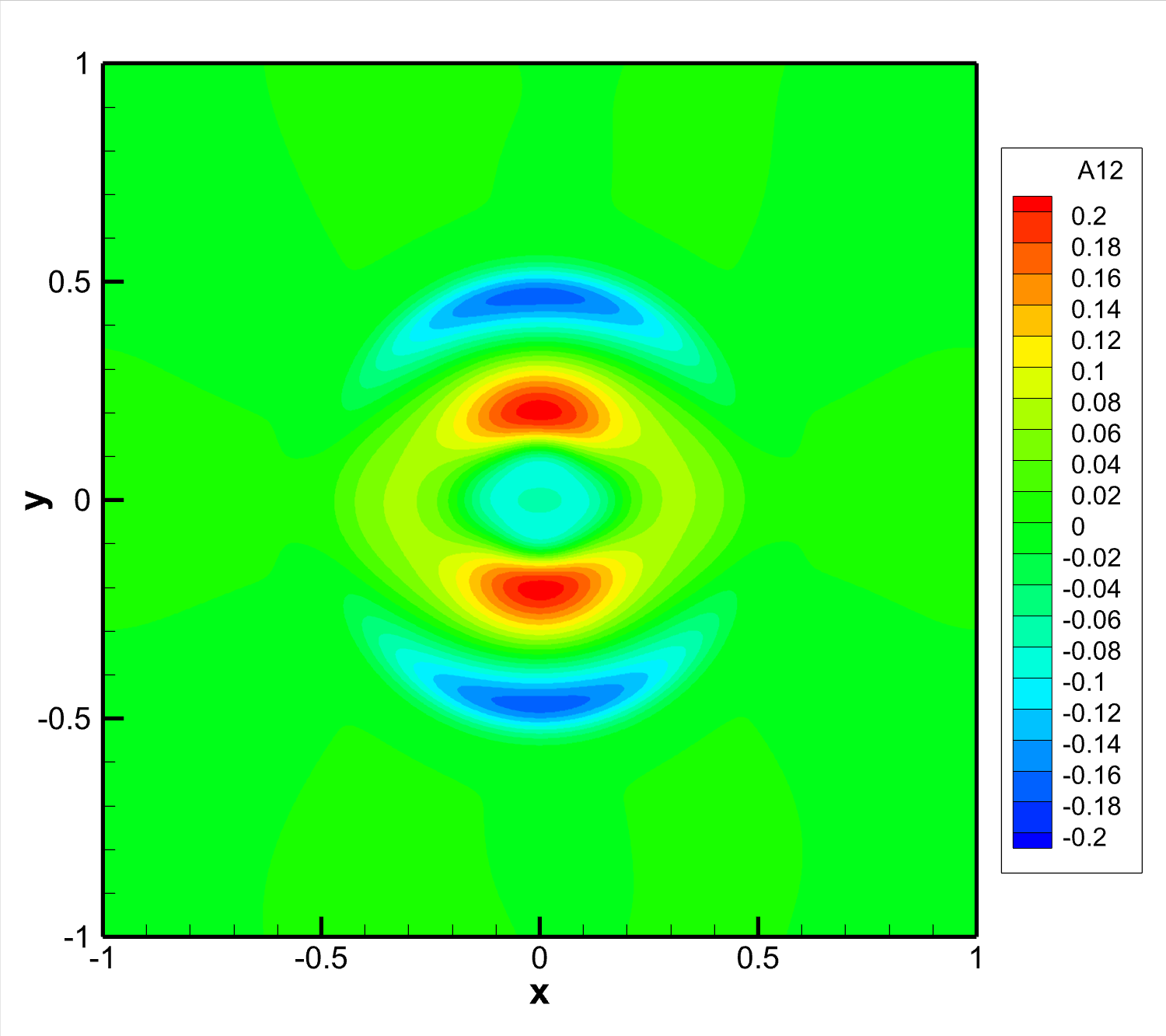}   
		\end{tabular}
		\caption{Solid rotor problem at time $t=0.3$. Contour colors of the horizontal velocity $u$ (left), component $A_{11}$ (center) and component $A_{12}$ (right) of the inverse deformation gradient. }
		\label{fig.solidrotor}
	\end{center}
\end{figure}

\begin{figure}[!htbp]
	\begin{center}
		\begin{tabular}{cc}
			\includegraphics[width=0.45\textwidth]{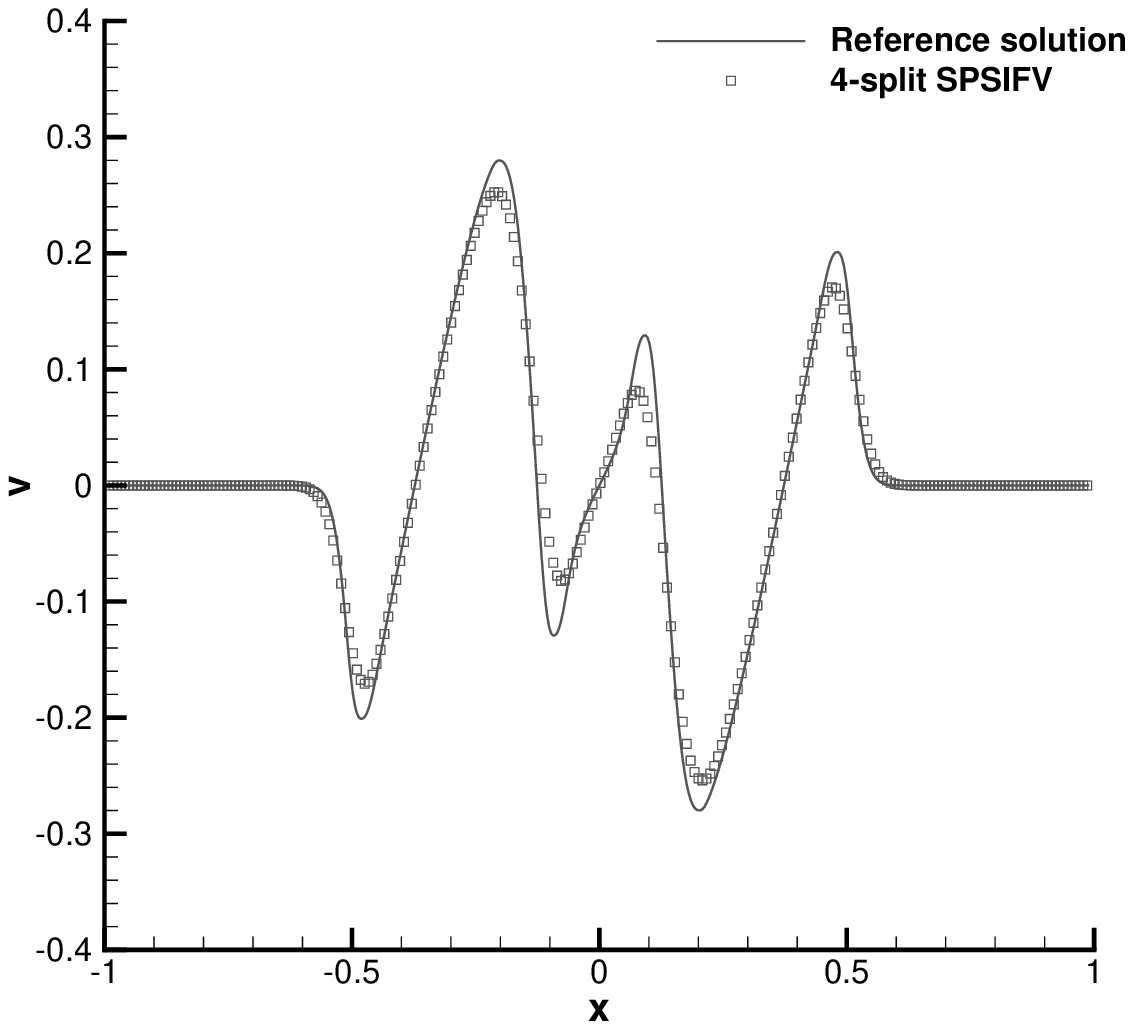}  &
			\includegraphics[width=0.45\textwidth]{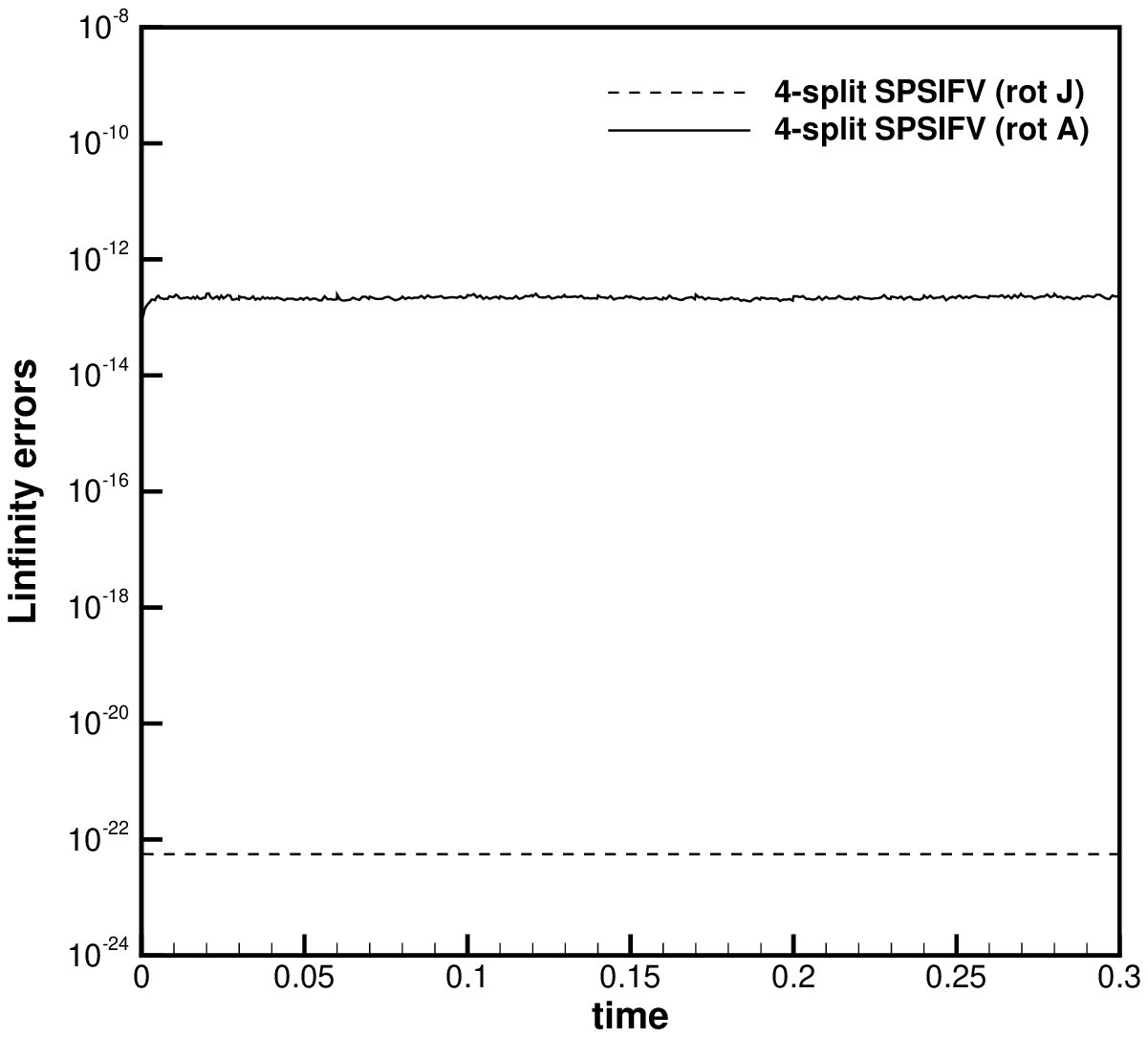} 
		\end{tabular}
		\caption{Solid rotor problem. 1D cut through the solution for velocity component $v$ at $y=0$ at time $t=0.3$ and comparison with the numerical reference solution (left). Time series of the curl errors of $\mathbf{A}$ and $\mathbf{J}$ in $L^{\infty}$ norm (right). }
		\label{fig.solidrotor.cutcurl}
	\end{center}
\end{figure}

\begin{figure}[!htbp]
	\begin{center}
		\begin{tabular}{cc} 
			\includegraphics[width=0.4\textwidth]{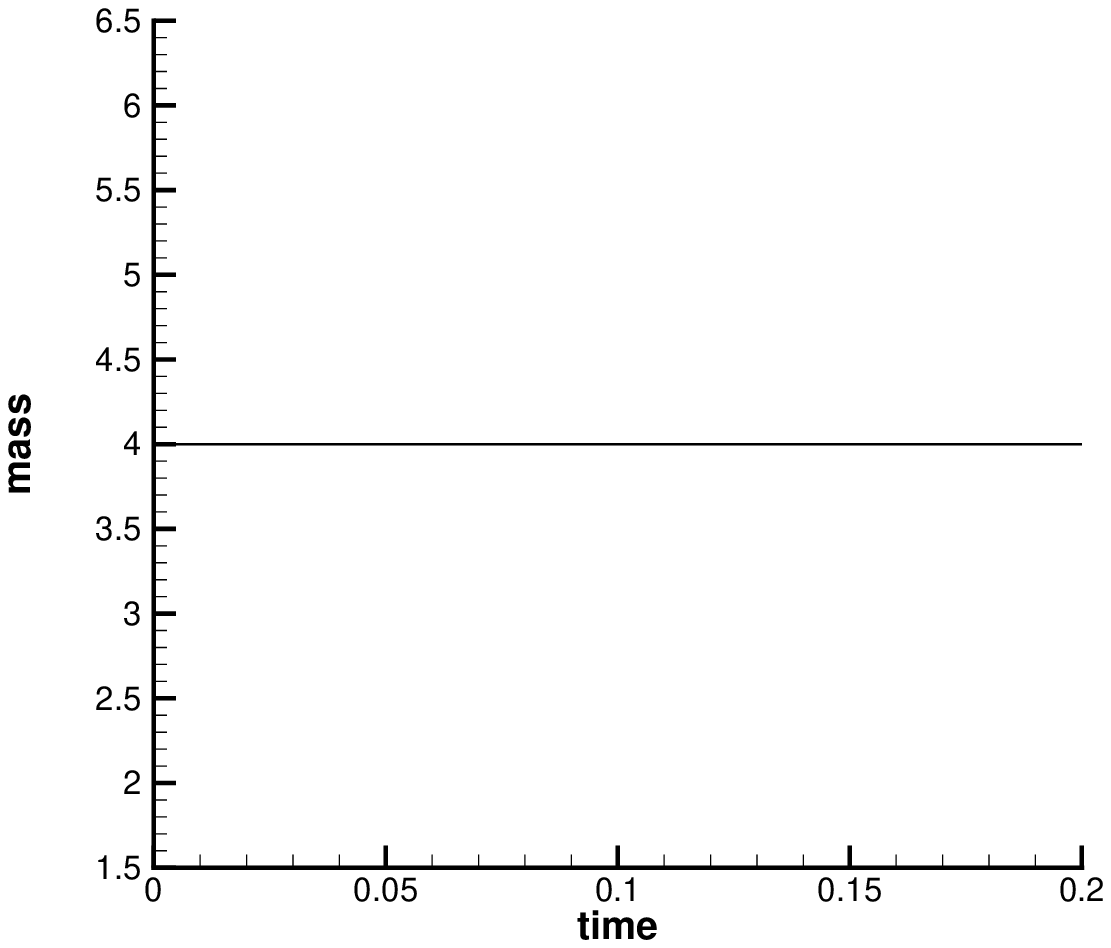}  &  
			\includegraphics[width=0.4\textwidth]{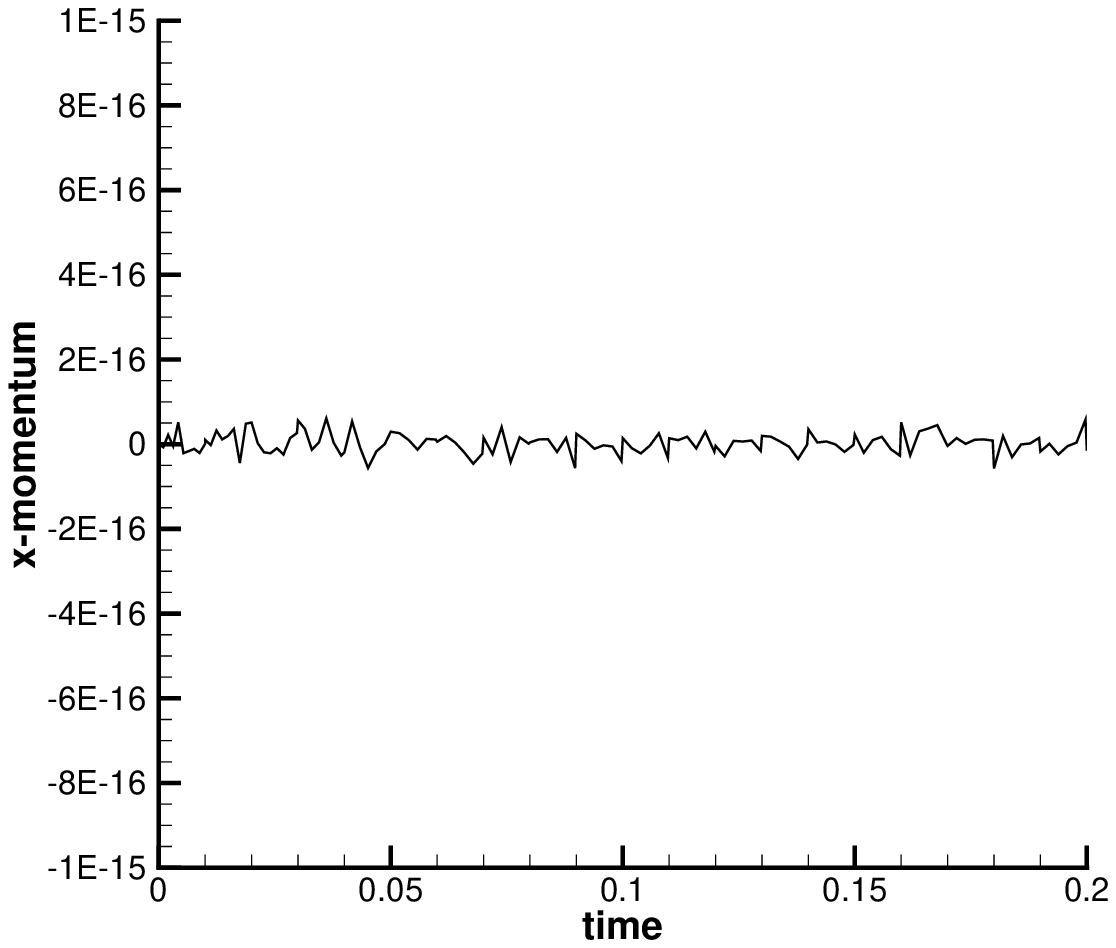}  \\ 
			\includegraphics[width=0.4\textwidth]{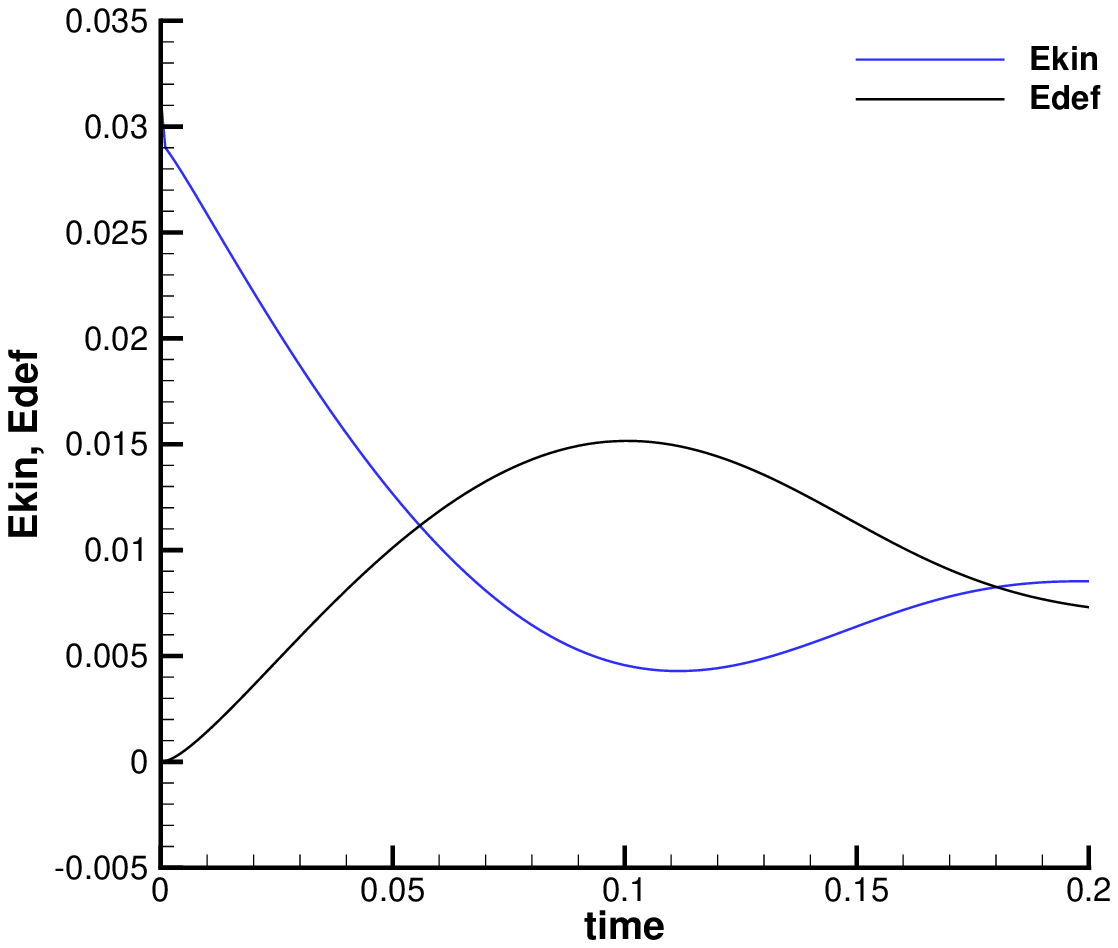}  & 
			\includegraphics[width=0.4\textwidth]{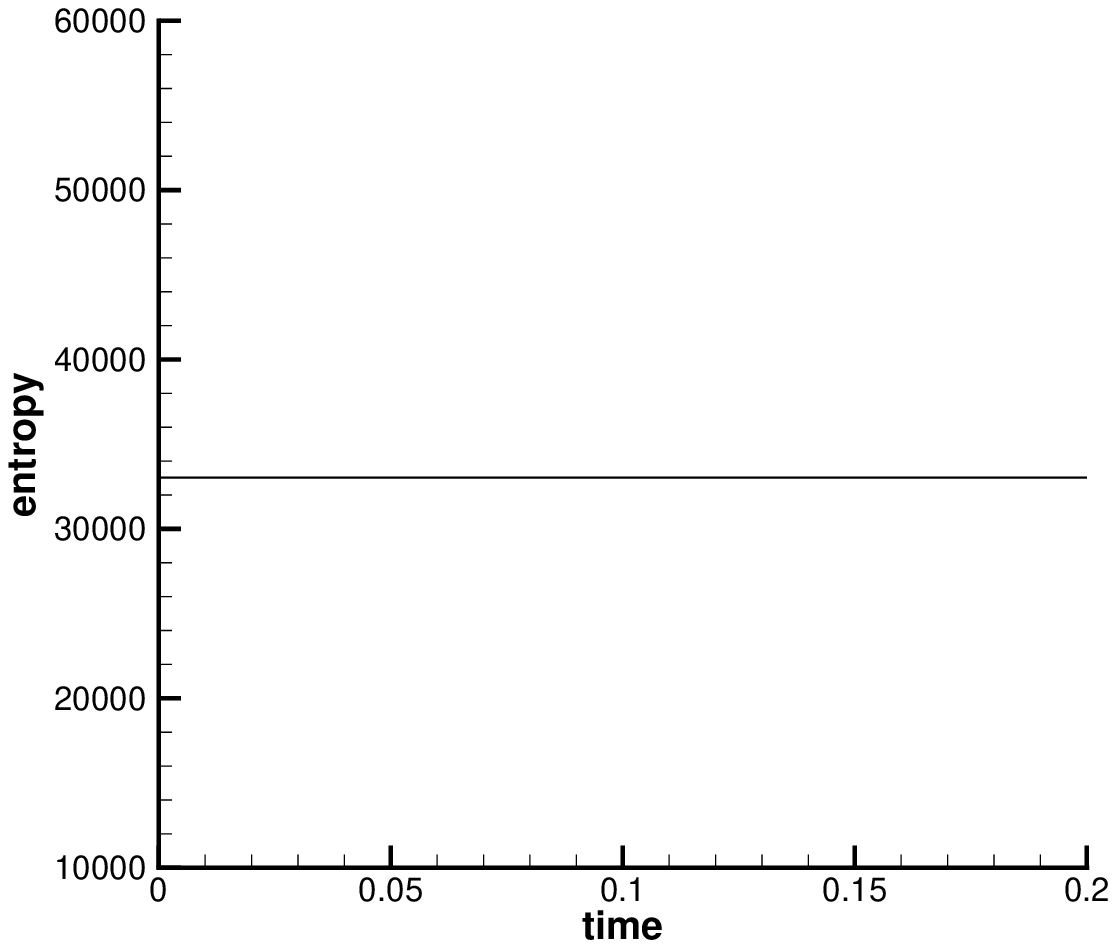}    
		\end{tabular} 
		\caption{\textcolor{black}{Solid rotor problem. Time evolution of mass, $x$-momentum, kinetic and deformation energy and entropy.}  } 
		\label{fig.solidrotor.energies}
	\end{center}
\end{figure}

The aim of this test case is to numerically verify the discrete curl-free property of our new 4-split semi-implicit finite volume scheme. For the setup of the test case we follow \cite{SIGPR}, solving \eqref{eqn.GPR} with $\tau_1  = 10^{20}$, but here we set $\tau_2=10^{-14}$, i.e. the strain relaxation source term is absent, while heat conduction is considered in the stiff Fourier limit. The computational domain 
is given by $\Omega = [-1,+1]^2$. The initial condition is
$\rho = 1$, $p = 10^5$, $\mathbf{A} = \mathbf{I}$ and $\mathbf{J}=0$ and the initial velocity field is $u = -y/R$, $v = +x/R$ and $w=0$ within the circular region $r \leq R$, where $r = \left\| \mathbf{x} \right\|_2$ and $R=0.2$, while $\mathbf{v}=0$ for $r > R$. The remaining parameters are $\gamma = 1.4$, $c_v=1004/\gamma$, $c_s = 1.0$ and
$c_h = 100$. With this setting we are in the low Mach regime for the acoustic Mach number $M_a$ and for the heat Mach number $M_h$. The test case is solved with the new structure-preserving semi-implicit 4-split finite volume scheme until a final time of $t=0.3$. The computational mesh uses $512 \times 512$ cells. In order to obtain a numerical reference solution, the same problem is solved on the same grid with the second order version of the 2-split structure-preserving scheme introduced in \cite{SIGPR}, but in the absence of heat conduction, i.e. setting $c_h=0$, since the 2-split method presented in \cite{SIGPR} is \textit{not} able to deal with low heat Mach numbers $M_h$. The results obtained with the new 4-split scheme are depicted in Fig. \ref{fig.solidrotor}, where the contour colors of the velocity component $u$ and the distortion field components $A_{11}$ and $A_{12}$ are shown. A comparison with the reference solution is provided in the 1D cut shown in the left panel of Fig. \ref{fig.solidrotor.cutcurl}, while in the right panel of Fig. \ref{fig.solidrotor.cutcurl} we show the the time series of the curl errors of $\mathbf{A}$ and $\mathbf{J}$ in the $L^\infty$ norm. One can observe a good agreement with the reference solution and that the new structure-preserving semi-implicit 4-split scheme presented in this paper produces curl errors that remain at the order of machine precision, as expected. \textcolor{black}{A standard MUSCL-Hancock scheme produces instead larger curl errors that are of the order of the mesh spacing, see \cite{SIGPR} for comparison. In Figure \ref{fig.solidrotor.energies} we report the temporal evolution of mass, $x$-momentum, entropy and of the kinetic energy as well as of the energy stored in the deformation of the solid. For this smooth problem, entropy is preserved, as expected. }

\subsection{Circular explosion problem}

As last test case we solve two circular explosion problems (EP1 and EP2), one in the quasi inviscid fluid limit of the model, one in the solid limit. The computational domain is $\Omega=[-1;1]^2$ and the initial condition is given by
\begin{equation}
	\mathbf{Q}(x,y,0)=\begin{cases}
		\mathbf{Q}_{in} \quad \text{if} \quad r\leq R \\
		\mathbf{Q}_{out} \quad \text{if} \quad r>R, \\
	\end{cases}
\end{equation}
with $\mathbf{Q}_{in}$ and $\mathbf{Q}_{out}$ the inner and outer states, respectively, and $r=\sqrt{x^2+y^2}$ is the radial coordinate. The initial discontinuity is placed at a radius of $R=0.5$.

\begin{figure}[!htbp]
	\begin{center}
		\begin{tabular}{cc} 
			\includegraphics[trim= 10 10 10 10,clip,width=0.47\textwidth]{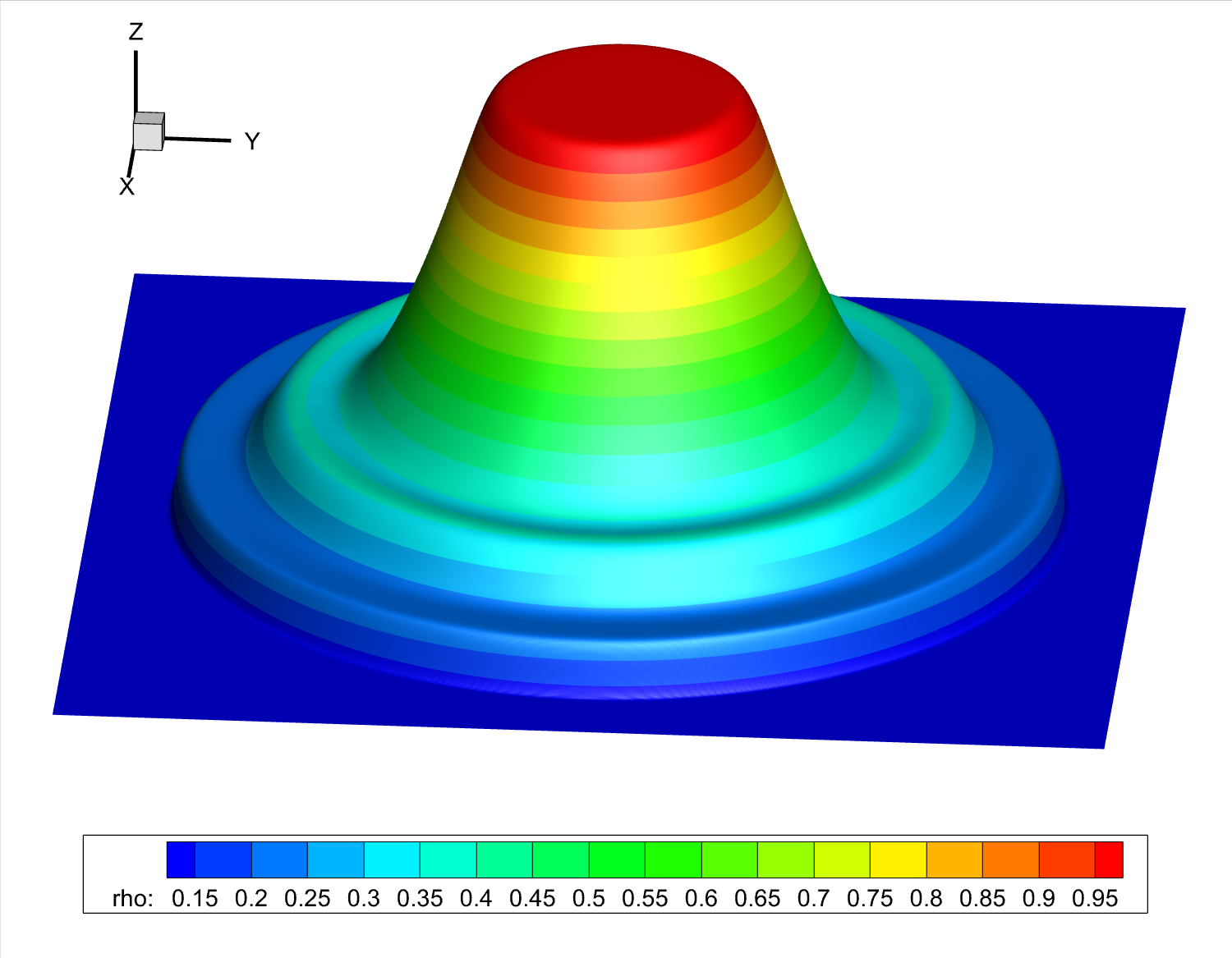}  &
			\includegraphics[trim= 10 10 10 10,clip,width=0.47\textwidth]{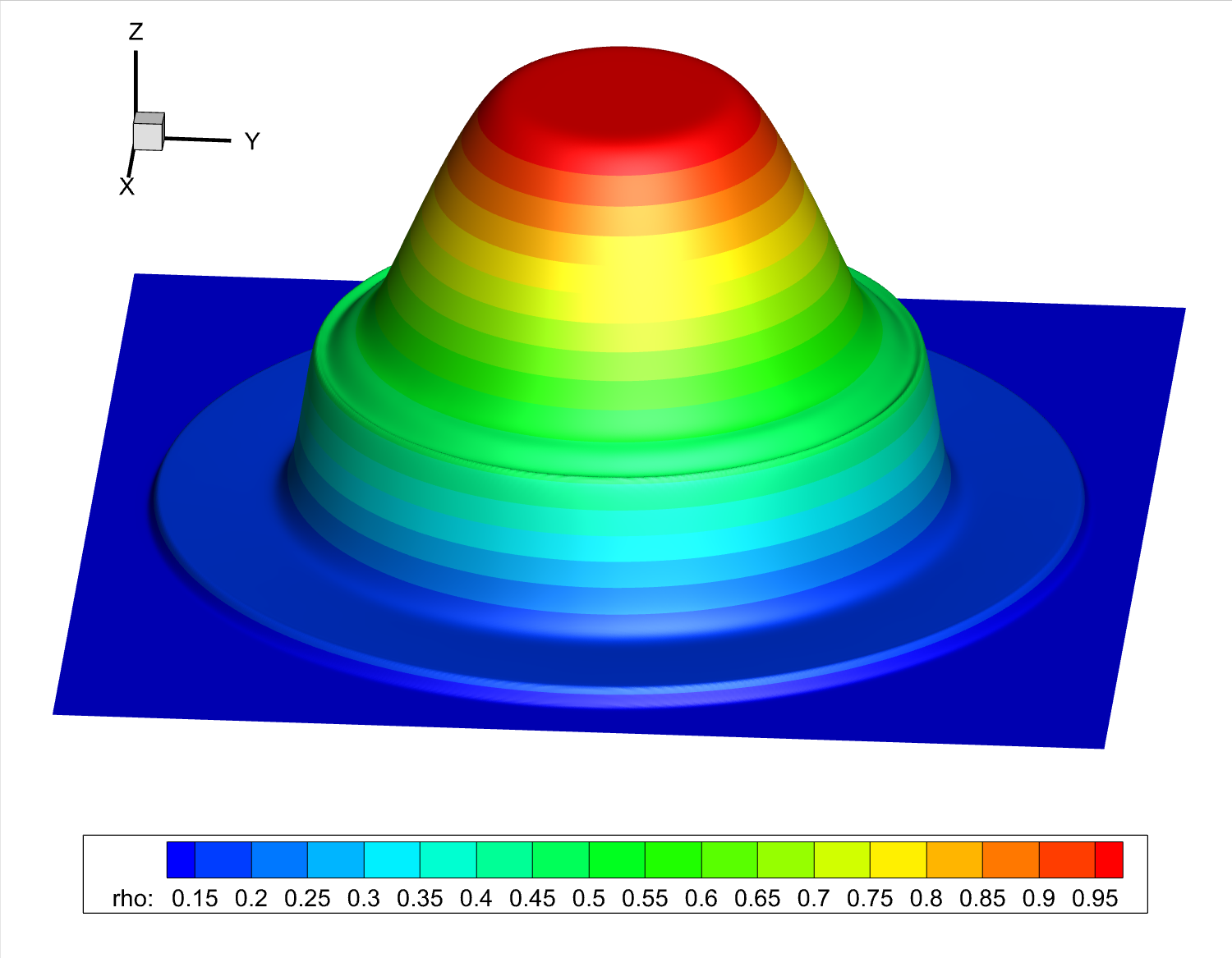}   
		\end{tabular} 
		\caption{2D explosion problems EP1 and EP2 at their final times $t=0.2$ and $t=0.15$, respectively. 3D density contour color plot for EP1 (left) and EP2 (right).  } 
		\label{fig.ep2dfluidsolid}
	\end{center}
\end{figure}

\begin{figure}[!htbp]
	\begin{center}
		\begin{tabular}{ccc} 
			\includegraphics[width=0.3\textwidth]{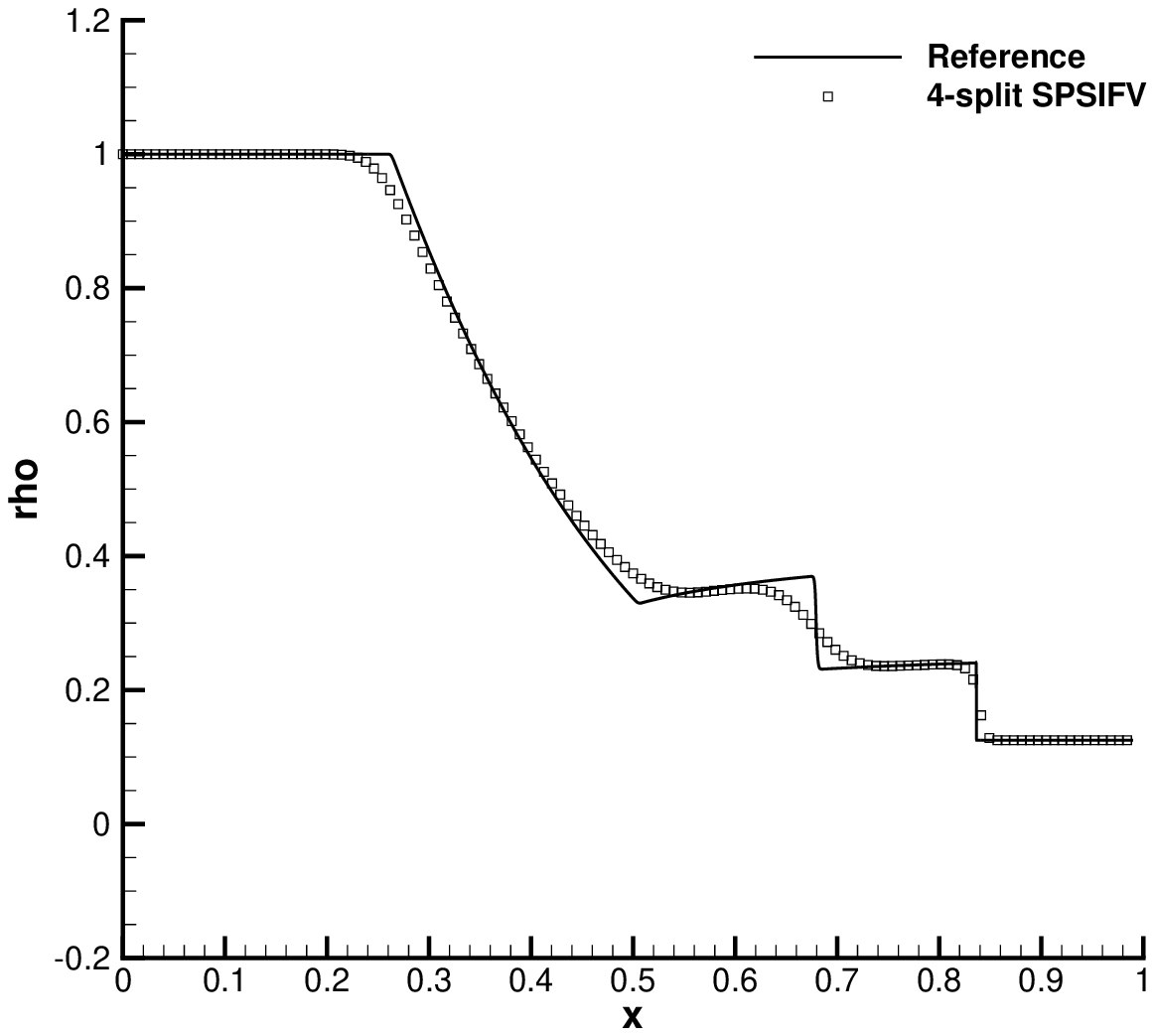}  &  
			\includegraphics[width=0.3\textwidth]{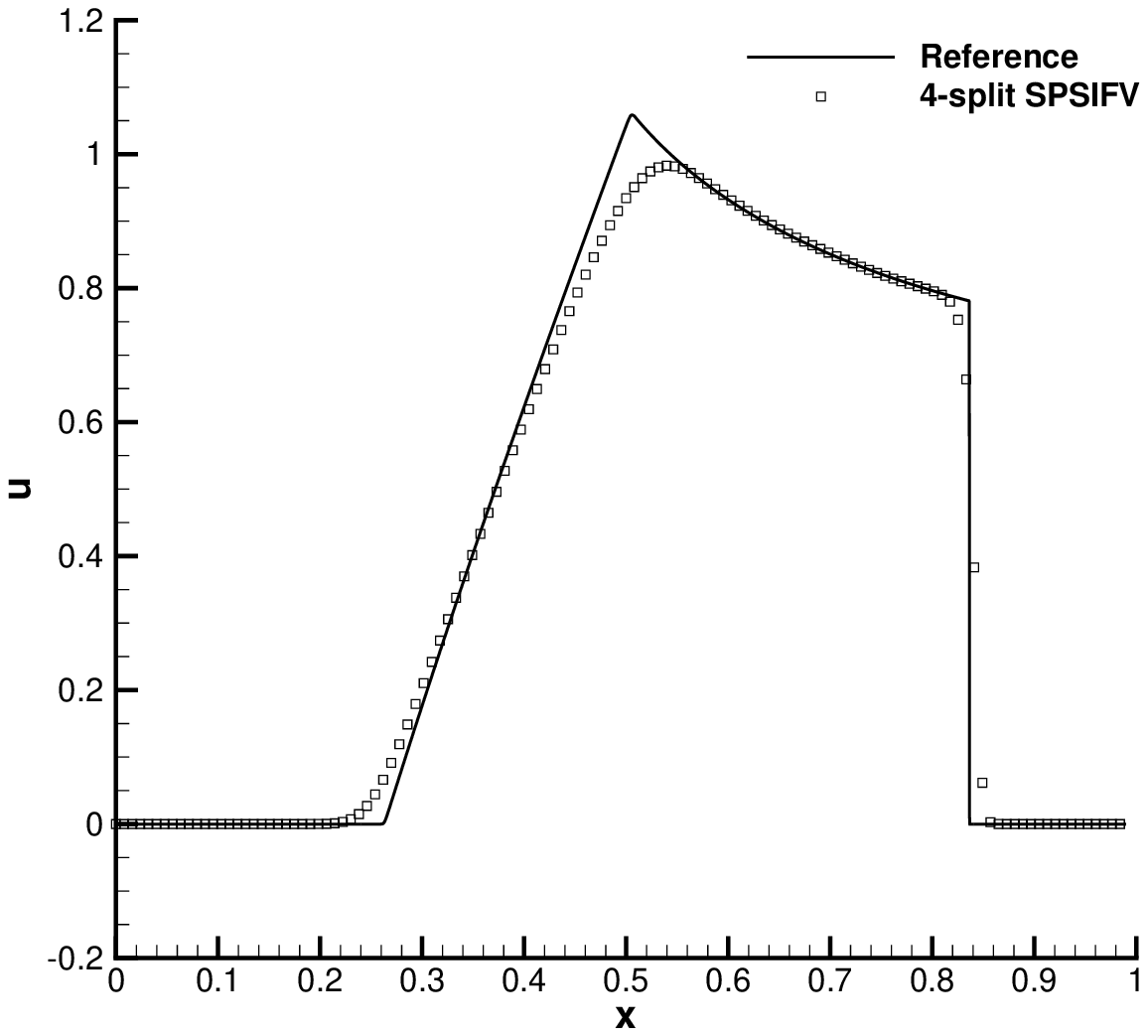}  & 
			\includegraphics[width=0.3\textwidth]{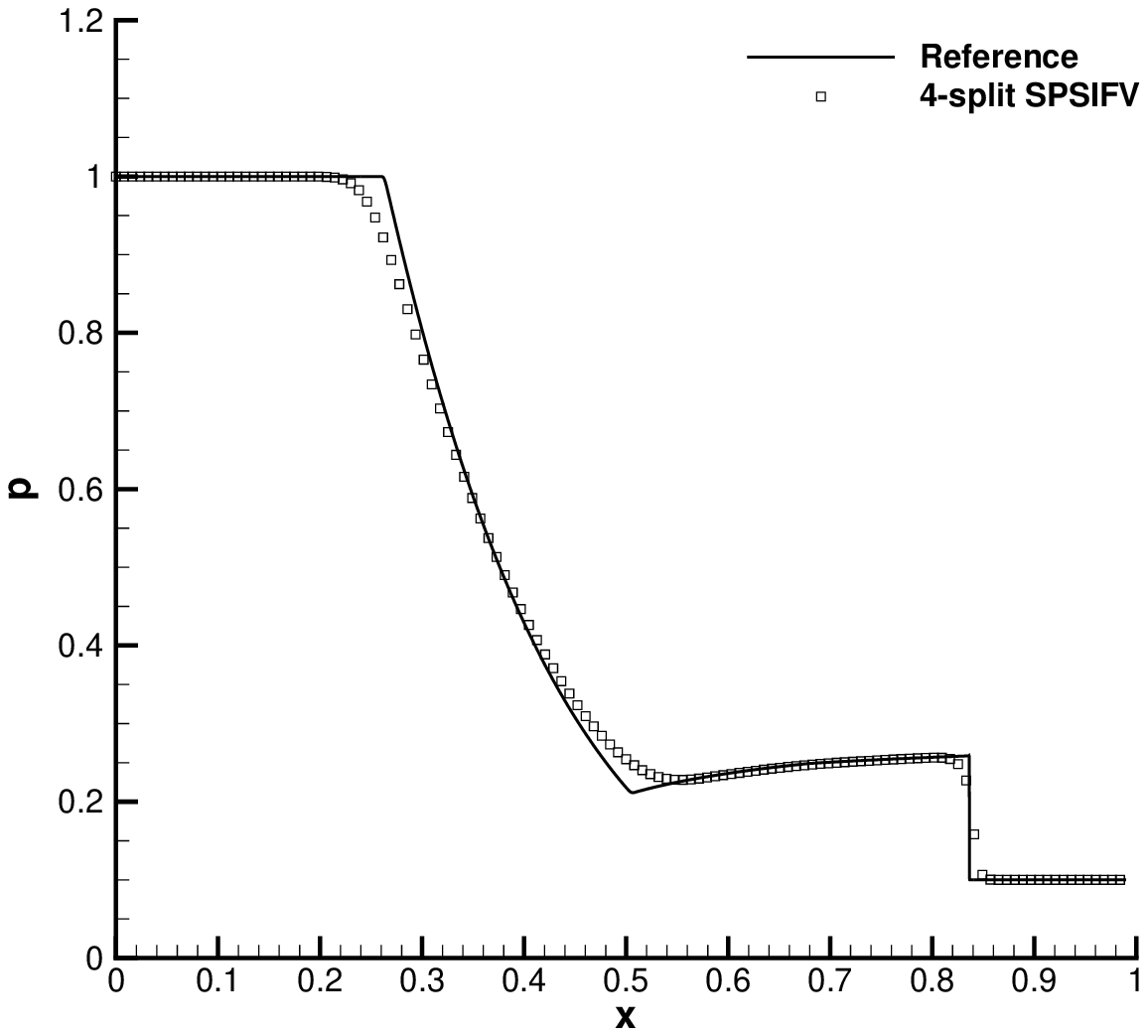}  \\ 
		\end{tabular} 
		\caption{2D explosion problem EP1 at time $t=0.2$. 1D cuts along the $x$ axis providing a comparison of the Euler reference solution (solid line) against the numerical solution of the GPR model obtained 
			with the new 4-split SIFV scheme (square symbols) in the stiff relaxation limit (fluid limit, $\tau_1=10^{-8}$, $\tau_2=10^{-10}$) for density (left), radial velocity (center) and pressure (right).   } 
		\label{fig.ep2dfluid}
	\end{center}
\end{figure}

\begin{figure}[!htbp]
	\begin{center}
		\begin{tabular}{cc} 
			\includegraphics[width=0.4\textwidth]{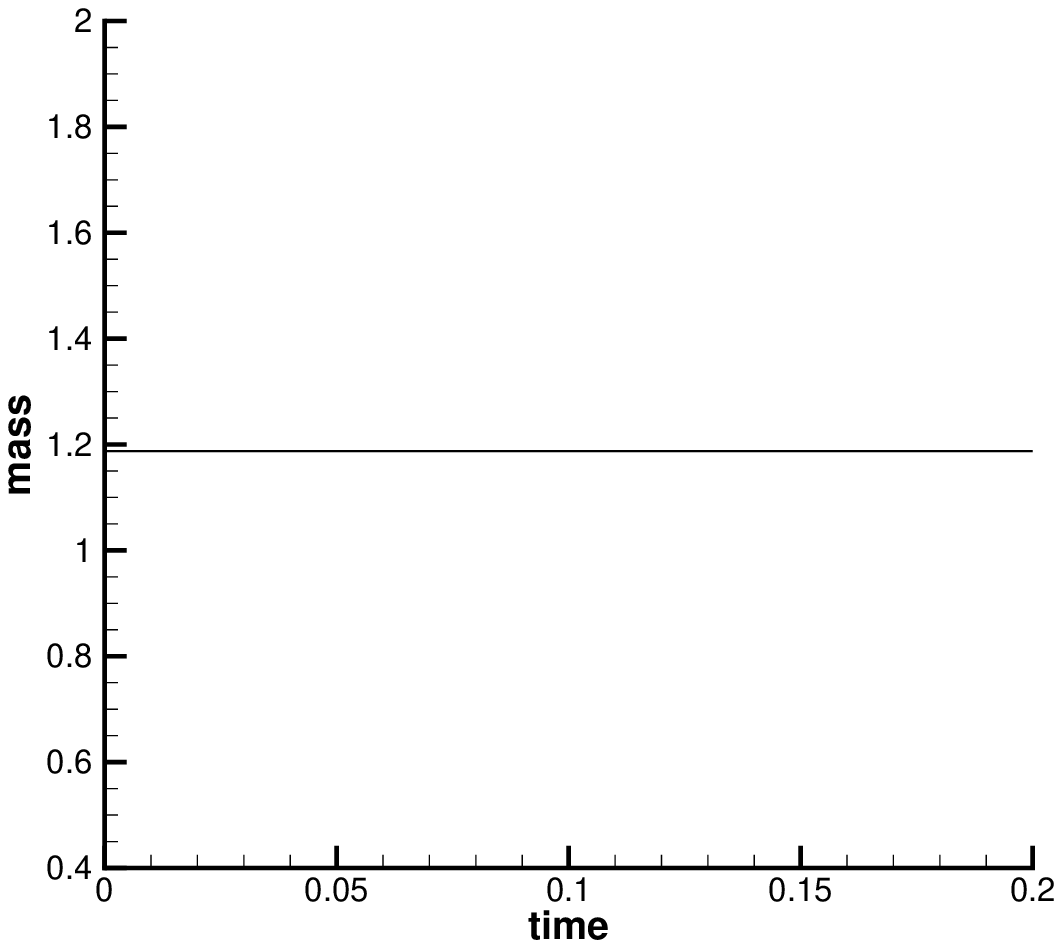}  &  
			\includegraphics[width=0.4\textwidth]{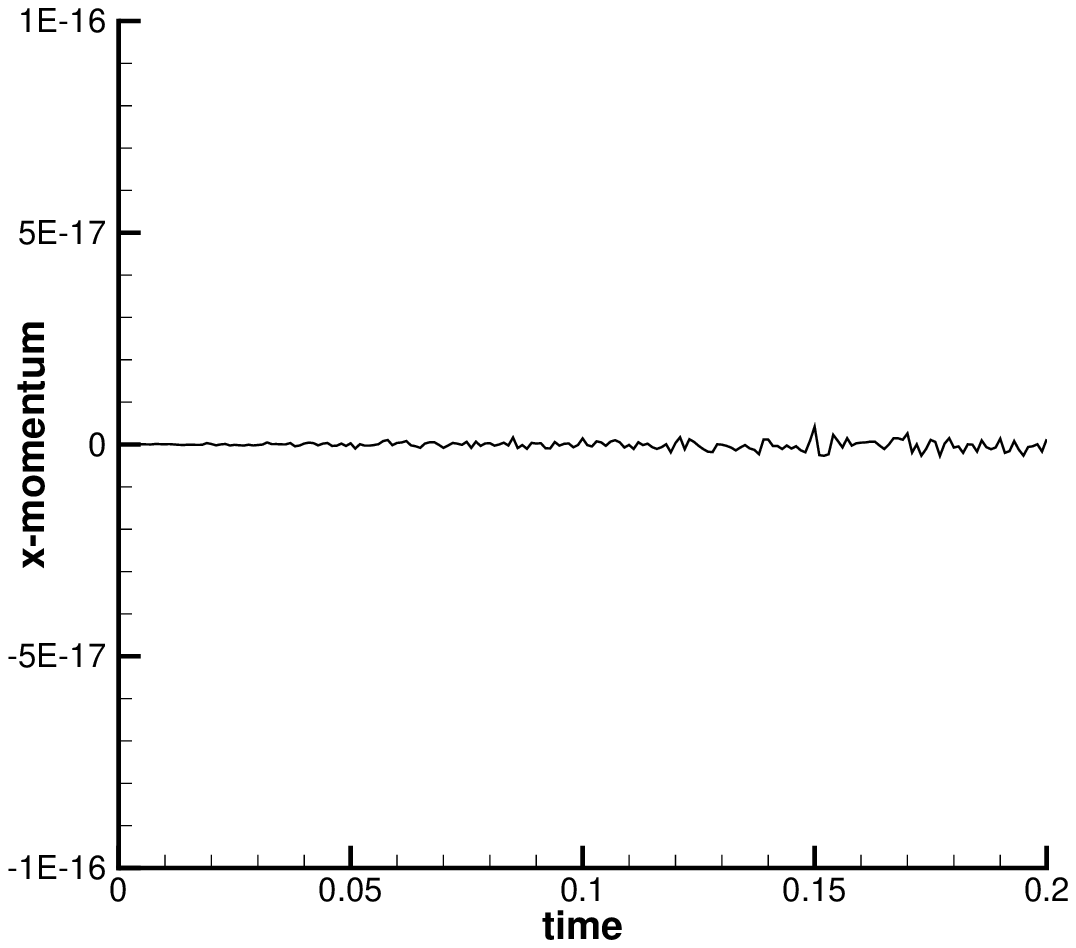}  \\ 
			\includegraphics[width=0.4\textwidth]{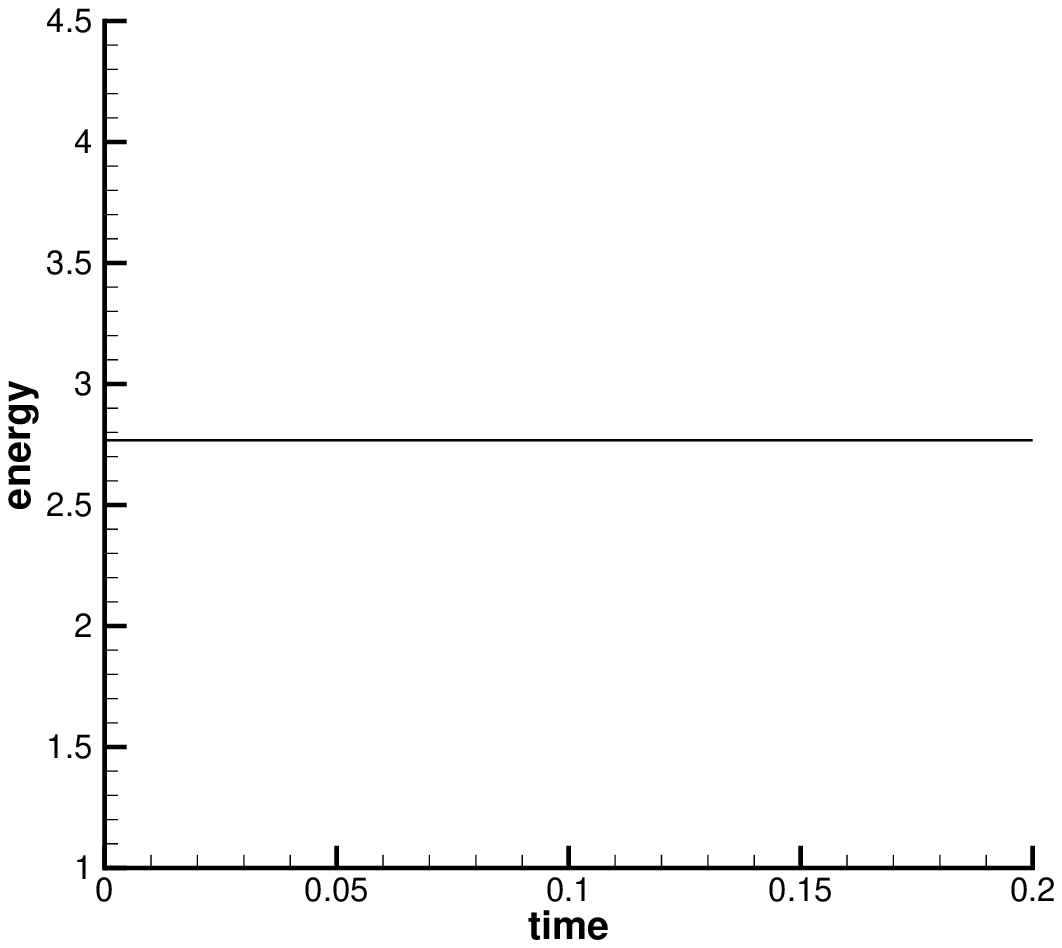}  & 
			\includegraphics[width=0.4\textwidth]{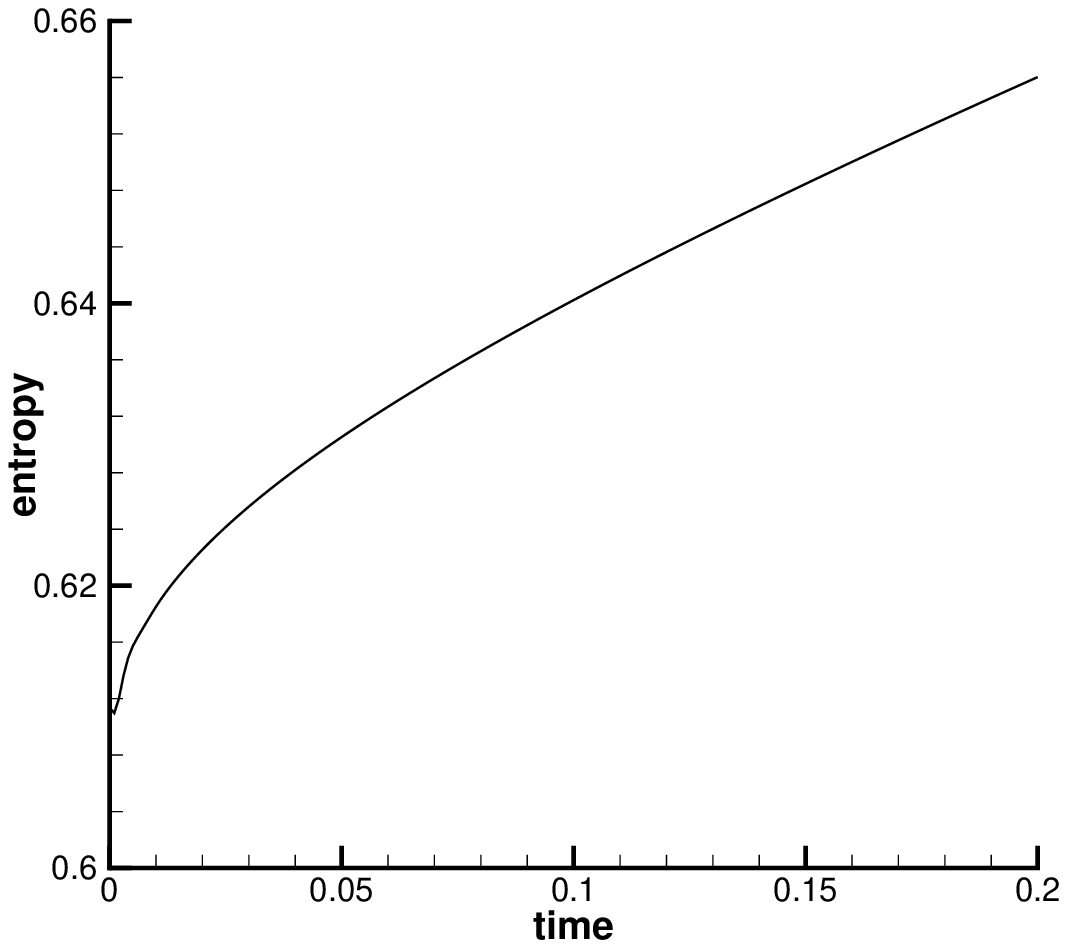}    
		\end{tabular} 
		\caption{\textcolor{black}{2D explosion problem EP1. Time evolution of mass, $x$-momentum, total energy and entropy.}  } 
		\label{fig.ep2dfluid.cons}
	\end{center}
\end{figure}

\paragraph{EP1: Fluid limit of the model} We first solve the governing PDE system in the quasi inviscid fluid limit of the model, i.e. $\tau_1 \ll 1$ and $\tau_2 \ll 1$. For the inner state we set the density and the pressure to
$\rho_{in}=1$ and $p_{in}=1$, while in the outer state we impose $\rho_{out}=0.125$ and $p_{out}=0.1$. The initial velocity is globally set to $\mathbf{v}=0$,
the initial thermal impulse vector is equal to $\J=0$ and the distortion field is initialized as $\A=\mathbf{I}$. The other parameters are chosen as  $\gamma=1.4$, $c_v=2.5$, $c_s=1$, $c_h=1$, $\rho_0=1$, $\tau_1 = 10^{-8}$ and $\tau_2 = 10^{-10}$.
The mesh is composed by $500 \times 500$ control volumes and the final time is set to $t=0.2$. Note that with the chosen parameters, all Mach numbers are of the order of unity. 
The reference solution is obtained by solving a 1D radial Euler system with source terms using a classical second-order MUSCL-Hancock scheme on a very fine mesh, as described in detail in \cite{toro-book,SIGPR}. 
In the left panel of Fig. \ref{fig.ep2dfluidsolid} the 3D density contours are shown at the final time. In Fig. \ref{fig.ep2dfluid} a comparison of the reference solution with the numerical results obtained at the aid of the new semi-implicit 4-split scheme is provided. One observes overall a good agreement. \textcolor{black}{In Fig. \ref{fig.ep2dfluid.cons} we present the time evolution of mass, $x$-momentum, total energy and entropy integrated over the domain. As expected, mass momentum and total energy are perfectly conserved. Our numerical results show that physical entropy is monotonically increasing in time, as it should, although we were not able to prove a discrete entropy inequality for our scheme. }

\begin{figure}[!htbp]
	\begin{center}
		\begin{tabular}{cc} 
			\includegraphics[width=0.45\textwidth]{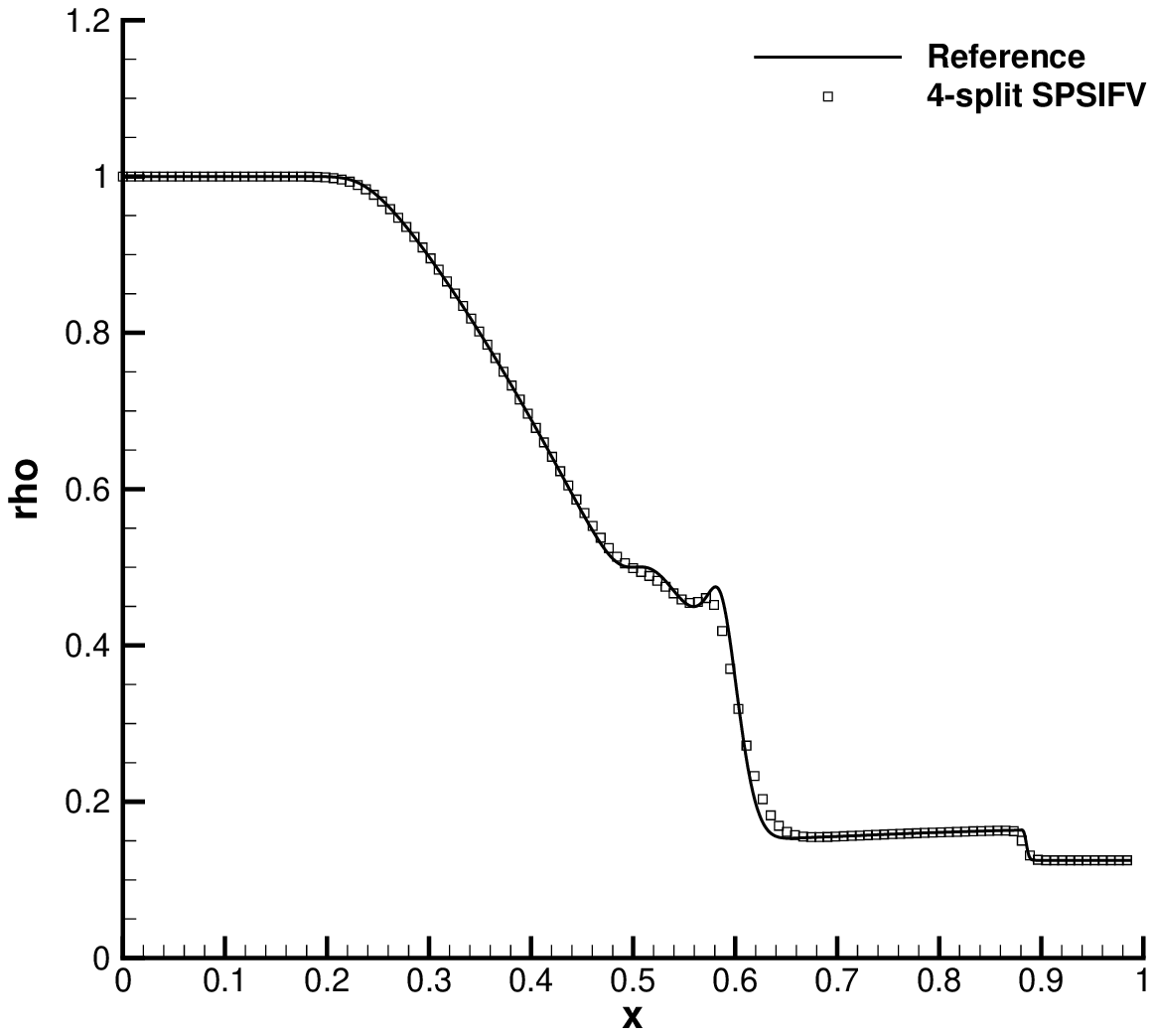}  &
			\includegraphics[width=0.45\textwidth]{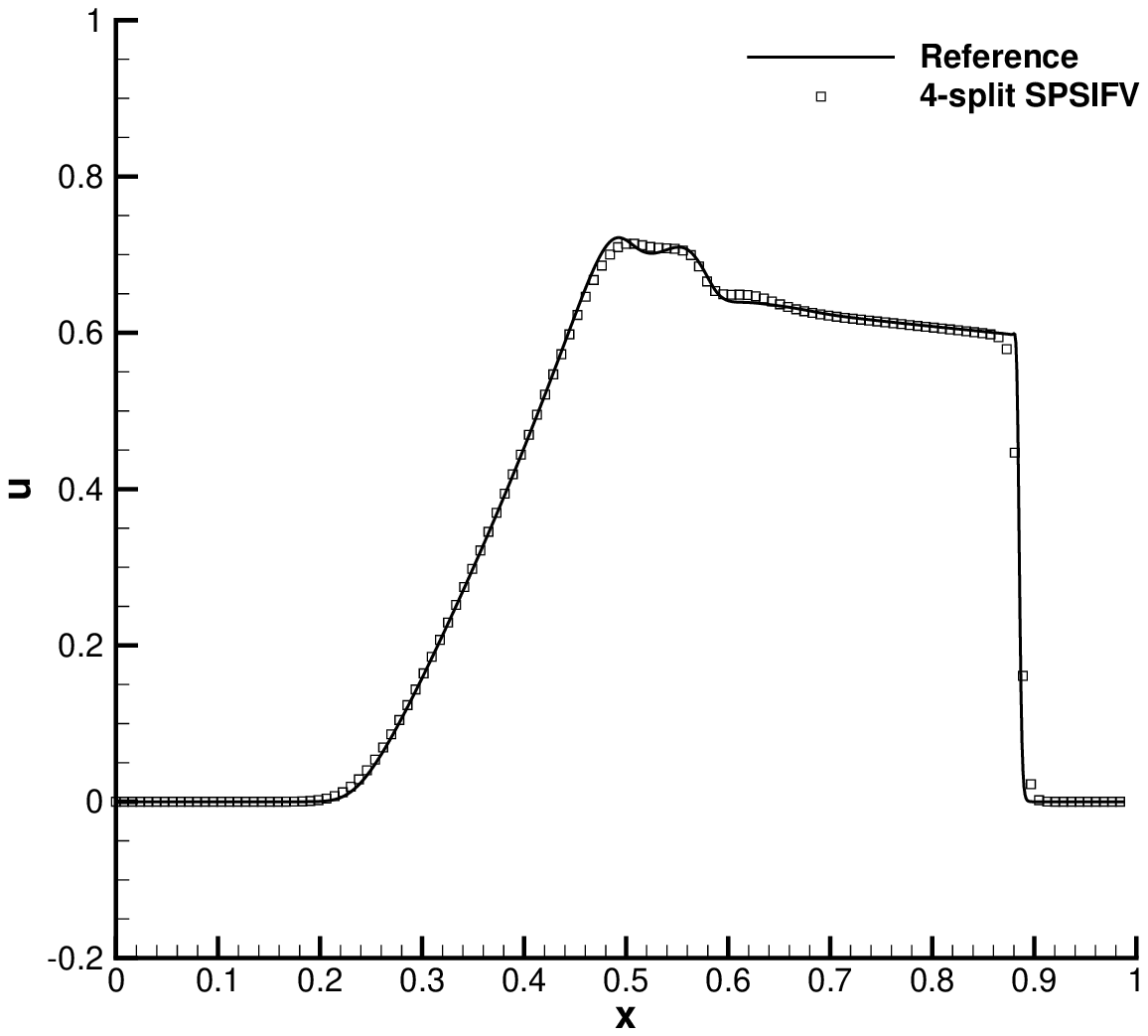}   \\ 
			\includegraphics[width=0.45\textwidth]{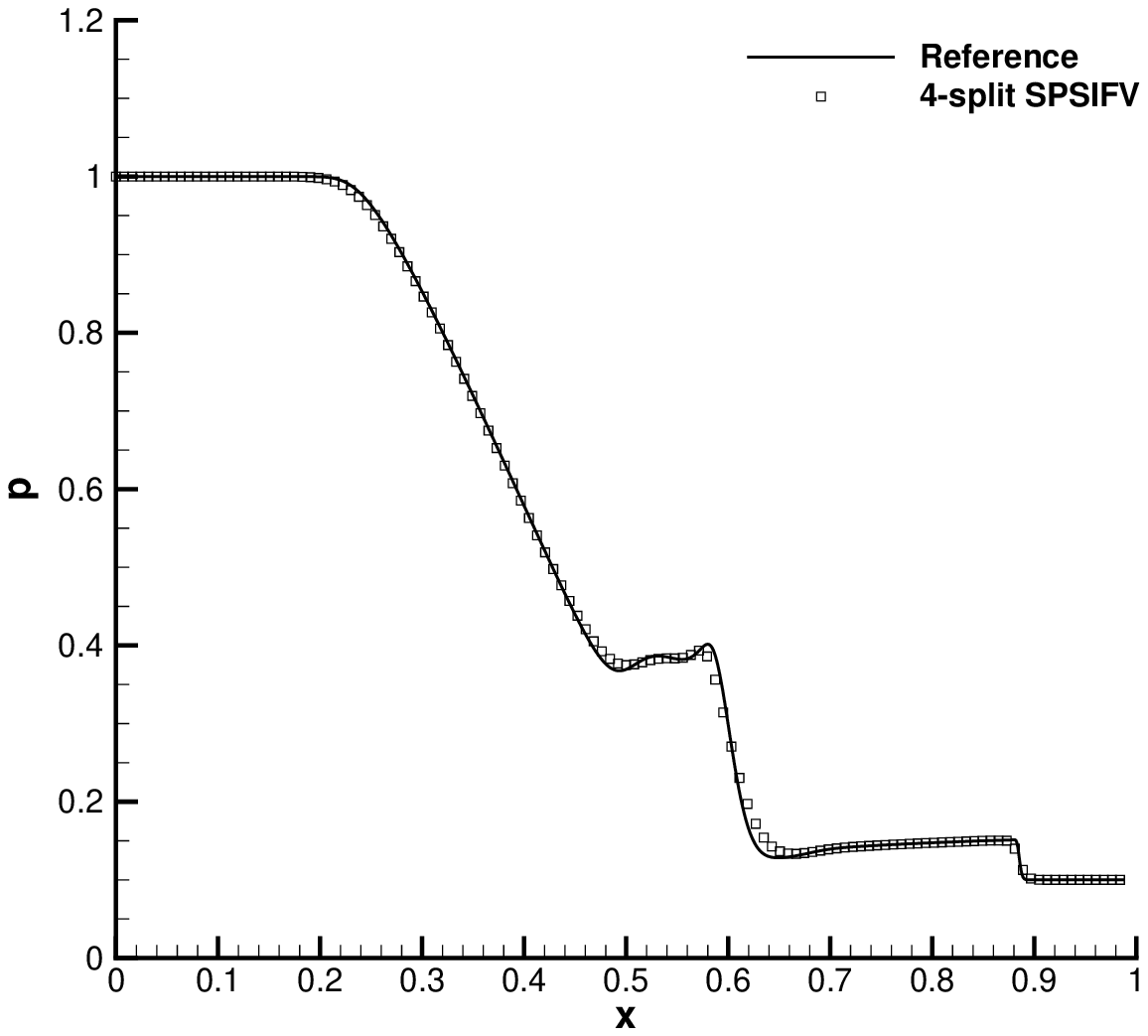}  & 
			\includegraphics[width=0.45\textwidth]{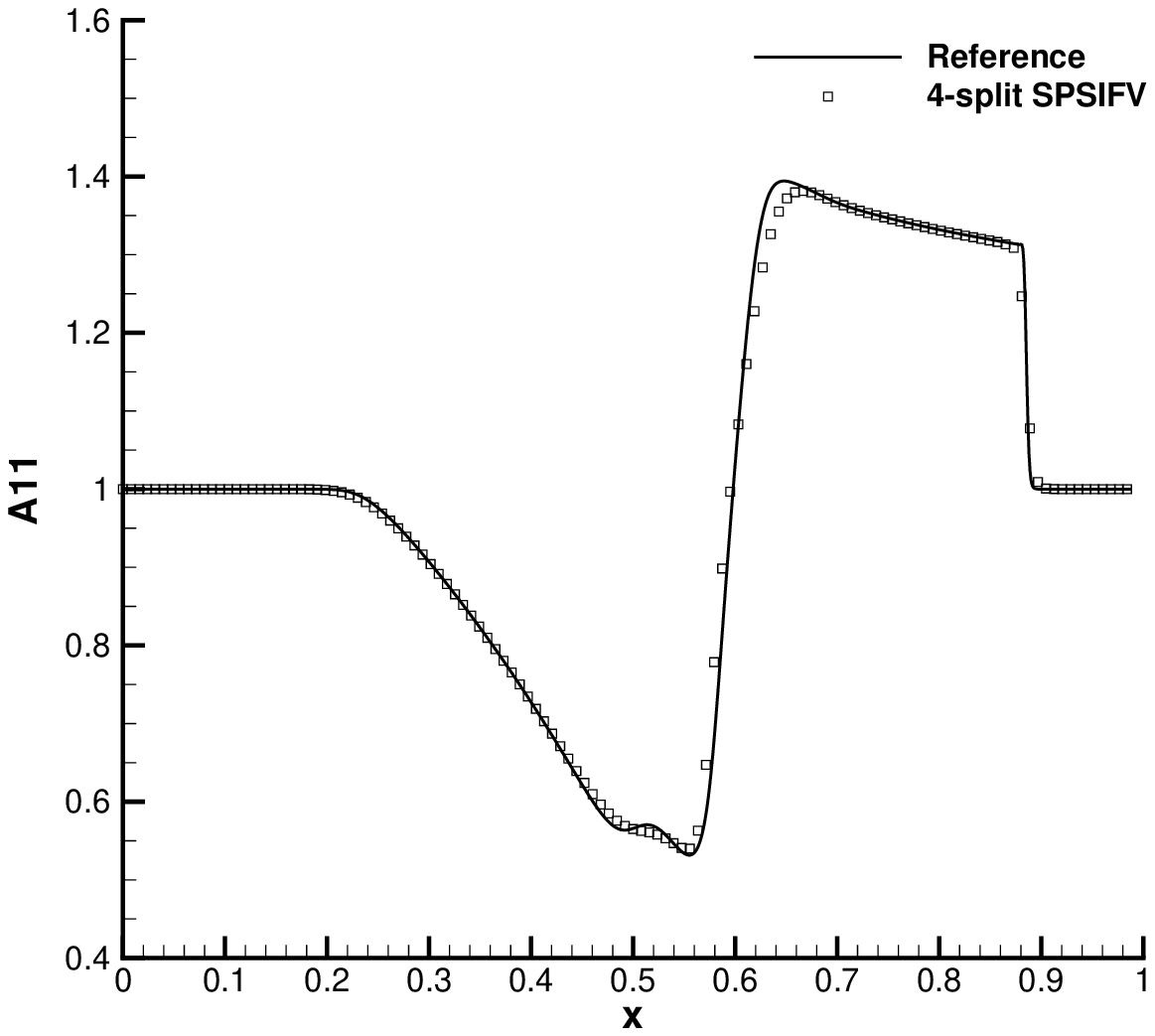}  
		\end{tabular} 
		\caption{Cut along the $x$ axis for the 2D explosion problem EP2 at time $t=0.15$. Comparison of a fine grid reference solution (solid line) against the numerical solution of the homogeneous GPR model (solid limit, $\tau_1=\tau_2=10^{20}$)  obtained with the new 4-split SIFV scheme (square symbols). Density $\rho$ (top left),
velocity component $u$ (top right),  pressure $p$ (bottom left) and distortion field component $A_{11}$ (bottom right).   } 
		\label{fig.ep2dsolid}
	\end{center}
\end{figure}

\begin{figure}[!htbp]
	\begin{center}
		\begin{tabular}{cc} 
			\includegraphics[width=0.4\textwidth]{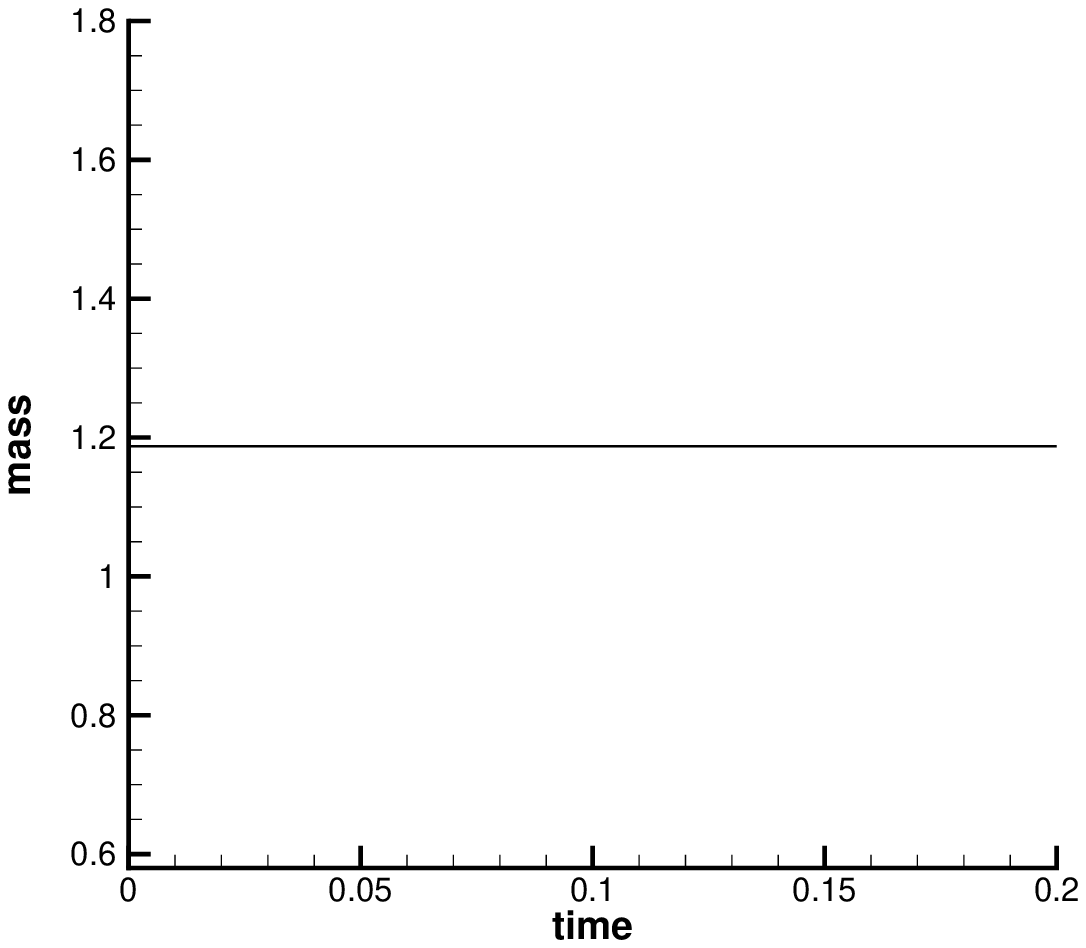}  &  
			\includegraphics[width=0.4\textwidth]{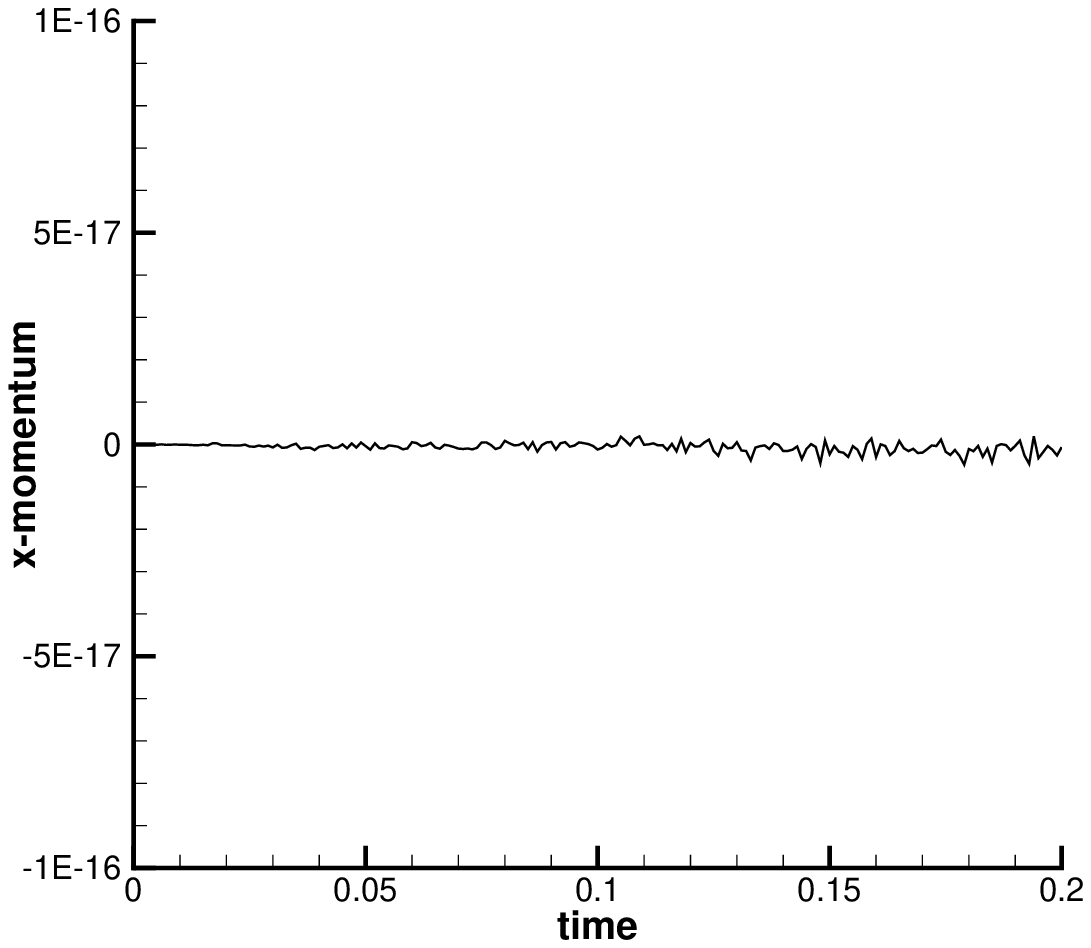}  \\ 
			\includegraphics[width=0.4\textwidth]{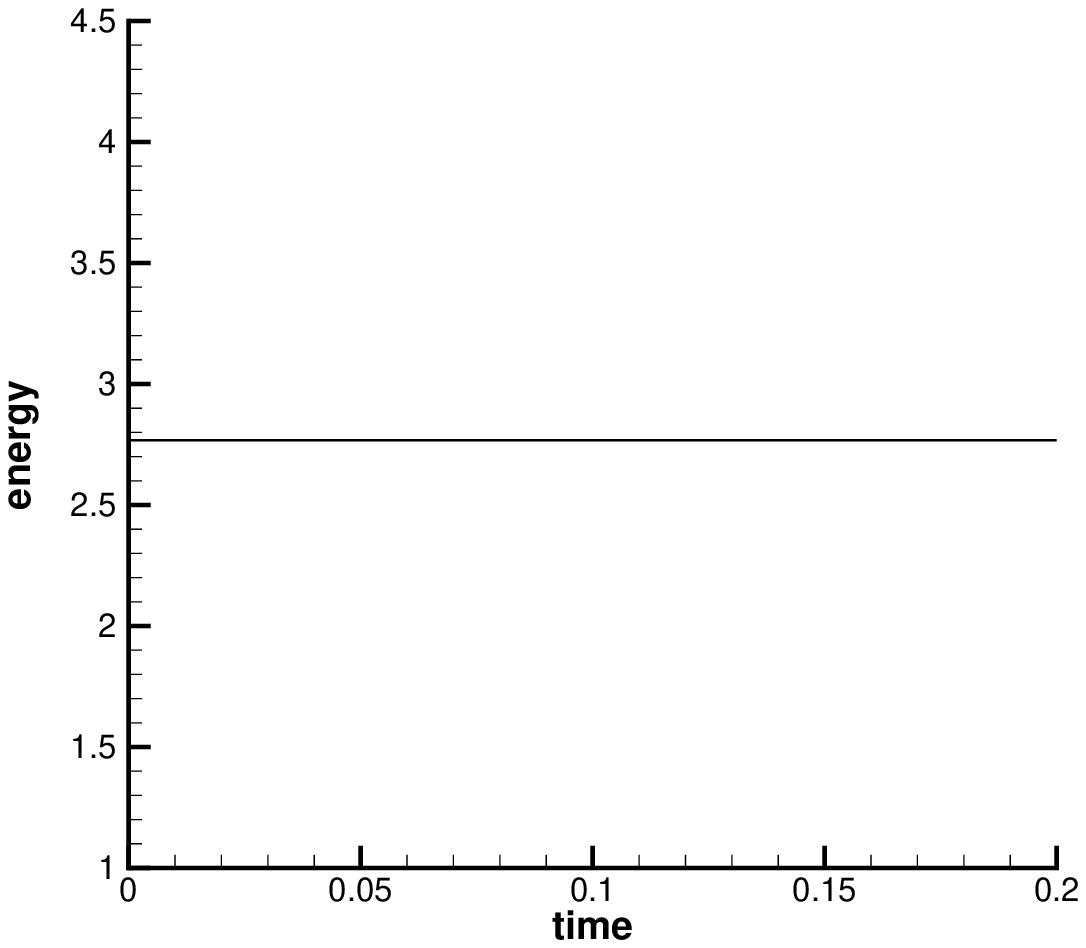}  & 
			\includegraphics[width=0.4\textwidth]{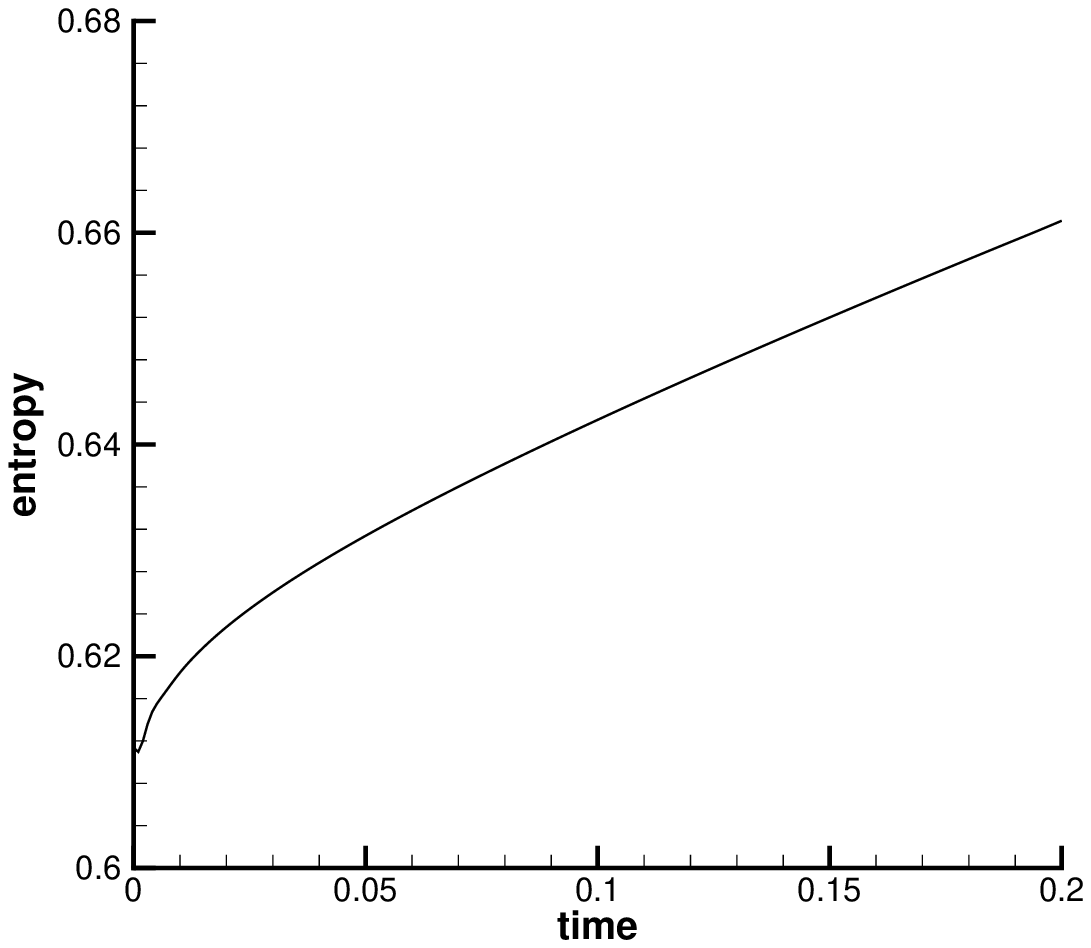}    
		\end{tabular} 
		\caption{\textcolor{black}{2D explosion problem EP2. Time evolution of mass, $x$-momentum, total energy and physical entropy.}  } 
		\label{fig.ep2dsolid.cons}
	\end{center}
\end{figure}

\paragraph{EP2: Solid limit of the model} We now solve the governing PDE system in the solid limit,
i.e. 
$\tau_1 \to \infty$ and $\tau_2 \to \infty$, with the same initial conditions and parameters as described above for EP1, apart from the relaxation times, which are set to $\tau_1 = \tau_2 = 10^{20}$. The mesh is again composed of $500 \times 500$ cells and the final time is now set to $t=0.15$. We solve the problem twice, once with the new 4-split SPSIFV scheme and another time with the thermodynamically compatible HTC scheme presented in \cite{HTCGPR,HTCAbgrall}, in order to obtain a numerical reference solution. A 3D density contour plot at the final time is depicted in the right panel of Fig. \ref{fig.ep2dfluidsolid}. The comparison of the computational results between the new 4-split SIFV scheme and the numerical reference solution is shown in via 1D cuts along the $x$-axis in Fig. \ref{fig.ep2dsolid}. One can again note a good agreement between the two solutions.  
\textcolor{black}{In Fig. \ref{fig.ep2dsolid.cons} we show the time evolution of mass, $x$-momentum, total energy and entropy integrated over the domain. Again, mass momentum and total energy are perfectly conserved, while physical entropy is monotonically increasing in time, as it should be. However, we have to stress again that so far we were not able to prove an entropy inequality for the scheme presented in this paper. }

% % % % % % % % % % % % % % % % % % % % % % % % % % % % % %
%                  Conclusions                            %
% % % % % % % % % % % % % % % % % % % % % % % % % % % % % %

\section{Conclusions}\label{sec.conclusions}
In this paper, we have designed and implemented a new structure-preserving semi-implicit 4-split scheme for the solution of the unified first order model of continuum mechanics originally introduced in \cite{PeshRom2014,GPRmodel}. The fluxes are split into four different contributions, leading to  i) convective, ii) heat, iii) mechanical and iv) pressure sub-systems. The convection terms are discretized explicitly, while the remaining sub-systems are treated implicitly in a sequential manner. In this way, the maximum admissible time step is limited by a CFL-type stability condition that only depends on the bulk velocity of the medium and not on the acoustic, shear and heat wave speeds. This makes the resulting method very efficient for the simulation of low Mach number flows or solid mechanics with high shear speeds. The implicit discretization of the distortion field and thermal impulse equations allows the stiff relaxation limits of the model to be properly retrieved, thus ensuring the asymptotic preserving property. More specifically, the novel scheme approaches the discrete Navier-Stokes stress tensor and the Fourier law when the stiffness parameters in the source terms go to zero. Likewise, the implicit solution of the pressure sub-system allows the scheme to be consistent with the low Mach number limit of the model, given by the incompressible equations. Mimetic finite difference operators on a vertex-staggered Cartesian mesh permit to respect the curl-free involutions embedded in the governing equations of the distortion tensor and the thermal impulse in the case of linear or vanishing source terms, while a robust shock-capturing finite volume discretization is adopted for the nonlinear convective terms. Several benchmarks for both fluid and solid mechanics have been shown, assessing the accuracy and the effectiveness of the novel scheme. The suite of test cases spans a wide range of Mach numbers as well as different regimes of the material, from ideal fluids to elastic solids, demonstrating the wide applicability of the proposed approach.

Further research will be devoted to the adoption of unstructured meshes with the aim of tackling more complex geometries, while preserving additional structure of the system exactly, see e.g. \cite{BoscheriGPRGCL,HTCLagrange,HTCLagrangeGPR}, as well as the incorporation of electromagnetic fields, following the general lines of \cite{GPRmodelMHD}. To improve the accuracy in space, high order extension is also foreseen following the ideas recently forwarded in \cite{CompatibleDG1}, while higher order in time may
be achieved according to the lines presented in \cite{Spiteri3Split}. 
\textcolor{black}{Last but not least, in future work we will try to find a compatible semi-implicit discretization that also provably satisfies a cell entropy inequality.}

% % % % % % % % % % % % % % % % % % % % % % % % % % % % % %
%                Acknowledgment                          %
% % % % % % % % % % % % % % % % % % % % % % % % % % % % % %

\section*{Acknowledgments}

	M.D. was financially supported by the Italian Ministry of University
	and Research (MUR) in the framework of the PRIN 2022 project \textit{High order structure-preserving semi-implicit schemes for hyperbolic equations} and via the  Departments of Excellence  Initiative 2018--2027 attributed to DICAM of the University of Trento (grant L. 232/2016).
	M.D. was also funded by the Fondazione Caritro via the project SOPHOS.
	A.T. acknowledges the financial support of the Agence Nationale de la Recherche (ANR) via the project DELFIN  (project number ANR-25-CE46-7729-01), and the MIUR Departments of Excellence Initiative 2018--2027 via the visiting program of DICAM of the University of Trento (grant L. 232/2016). 
	W.B. received financial support by the Italian Ministry of University and Research (MUR) with the PRIN Project 2022 No.2022N9BM3N and by the "Institut des Mathématiques pour la Planète Terre" (France).

    \vspace{3mm}

	M.D., A.T., M.T. and W.B. are members of the Gruppo Nazionale per il Calcolo Scientifico dell'Istituto Nazionale di Alta Matematica (GNCS-INdAM).

%
%
%
%
%% % % % % % % % % % % % % % % % % % % % % % % % % % % % % %
%% % % % % % % % % % % % % % % % % % % % % % % % % % % % % %
%%              Bibliography
%% % % % % % % % % % % % % % % % % % % % % % % % % % % % % %
%% % % % % % % % % % % % % % % % % % % % % % % % % % % % % %
\bibliographystyle{elsarticle-num}
\bibliography{./biblio}

\begin{thebibliography}{10}
\expandafter\ifx\csname url\endcsname\relax
  \def\url#1{\texttt{#1}}\fi
\expandafter\ifx\csname urlprefix\endcsname\relax\def\urlprefix{URL }\fi
\expandafter\ifx\csname href\endcsname\relax
  \def\href#1#2{#2} \def\path#1{#1}\fi

\bibitem{MunzPark}
J.~Park, C.~Munz, Multiple pressure variables methods for fluid flow at all
  mach numbers, Int. J. Numer. Meth. Fl. 49 (2005) 905--931.

\bibitem{KleinMach}
R.~Klein, Semi-implicit extension of a godunov-type scheme based on low mach
  number asymptotics {I}: one-dimensional flow, J. Comput. Phys. 121 (1995)
  213--237.

\bibitem{Casulli1984}
V.~Casulli, D.~Greenspan, Pressure method for the numerical solution of
  transient, compressible fluid flows, Int. J. Numer. Meth. Fl. 4~(11) (1984)
  1001--1012.

\bibitem{Casulli1990}
V.~Casulli, Semi-implicit finite difference methods for the two-dimensional
  shallow water equations, J. Comput. Phys. 86 (1990) 56--74.

\bibitem{BosFil2016}
S.~Boscarino, F.~Filbet, G.~Russo, {High Order Semi-implicit Schemes for Time
  Dependent Partial Differential Equations}, J. Sci. Comput. 68 (2016)
  975--1001.

\bibitem{BP2021}
W.~Boscheri, L.~Pareschi, High order pressure-based semi-implicit imex schemes
  for the 3d navier-stokes equations at all mach numbers, J. Comput. Phys. 434
  (2021) 110206.

\bibitem{AscRuuSpi}
U.~M. Ascher, S.~J. Ruuth, R.~J. Spiteri, {Implicit-explicit {R}unge-{K}utta
  methods for time-dependent partial differential equations}, Appl. Numer.
  Math. 25 (1982) 151--167.

\bibitem{BP2017}
S.~Boscarino, L.~Pareschi, {On the asymptotic properties of IMEX Runge–Kutta
  schemes for hyperbolic balance laws}, J. Comput. Appl. Math. 316 (2017)
  60--73.

\bibitem{BosRus}
S.~Boscarino, G.~Russo, {On a class of uniformly accurate {IMEX}
  {R}unge-{K}utta schemes and applications to hyperbolic systems with
  relaxation}, SIAM J. Sci. Comput. 31 (2009) 1926--1945.

\bibitem{PR_IMEX}
L.~Pareschi, G.~Russo, {Implicit-explicit Runge-Kutta schemes and applications
  to hyperbolic systems with relaxation}, J. Sci. Comput. 25 (2005) 129--155.

\bibitem{BPR2017}
S.~Boscarino, L.~Pareschi, G.~Russo, {A unified {IMEX} {R}unge-{K}utta approach
  for hyperbolic systems with multiscale relaxation}, SIAM J. Numer. Anal.
  55~(4) (2017) 2085--2109.

\bibitem{Hofer}
S.~Osher, F.~Solomon, {A partially implicit method for large stiff systems of
  Ode's with only few equations introducing small time-constants}, SIAM J.
  Numer. Anal. 13 (1976) 645--663.

\bibitem{JINAP1999}
S.~Jin, {Efficient asymptotic-preserving ({AP}) schemes for some multiscale
  kinetic equations}, SIAM J. Sci. Comput. 21~(2) (1999) 441--454.

\bibitem{JP2001}
S.~Jin, L.~Pareschi, {Asymptotic-preserving ({AP}) schemes for multiscale
  kinetic equations: a unified approach}, in: G.~W. H.~Freist\"uhler (Ed.),
  Hyperbolic problems: theory, numerics, applications, Vol. 141 of Internat.
  Ser. Numer. Math., Birkh\"{a}user, Basel, 2001, pp. 573--582.

\bibitem{KLARAP1999}
A.~Klar, {An asymptotic preserving numerical scheme for kinetic equations in
  the low {M}ach number limit}, SIAM J. Numer. Anal. 36~(5) (1999) 1507--1527.

\bibitem{Degond2}
F.~Cordier, P.~Degond, A.~Kumbaro, {An Asymptotic-Preserving all-speed scheme
  for the Euler and Navier-Stokes equations}, J. Comp. Phys. 231 (2012)
  5685--5704.

\bibitem{DumbserCasulli2016}
M.~Dumbser, V.~Casulli, {A conservative, weakly nonlinear semi-implicit finite
  volume method for the compressible Navier-Stokes equations with general
  equation of state}, Appl. Math. Comput. 272 (2016) 479--497.

\bibitem{TavelliDumbser2017}
M.~Tavelli, M.~Dumbser, {A pressure-based semi-implicit space-time
  discontinuous Galerkin method on staggered unstructured meshes for the
  solution of the compressible Navier-Stokes equations at all Mach numbers}, J.
  Comput. Phys. 341 (2017) 341--376.

\bibitem{Avgerinos2019}
S.~Avgerinos, F.~Bernard, A.~Iollo, G.~Russo, {Linearly implicit all Mach
  number shock capturing schemes for the Euler equations}, J. Comput. Phys. 393
  (2019) 278--312.

\bibitem{Thomann2020}
A.~Thomann, M.~Zenk, G.~Puppo, C.~Klingenberg, An all speed second order {IMEX}
  relaxation scheme for the {E}uler equations, Commun. Comput. Phys. 28~(2)
  (2020) 591--620.

\bibitem{dimarco2012}
G.~Dimarco, L.~Pareschi, High order asymptotic-preserving schemes for the
  boltzmann equation, Comptes Rendus. Math{\'e}matique 350~(9-10) (2012)
  481--486.

\bibitem{BD2021_FVBoltz}
W.~Boscheri, G.~Dimarco, High order finite volume schemes with imex time
  stepping for the boltzmann model on unstructured meshes, Comput. Method.
  Appl. M. 387 (2021) 114180.

\bibitem{xiong2015_BGK_DG}
T.~Xiong, J.~Jang, F.~Li, J.-M. Qiu, High order asymptotic preserving nodal
  discontinuous galerkin imex schemes for the bgk equation, J. Comput. Phys.
  284 (2015) 70--94.

\bibitem{bisi2022}
M.~Bisi, W.~Boscheri, G.~Dimarco, M.~Groppi, G.~Martal{\`o}, A new mixed
  boltzmann-bgk model for mixtures solved with an imex finite volume scheme on
  unstructured meshes, Appl. Math. Comput. 433 (2022) 127416.

\bibitem{BosTav22}
W.~Boscheri, M.~Tavelli, High order semi-implicit schemes for viscous
  compressible flows in 3d, Appl. Math. Comput. 434 (2022) 127457.

\bibitem{ToroVazquez}
E.~Toro, M.~V\'azquez-Cend\'on, {Flux splitting schemes for the Euler
  equations}, Computers and Fluids 70 (2012) 1--12.

\bibitem{Fambri_3splitMHD}
F.~Fambri, A novel structure preserving semi-implicit finite volume method for
  viscous and resistive magnetohydrodynamics, Int. J. Numer. Meth. Fl. 93~(12)
  (2021) 3447--3489.

\bibitem{dematte2024}
R.~Dematt{\'e}, A.~Farmakalides, S.~Millmore, N.~Nikiforakis, An all mach
  number scheme for visco-resistive magnetically-dominated mhd flows, J.
  Comput. Phys. 514 (2024) 113229.

\bibitem{BosTho_3splitMHD}
W.~Boscheri, A.~Thomann, A structure-preserving semi-implicit imex finite
  volume scheme for ideal magnetohydrodynamics at all mach and alfv{\'e}n
  numbers, J. Sci. Comput. 100~(3) (2024) 67.

\bibitem{PeshRom2014}
I.~Peshkov, E.~Romenski, A hyperbolic model for viscous {{N}ewtonian} flows,
  Continuum Mech. Thermodyn. 28 (2016) 85--104.

\bibitem{GPRmodel}
M.~Dumbser, I.~Peshkov, E.~Romenski, O.~Zanotti, {High order ADER schemes for a
  unified first order hyperbolic formulation of continuum mechanics: Viscous
  heat--conducting fluids and elastic solids}, J. Comput. Phys. 314 (2016)
  824--862.

\bibitem{SIGPR}
W.~Boscheri, M.~Dumbser, M.~Ioriatti, I.~Peshkov, E.~Romenski, {A
  structure-preserving staggered semi-implicit finite volume scheme for
  continuum mechanics}, J. Comput. Phys. 424 (2021) 109866.

\bibitem{ChiocchettiSI}
S.~Chiocchetti, M.~Dumbser, {An exactly curl-free staggered semi-implicit
  finite volume scheme for a first order hyperbolic model of viscous two-phase
  flows with surface tension}, J. Sci. Comput. 94 (2023) 24.

\bibitem{GodRom2003}
S.~Godunov, E.~Romenski, Elements of continuum mechanics and conservation laws,
  Kluwer Academic/Plenum Publishers, 2003.

\bibitem{MaxwellGLM}
M.~Dumbser, A.~Lucca, I.~Peshkov, O.~Zanotti, {Variational derivation and
  compatible discretizations of the Maxwell-GLM system}, Proceedings of the
  Royal Society A 481 (2025) 20240864.

\bibitem{CompatibleDG1}
R.~Abgrall, M.~Dumbser, P.-H. Maire, A simple and general framework for the
  construction of exactly div-curl-grad compatible discontinuous galerkin
  finite element schemes on unstructured simplex meshes, J. Comput. Phys. 541
  (2025) 114340.

\bibitem{Yee66}
K.~Yee, {Numerical solution of initial voundary value problems involving
  Maxwell equation in isotropic media}, IEEE Trans. Antenna Propagation 14
  (1966) 302--307.

\bibitem{DeVore}
C.~DeVore, {Flux-corrected transport techniques for multidimensional
  compressible magnetohydrodynamics}, J. Comput. Phys. 92 (1991) 142--160.

\bibitem{BalsaraSpicer1999}
D.~Balsara, D.~Spicer, A staggered mesh algorithm using high order godunov
  fluxes to ensure solenoidal magnetic fields in magnetohydrodynamic
  simulations, J. Comput. Phys. 149 (1999) 270--292.

\bibitem{Balsara2004}
D.~Balsara, Second-order accurate schemes for magnetohydrodynamics with
  divergence-free reconstruction, Astrophys. J. Suppl. Ser. 151 (2004)
  149--184.

\bibitem{GardinerStone}
T.~Gardiner, J.~Stone, {An unsplit Godunov method for ideal MHD via constrained
  transport}, J. Comput. Phys. 205 (2005) 509--539.

\bibitem{balsarahlle2d}
D.~Balsara, {Multidimensional HLLE Riemann solver: Application to Euler and
  magnetohydrodynamic flows}, J. Comput. Phys. 229 (2010) 1970--1993.

\bibitem{ADERdivB}
D.~Balsara, M.~Dumbser, {Divergence-free MHD on unstructured meshes using high
  order finite volume schemes based on multidimensional {R}iemann solvers}, J.
  Comput. Phys. 299 (2015) 687--715.

\bibitem{HymanShashkov1997}
J.~Hyman, M.~Shashkov, {Natural discretizations for the divergence, gradient,
  and curl on logically rectangular grids}, Comput. Math. Appl. 33 (1997)
  81--104.

\bibitem{JeltschTorrilhon2006}
R.~Jeltsch, M.~Torrilhon, {On curl--preserving finite volume discretizations
  for shallow water equations}, BIT Numerical Mathematics 46 (2006) S35--S53.

\bibitem{Torrilhon2004}
M.~Torrilhon, M.~Fey, {Constraint-preserving upwind methods for
  multidimensional advection equations}, SIAM J. Numer. Anal. 42 (2004)
  1694--1728.

\bibitem{Margolin2000}
G.~Margolin, M.~Shashkov, P.~Smolarkiewicz, {A discrete operator calculus for
  finite difference approximations}, Comput. Method. Appl. M. 187 (2000)
  365--383.

\bibitem{Lipnikov2014}
K.~Lipnikov, G.~Manzini, M.~Shashkov, {Mimetic finite difference method}, J.
  Comput. Phys. 257 (2014) 1163--1227.

\bibitem{Carney2013}
T.~Carney, N.~Morgan, S.~Sambasivan, M.~Shashkov, {A cell--centered Lagrangian
  Godunov--like method for solid dynamics}, Computers and Fluids 83 (2013)
  33--47.

\bibitem{Nedelec1}
J.~N\'ed\'elec, {Mixed finite elements in R3}, Numer. Math. 35 (1980) 315--341.

\bibitem{Nedelec2}
J.~N\'ed\'elec, {A new family of mixed finite elements in R3}, Numer. Math. 50
  (1986) 57--81.

\bibitem{Cantarella}
J.~Cantarella, D.~DeTurck, H.~Gluck, {Vector calculus and the topology of
  domains in 3-space}, Am. Math. Mon. 109 (2002) 409--442.

\bibitem{Hiptmair}
R.~Hiptmair, {Finite elements in computational electromagnetism}, Acta Numer.
  11 (2002) 237--339.

\bibitem{Monk}
P.~Monk, {Finite Element Methods for Maxwell's Equations}, Oxford University
  Press, Oxford, 2003.

\bibitem{Arnold}
D.~Arnold, R.~Falk, R.~Winther, {Finite element exterior calculus, homological
  techniques, and applications}, Acta Numer. 15 (2006) 1--155.

\bibitem{Alonso2015}
A.~A. Rodriguez, A.~Valli, Finite element potentials, Appl. Numer. Math. 95
  (2015) 2--14.

\bibitem{CAMPOSPINTO2016}
M.~{Campos Pinto}, M.~Mounier, E.~Sonnendrücker, Handling the divergence
  constraints in maxwell and vlasov–maxwell simulations, Appl. Math. Comput.
  272 (2016) 403--419, recent Advances in Numerical Methods for Hyperbolic
  Partial Differential Equations.

\bibitem{Zampa1}
E.~Zampa, S.~Busto, M.~Dumbser, {A divergence-free hybrid finite volume /
  finite element scheme for the incompressible MHD equations based on
  compatible finite element spaces}, Applied Numerical Mathematics 198 (2024)
  346--374.

\bibitem{Zampa2}
E.~Zampa, M.~Dumbser, {An asymptotic-preserving and exactly mass-conservative
  semi-implicit scheme for weakly compressible flows based on compatible finite
  elements}, Journal of Computational Physics 521 (2025) 113551.

\bibitem{SPDG2023}
W.~Boscheri, G.~Dimarco, L.~Pareschi, Locally structure-preserving div-curl
  operators for high order discontinuous galerkin schemes, J. Comput. Phys. 486
  (2023) 112130.

\bibitem{Ern2023}
A.~Ern, J.-L. Guermond, The discontinuous galerkin approximation of the
  grad-div and curl-curl operators in first-order form is involution-preserving
  and spectrally correct, SIAM J. Numer. Anal. 61~(6) (2023) 2940--2966.

\bibitem{perrier2024}
V.~Perrier, Development of discontinuous galerkin methods for hyperbolic
  systems that preserve a curl or a divergence constraint, arXiv preprint
  arXiv:2405.04347 (2024).

\bibitem{bonelle2015}
J.~Bonelle, A.~Ern, Analysis of compatible discrete operator schemes for the
  stokes equations on polyhedral meshes, IMA Journal of numerical analysis
  35~(4) (2015) 1672--1697.

\bibitem{DiPietro2023}
D.~A. Di~Pietro, J.~Droniou, An arbitrary-order discrete de rham complex on
  polyhedral meshes: Exactness, poincar{\'e} inequalities, and consistency,
  Lond. Math. S. 23~(1) (2023) 85--164.

\bibitem{Barsukow2024}
W.~Barsukow, R.~Loub{\`e}re, P.-H. Maire, A node-conservative vorticity
  preserving finite volume method for linear acoustics on unstructured grids,
  Math. Comput. (2024).

\bibitem{Sidilkover2025}
D.~Sidilkover, {Spurious vorticity in Eulerian and Lagrangian methods}, J.
  Comput. Phys. 520 (2025) 113510.

\bibitem{Rom1998}
E.~Romenski, Hyperbolic systems of thermodynamically compatible conservation
  laws in continuum mechanics, Math. Comput. Modell. 28(10) (1998) 115--130.

\bibitem{HTCGPR}
S.~Busto, M.~Dumbser, I.~Peshkov, E.~Romenski, {On thermodynamically compatible
  finite volume schemes for continuum mechanics}, SIAM J. Sci. Comput. 44
  (2022) A1723--A1751.

\bibitem{Maire2007}
P.-H. Maire, R.~Abgrall, J.~Breil, J.~Ovadia, A cell-centered {Lagrangian}
  scheme for two-dimensional compressible flow problems, SIAM J. Sci. Comput.
  29 (2007) 1781--1824.

\bibitem{Despres2009}
G.~Carr\'e, S.~D. Pino, B.~Despr\'es, E.~Labourasse, {A cell-centered
  Lagrangian hydrodynamics scheme on general unstructured meshes in arbitrary
  dimension.}, J. Comput. Phys. 228 (2009) 5160--5183.

\bibitem{Maire2011}
P.-H. Maire, A high-order one-step sub-cell force-based discretization for
  cell-centered lagrangian hydrodynamics on polygonal grids, Computers and
  Fluids 46(1) (2011) 341--347.

\bibitem{Maire2020}
P.-H. Maire, I.~Bertron, R.~Chauvin, B.~Rebourcet, {Thermodynamic consistency
  of cell-centered Lagrangian schemes}, Computers and Fluids 203 (2020) 104527.

\bibitem{HTCLagrange}
W.~Boscheri, M.~Dumbser, P.-H. Maire, {A new thermodynamically compatible
  finite volume scheme for Lagrangian gas dynamics}, SIAM J. Sci. Comput. 44
  (2024) A1723--A1751.

\bibitem{HTCLagrangeGPR}
W.~Boscheri, M.~Dumbser, R.~Loub\`ere, P.-H. Maire, {A structure-preserving and
  thermodynamically compatible cell-centered Lagrangian finite volume scheme
  for continuum mechanics}, SIAM J. Num. Anal.To appear.

\bibitem{BoscheriGPRGCL}
W.~Boscheri, R.~Loub\`ere, J.~Braeunig, P.-H. Maire, A geometrically and
  thermodynamically compatible finite volume scheme for continuum mechanics on
  unstructured polygonal meshes, J. Comput. Phys. 507 (2024) 112957.

\bibitem{KlaMaj}
S.~Klainermann, A.~Majda, Singular limits of quasilinear hyperbolic systems
  with large parameters and the incompressible limit of compressible fluid,
  Comm. Pure Appl. Math. 34 (1981) 481--524.

\bibitem{KlaMaj82}
S.~Klainermann, A.~Majda, Compressible and incompressible fluids, Comm. Pure
  Appl. Math. 35 (1982) 629--651.

\bibitem{Klein2001}
R.~Klein, N.~Botta, T.~Schneider, C.~Munz, S.Roller, A.~Meister, L.~Hoffmann,
  T.~Sonar, Asymptotic adaptive methods for multi--scale problems in fluid
  mechanics, J. of Eng. Math. 39 (2001).

\bibitem{MunzDumbserRoller}
C.~Munz, M.~Dumbser, S.~Roller, Linearized acoustic perturbation equations for
  low {Mach} number flow with variable density and temperature, J. Comput.
  Phys. 224 (2007) 352--364.

\bibitem{toro-book}
E.~Toro, {Riemann} Solvers and Numerical Methods for Fluid Dynamics, Springer,
  2009.

\bibitem{Ghia1982}
U.~Ghia, K.~N. Ghia, C.~T. Shin, {High-Re solutions for incompressible flow
  using Navier-Stokes equations and multigrid method}, J. Comput. Phys. 48
  (1982) 387--411.

\bibitem{HTCAbgrall}
R.~Abgrall, S.~Busto, M.~Dumbser, {A simple and general framework for the
  construction of thermodynamically compatible schemes for computational fluid
  and solid mechanics}, Appl. Math. Comput. 440 (2023) 127629.

\bibitem{GPRmodelMHD}
M.~Dumbser, I.~Peshkov, E.~Romenski, O.~Zanotti, {H}igh order {ADER} schemes
  for a unified first order hyperbolic formulation of {N}ewtonian continuum
  mechanics coupled with electro--dynamics, J. Comput. Phys. 348 (2017)
  298--342.

\bibitem{Spiteri3Split}
R.~Spiteri, A.~Tavassoli, S.~Wei, A.~Smolyakov, {Beyond Strang: a Practical
  Assessment of Some Second-Order 3-Splitting Methods}, Commun. Appl. Math.
  Comput. 7 (2025) 95–114.

\end{thebibliography}

% % % % % % % % % % % % % % % % % % % % % % % % % % % % % %
% % % % % % % % % % % % % % % % % % % % % % % % % % % % % %
%                   Appendix                              %
% % % % % % % % % % % % % % % % % % % % % % % % % % % % % %
% % % % % % % % % % % % % % % % % % % % % % % % % % % % % %
%\appendix

\end{document}